\documentclass[titlepage,11pt]{article}

\usepackage{latexsym}
\usepackage{amsfonts}
\usepackage{amsmath}
\usepackage{tikz}
\usepackage{float}
\usepackage{lmodern}

\newcommand\blackslug{\hbox{\hskip 1pt \vrule width 4pt height 8pt depth 1.5pt
        \hskip 1pt}}
\newcommand\bbox{\hfill \quad \blackslug \medbreak}
\def\dd{\hbox{-}}
\def\cc{\hbox{-}\cdots\hbox{-}}
\def\ll{,\ldots,}
\def\cupcup{\cup\cdots\cup}

\title{Even-hole-free graphs still have bisimplicial vertices}
\author{Maria Chudnovsky\thanks{This material is based upon work supported in part by the U. S. Army
Research Office under grant   number W911NF-16-1-0404, and by  NSF grant DMS-1763817.}
and
Paul Seymour\thanks{Partially supported by NSF grant  DMS-1800053 and AFOSR grant A9550-19-1-0187.}\\
Princeton University, Princeton, NJ 08544}

\date{June 5, 2019; revised \today}

\newtheorem{thm}{}[section]

\newcommand{\Proof}{\noindent{\bf Proof.}\ \ }

\begin{document}
\maketitle

\begin{abstract}
A {\em hole} in a graph is an induced subgraph which is a cycle of length at 
least four. A hole is called {\em even} if it has an even number of vertices. 
An  {\em even-hole-free}  graph is a graph with no even holes. 
A vertex of a graph is {\em bisimplicial} 
if the set of its neighbours is the union of two cliques. 

In an earlier paper~\cite{bisimplicial}, Addario-Berry, Havet and Reed, with the authors, 
claimed to prove a conjecture of Reed, that
every even-hole-free graph has a bisimplicial  vertex, but we have recently been shown that the ``proof'' has a serious error.
Here we give a proof using a different approach.

\end{abstract}

\section{Introduction}

All graphs in this paper are finite and simple.  Let $G$ be a graph.   
A {\em clique} in $G$ is a set of pairwise adjacent vertices. 
A vertex is {\em bisimplicial (in $G$)} if its 
neighbourhood is the union of two cliques. 
A {\em hole} in a graph is an induced subgraph that is a cycle of length at 
least four. 
A hole is  {\em even} if it has even length
and {\em odd} otherwise.
A graph is {\em even-hole-free} if it contains no even hole. 
The following was conjectured in~\cite{RAR}:
\begin{thm}
\label{conj}
Every non-null even-hole-free graph has a bisimplicial vertex.
\end{thm}
Louigi Addario-Berry, 
Fr\'ed\'eric Havet and Bruce Reed, with the authors, published a ``proof'' in~\cite{bisimplicial}. However,
there is a major error in this proof, pointed out to us recently by Rong Wu. 
The flawed proof is for a result (theorem 3.1 of that paper) that is fundamental
to much of the remainder of the paper, and we have not been able to 
fix the error (although we still believe 3.1 of that paper to be true). 
The error in~\cite{bisimplicial} is in the very last line of the proof of theorem 3.1 of that paper: 
we say ``it follows that $N_G(v)=N_{G'}(v)$, and so $v$ is
bisimplicial in $G$''; and this is not correct, since cliques of $G'$ may not be cliques of $G$.

In this paper we give a different proof of \ref{conj}. For inductive purposes
we prove something a little stronger, namely:

\begin{thm}\label{mainthm}
Let $G$ be even-hole-free, and let $K$ be a clique of $G$ 
with $|K|\le 2$. Let $M$ be the set of vertices in $V(G)\setminus V(K)$ with no neighbour in  $V(K)$. If $M\ne \emptyset$, some vertex in $M$
is bisimplicial.
\end{thm}

The proof is via two decomposition theorems for even-hole-free graphs. Most of the paper is concerned with 
proving these decomposition theorems,
and at the end we give the application to finding a bisimplicial vertex.

\section{Preliminaries, and a sketch of the proof}

Before we can outline the proof we need more definitions.
Let
$S$ be a subset of $V(G)$. We denote by $G[S]$ the subgraph of $G$ induced on
$S$, and by $G\setminus S$ the subgraph of $G$ induced on  $V(G) \setminus S$.
We say $S\subseteq V(G)$ is {\em connected} if $G[S]$ is connected.
The {\em neighbourhood} of  $S$, denoted by  $N_G(S)$ (or $N(S)$ when there is
no risk of confusion), is the set of all vertices of
$V(G) \setminus S$  with a neighbour in  $S$, and $N[S]$ means $N(S)\cup S$. If $S=\{v\}$, we write $N_G(v)$
instead of  $N_G(\{v\})$; for an induced subgraph $H$ of $G$, we define $N(H)$ to be $N(V(H))$, and so on.
A subgraph $S$ is {\em dominating} in $G$ if $N[S]=V(G)$, and {\em non-dominating} otherwise.

Two disjoint subsets $A,B$ 
of $V(G)$ are {\em complete} to each other if every vertex of $A$ is adjacent
to every vertex of $B$, and {\em anticomplete} to each other if no vertex of
$A$ is adjacent to any vertex of $B$. If $A=\{a\}$, we write ``$a$ is
complete (anticomplete) to $B$'' instead of ``$\{a\}$   is complete
(anticomplete) to $B$''.

The {\em length} of a path is
the number of edges in it. A path is called {\em even} if its length is even,
and {\em odd} otherwise. Let the vertices of $P$ be $p_1,\ldots, p_k$ in
order. Then $p_1,p_k$ are called the {\em ends} of $P$ (sometimes we say
$P$ is {\em from $p_1$ to $p_k$} or {\em between $p_1$ and $p_k$}), and the
set $V(P) \setminus \{p_1, p_k\}$ is the {\em interior} of $P$ and is
denoted by $P^*$.
For $1 \leq i < j  \leq k$ we will write $p_i \dd P \dd p_j$ or
$p_j \dd P \dd p_i$ to mean the subpath of $P$ between $p_i$ and $p_j$.
More generally, if $S$ is an induced subgraph of a graph $G$, and $u,v$ both have neighbours in $V(S)$,
we denote by $u\dd S\dd v$ some induced path between $u,v$ with interior in $V(S)$. (Here $u,v$ might or might not belong to $V(S)$.)
If $H$ is a cycle, and $a,b$ and $c$ are three vertices of $H$ such
that $a$ is adjacent to $b$, then $a \dd b \dd H \dd c$ is a path, consisting
of $a$, and the subpath of $H\setminus \{a\}$ between $b$ and $c$.
A {\em triangle} is a set of three vertices, pairwise adjacent, and we use the same word for the subgraph induced on a triangle.

Here are some types of graph that we will need:
\begin{itemize}
\item A {\em theta} with {\em ends} $s,t$ is a graph that is the union of three paths $R_1,R_2,R_3$, each with the same pair of ends $s,t$,
each of length more than one, and pairwise vertex-disjoint except for their ends.
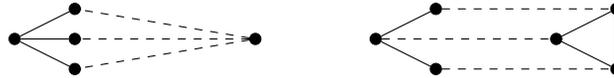
\begin{figure}[H]
\centering

\begin{tikzpicture}[scale=.8,auto=left]
\tikzstyle{every node}=[inner sep=1.5pt, fill=black,circle,draw]
\node (z) at (0,0) {};
\node (a) at (1,.5) {};
\node (b) at (1,0) {};
\node (c) at (1,-.5) {};
\node (d) at (4,0) {};

\foreach \from/\to in {z/a,z/b, z/c}
\draw [-] (\from) -- (\to);
\foreach \from/\to in {d/a,d/b, d/c}
\draw [dashed] (\from) -- (\to);

\node (a) at (6,0) {};
\node (c1) at (7,.5) {};
\node (c3) at (7,-.5) {};
\node (b1) at (10,.5) {};
\node (b2) at (9,0) {};
\node (b3) at (10,-.5) {};

\foreach \from/\to in {a/c1, a/c3,b1/b2,b1/b3,b2/b3}
\draw [-] (\from) -- (\to);
\foreach \from/\to in {c1/b1, a/b2,c3/b3}
\draw [dashed] (\from) -- (\to);

\end{tikzpicture}

\caption{A theta and a pyramid (dashed lines mean paths of arbitrary positive length)} \label{theta}
\end{figure}

\item A {\em pyramid} with {\em apex} $a$ and {\em base} $\{b_1,b_2,b_3\}$ is a graph $P$ in which
\begin{itemize}
\item $a,b_1,b_2,b_3$ are distinct vertices, and $\{b_1,b_2,b_3\}$ is a triangle,
\item $P$ is the union of this triangle and three paths $R_1,R_2, R_3$, where $R_i$ has ends $a,b_i$ for $i = 1,2,3$, and
\item $R_1,R_2,R_3$ are pairwise vertex-disjoint except for their common end, and at least two of $R_1,R_2,R_3$
have length at least two.
\end{itemize}
\item A {\em near-prism} with {\em bases} $\{a_1,a_2,a_3\}, \{b_1,b_2,b_3\}$ is a graph $P$ in which
\begin{itemize}
\item $\{a_1,a_2,a_3\}$ and $\{b_1,b_2,b_3\}$ are triangles, and $\{a_1,a_2,a_3\}\cap  \{b_1,b_2,b_3\}=\{a_3\}\cap \{b_3\}$ 
(that is, the triangles are disjoint except that possibly $a_3=b_3$).
\item $P$ is the union of these two triangles and three paths
$R_1,R_2,R_3$, where $R_i$ has ends $a_i,b_i$
for $i = 1,2,3$ (and so $R_3$ has length zero if $a_3=b_3$).
\item $R_1,R_2,R_3$ are pairwise vertex-disjoint.
\end{itemize}
If $a_3\ne b_3$, $P$ is also called a {\em prism}.
\begin{figure}[H]
\centering

\begin{tikzpicture}[scale=.8,auto=left]
\tikzstyle{every node}=[inner sep=1.5pt, fill=black,circle,draw]
\node (a1) at (1,0) {};
\node (a2) at (0,.5) {};
\node (a3) at (0,-.5) {};
\node (b1) at (3,0) {};
\node (b2) at (4,.5) {};
\node (b3) at (4,-.5) {};

\foreach \from/\to in {a1/a2,a1/a3,a2/a3,b1/b2,b1/b3,b2/b3}
\draw [-] (\from) -- (\to);
\foreach \from/\to in {a1/b1,a2/b2,a3/b3}
\draw [dashed] (\from) -- (\to);

\node (a1) at (7,.5) {};
\node (a2) at (9,0) {};
\node (a3) at (10,-1) {};
\node (b1) at (13,.5) {};
\node (b2) at (11,0) {};

\foreach \from/\to in {a1/a2,a1/a3,a2/a3,b1/b2,b1/a3,b2/a3}
\draw [-] (\from) -- (\to);
\foreach \from/\to in {a1/b1,a2/b2}
\draw [dashed] (\from) -- (\to);

\end{tikzpicture}

\caption{Near-prisms} \label{nearprism}
\end{figure}
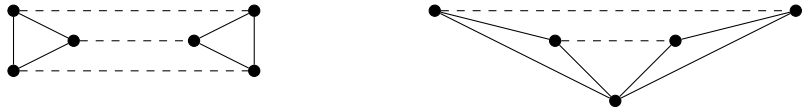

\item A {\em wheel} is a graph consisting of a hole $H$ and
a vertex $v \not \in V(H)$ with at least three neighbours
in $V(H)$, and if it has exactly three neighbours in $V(H)$ then no two of them are adjacent.
We call $v$ its {\em centre} and $H$ its {\em hole}. If $v$ has $k$ neighbours in $H$ we also call it a {\em $k$-wheel.}
If $k$ is even we call it an {\em even wheel}.
\end{itemize}
For a theta, pyramid or near-prism, we call $R_1,R_2,R_3$ its {\em constituent paths}.
It is easy to see that:
\begin{thm}
\label{subgraphs}
No even-hole-free graph contains a theta, a near-prism or an even wheel as an induced subgraph.
\end{thm}
Even-hole-free graphs can contain pyramids, however.
A pyramid is {\em short} if one of the three constituent paths has length one.

An {\em extended near-prism} is a graph obtained from a near-prism by adding one extra edge, as follows. 
Let $R_1,R_2,R_3$ be as in the definition
of a near-prism, and let $a\in R_1^*$ and $b\in R_2^*$; and add an edge $ab$. (It is important that $a,b$ do not belong 
to the triangles.)
We call $ab$ the {\em cross-edge} of the extended near-prism.

\begin{figure}[H]
\centering

\begin{tikzpicture}[scale=.8,auto=left]
\tikzstyle{every node}=[inner sep=1.5pt, fill=black,circle,draw]
\node (a2) at (1,0) {};
\node (a1) at (0,.5) {};
\node (a3) at (0,-.5) {};
\node (b2) at (3,0) {};
\node (b1) at (4,.5) {};
\node (b3) at (4,-.5) {};
\node (d1) at (2,.5) {};
\node (d2) at (2,0) {};

\foreach \from/\to in {a1/a2,a1/a3,a2/a3,b1/b2,b1/b3,b2/b3, d1/d2}
\draw [-] (\from) -- (\to);
\foreach \from/\to in {a1/d1,d1/b1,a2/d2, d2/b2,a3/b3}
\draw [dashed] (\from) -- (\to);

\node (a1) at (7,.5) {};
\node (a2) at (9,0) {};
\node (a3) at (10,-1) {};
\node (b1) at (13,.5) {};
\node (b2) at (11,0) {};
\node (d1) at (10,.5) {};
\node (d2) at (10,0) {};

\foreach \from/\to in {a1/a2,a1/a3,a2/a3,b1/b2,b1/a3,b2/a3, d1/d2}
\draw [-] (\from) -- (\to);
\foreach \from/\to in {a1/b1,a2/b2}
\draw [dashed] (\from) -- (\to);

\end{tikzpicture}

\caption{Extended near-prisms} \label{extnearprism}
\end{figure}
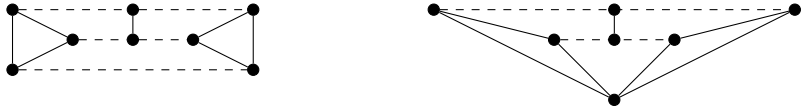

A vertex $a\in V(G)$ is {\em splendid} if
\begin{itemize}
\item $V(G)\setminus N[a]$ is connected;
\item every vertex in $N(a)$ has a neighbour in $V(G)\setminus N[a]$; and
\item there is no short pyramid with apex $a$ in $G$.
\end{itemize}

Now we can sketch the idea of the proof.
In order to prove \ref{mainthm}, we use induction on $|V(G)|$. From a result of \cite{bisimplicial} (that did not depend on theorem
3.1 of that paper, and so is still valid), we may assume that $G$ admits no ``full star cutset'' (defined later). It follows
that, with $K$ as in \ref{mainthm}, there is a splendid vertex $a\in V(G)\setminus N[K]$. We can assume that $a$ is not bisimplicial.
Now there are two possibilities:
\begin{itemize}
\item there is an extended near-prism in which $a$ belongs to the cross-edge;
\item there is a pyramid with apex $a$, but there is no
extended near-prism in which $a$ belongs to the cross-edge.
\end{itemize}
In both cases we use a decomposition theorems to find a smaller subgraph to which we 
can apply the inductive hypthesis and win. There are two different decompositions theorems. 
The first gives a decomposition of $G$ relative to an extended near-prism, and is fully general
(that is, it does not require any vertex to be splendid), and so may be useful in other applications. The second is more tailored
to our application, in that it needs $a$ to be splendid.

To apply these to find bisimplicial vertices, we use that
both theorems provide a choice of subgraphs (two in the first case, three in the second) 
that are separated from the remainder of the graph in a convenient way, and we can prove inductively that
all these subgraphs contain bisimplicial vertices of $G$; and in both cases these subgraphs are sufficiently widely separated that at 
least one of these bisimplicial
vertices has no neighbours in $K$.

The main part of the paper concerns proving the two decomposition theorems, and we use them to prove \ref{mainthm} in the final section.

\section {Some results from \cite{bisimplicial}.}
We will need to use some results of \cite{bisimplicial} that did not depend on the flawed theorem 3.1 of that paper.
A {\em cutset} in $G$ is a subset
$C$ of $V(G)$ such that $V(G) \setminus C$ is the union of two non-empty
sets, anticomplete to each  other. A {\em star cutset} is a cutset 
consisting of a vertex and some of its neighbours. If $v$ together with a 
subset of  $N(v)$ is a cutset, we say that $v$ is a {\em centre} of 
this star cutset. A star  cutset $C$ is called {\em full} if
it consists of a vertex and all its neighbours. We need the following, theorem 4.2 of~\cite{bisimplicial}:

\begin{thm} 
\label{nofullstar}
Let $G$ be an even-hole-free graph such that, for every even-hole-free graph $H$
 with fewer vertices than $G$, and every non-dominating clique $J$ of $H$ with $|J|\le 2$, there
is a bisimplicial vertex of $H$ in $V(H)\setminus N_H(J)$.
Assume that there exists a 
non-dominating clique $K$ with $|K|\le 2$  such  that no vertex 
of $V(G) \setminus N_G(K)$ is bisimplicial in $G$. Then $G$ does not admit a 
full star cutset.
\end{thm}
(Actually, theorem 4.2 in~\cite{bisimplicial} takes a stronger hypothesis than we give here,
requiring that the dubious theorem 3.1 of that paper holds for all graphs with fewer vertices than $G$; but fortunately its proof
in that paper does not use the extra hypothesis, so we can legitimately omit it.)
We will also need the following consequence of theorem 4.5 of \cite{bisimplicial}:

\begin{thm}
\label{bigvertex}
Let $G$ be even-hole-free, let $H$ be a hole in $G$, and let $a\notin V(H)$. If $G$ admits no full star cutset with centre $a$, 
then either 
\begin{itemize}
\item $a$ is complete or anticomplete to $V(H)$; or
\item $H[V(H)\cap N(a)]$ is a path; or
\item $a$ has exactly three neighbours in $H$, and two of them are adjacent.
\end{itemize}
\end{thm}

\section{Tree strip systems}

In this section and the next, 
we state and prove the decomposition theorem for even-hole-free graphs that contain an extended near-prism.

Here is an example of an even-hole-free graph, due to Conforti, Cornu\'ejols and Vu\v{s}kovi\'{c}~\cite{conforti}, and see
also~\cite{kristina}.
Start with a tree $T$ with $|V(T)|\ge 3$. 
(A {\em leaf} of $T$ means a vertex of degree exactly one, and a {\em leaf-edge} is an edge incident with a leaf.) 
Let $(A',B')$ be a bipartition of $T$. Since $|V(T)|\ge 3$, each leaf-edge
is incident with only one leaf; let $A$ be the set of leaf-edges incident with a leaf in $A'$, and define $B$
similarly.
Let $L(T)$ be the line-graph of $T$. Thus the vertex set of $L(T)$ is the edge set of $T$, and 
$A, B$ are disjoint subsets of $V(L(T))$.
Add to $L(T)$ two more vertices $a,b$ and the edge $ab$, and make $a$ complete to $A$ and $b$ complete to $B$,
forming a graph $H(T)$ say. Thus $H(T)$ has vertex set $E(T)\cup \{a,b\}$. This graph $H(T)$ is
even-hole-free, but it is helpful for our purposes to impose additional conditions on $T$. We will assume that
$T$ has at least three leaves, and 
for every $v\in V(T)$, there is at most one component $C$ of $V(T)\setminus v$ such that $A'\cap V(C)=\emptyset$,
and at most one such that $B'\cap V(C)=\emptyset$. (Note that every component $C$ of $V(T)\setminus v$ contains a leaf
of $T$ and therefore meets at least one of $A',B'$.) If this additional condition is satisfied, we say that
$H(T)$ is an {\em extended tree line-graph}, and $ab$ is its {\em cross-edge}. 
\begin{figure}[H]
\centering

\begin{tikzpicture}[scale=.8,auto=left]
\tikzstyle{every node}=[inner sep=1.5pt, fill=black,circle,draw]
\node (a) at (0,0) {};
\node (b) at (0,1) {};
\node (c) at (0,2) {};
\node (d) at (0,-1) {};
\node (e) at (2/3,-2) {};
\node (f) at (4/3,-3) {};
\node (g) at (1,0) {};
\node (h) at (-2/3,-2) {}; 
\node (i) at (-1,0) {};
\node (j) at (-2,0) {};
\node (k) at (-8/3,-1) {};
\node (l) at (-8/3,1) {};
\node (m) at (-10/3,2) {}; 
\node (n) at (2/3,3) {};
\node (o) at (-2/3,3) {};
\node (p) at (-4/3,4) {};

\foreach \from/\to in {a/b,a/g,a/d,a/i,b/c,c/n,c/o,d/e,d/h,e/f,g/i,i/j,j/k,j/l,l/m,o/p}
\draw [-] (\from) -- (\to);

\begin{scope}[shift ={(6,0)}]
\node (a) at (1/2,0) {};
\node (b) at (0,1/2) {};
\node (c) at (-1/2,0) {};
\node (d) at (0,-1/2) {};
\node (e) at (1/3,-3/2) {};
\node (f) at (-1/3,-3/2) {};
\node (g) at (1,-5/2) {};
\node (h) at (0,3/2) {};
\node (i) at (1/3,5/2) {};
\node (j) at (-1/3,5/2) {};
\node (k) at (-1,7/2) {};
\node (l) at (-3/2,0) {};
\node (m) at (-7/3,-1/2) {};
\node (n) at (-7/3,1/2) {};
\node (o) at (-3,3/2) {};
\foreach \from/\to in {a/b,a/c,a/d,b/c,b/d,b/h,c/d,c/l,d/e,d/f,e/f,e/g,h/i,h/j,i/j,j/k,l/m,l/n,m/n,n/o}
\draw [-] (\from) -- (\to);
\end{scope}
\begin{scope}[shift ={(12,0)}]
\node (a) at (1/2,0) {};
\node (b) at (0,1/2) {};
\node (c) at (-1/2,0) {};
\node (d) at (0,-1/2) {};
\node (e) at (1/3,-3/2) {};
\node (f) at (-1/3,-3/2) {};
\node (g) at (1,-5/2) {};
\node (h) at (0,3/2) {};
\node (i) at (1/3,5/2) {};
\node (j) at (-1/3,5/2) {};
\node (k) at (-1,7/2) {};
\node (l) at (-3/2,0) {};
\node (m) at (-7/3,-1/2) {};
\node (n) at (-7/3,1/2) {};
\node (o) at (-3,3/2) {};
\tikzstyle{every node}=[inner sep=1.5pt, fill=white,circle,draw]
\node (B) at (2,4) {};
\node (A) at (-2,4) {};
\foreach \from/\to in {a/b,a/c,a/d,b/c,b/d,b/h,c/d,c/l,d/e,d/f,e/f,e/g,h/i,h/j,i/j,j/k,l/m,l/n,m/n,n/o}
\draw [-] (\from) -- (\to);
\foreach \from/\to in {A/k, A/o, A/f,B/i, B/a, B/g, A/B}
\draw [thick, dotted] (\from) -- (\to);
\draw[thick, dotted] (B) to [bend right=30] (m);
\tikzstyle{every node}=[]
\draw (A) node [left]           {$a$};
\draw (B) node [right]           {$b$};
\end{scope}

\end{tikzpicture}

\caption{A tree $T$, $L(T)$, and $H(T)$. (The dotted lines are just edges.)} \label{extendedtree}
\end{figure}
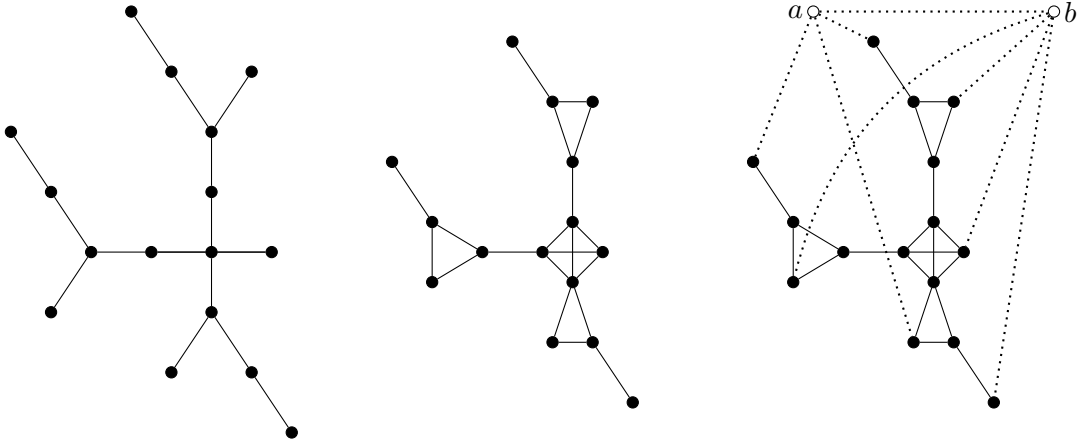

Every extended near-prism is an extended tree line-graph, where the corresponding tree has four leaves and exactly two
vertices of degree three. In the next few sections
we will be working with even-hole-free graphs $G$ that contain extended near-prisms, and therefore the graph also contains
an extended tree line-graph that is maximal (subject to keeping the cross-edge fixed); and examining how the remainder
of the graph attaches to this subgraph will lead us to the decomposition.

Sometimes we have different graphs with the same vertex set or edge set,
and we say {\em $G$-incident} to mean incident in $G$, and {\em $G$-adjacent} to mean adjacent in $G$, and so on.
A {\em branch-vertex} of a tree means a vertex of degree different from two (thus, leaves are branch-vertices).
A {\em branch} of a tree $T$ means a path $P$ of $T$ with distinct ends $u,v$, both branch-vertices, 
such that every vertex of $P^*$ has degree
two in $T$. 
Every edge of $T$ belongs to a unique branch.
A {\em leaf-branch} is a branch such that one of its ends is a leaf of $T$. In general, a {\em leaf-path}
of $T$ means a path of $T$ with one end a leaf of $T$ and the other end a vertex of $T$ that is not a leaf.

Let $T$ be a tree, and let $U$ be the set of branch-vertices of $T$;
and make a new tree $J$ with vertex set $U$  by making $u,v\in U$ $J$-adjacent if there
is a branch of $T$ with ends $u,v$. We call $J$ the {\em shape} of $T$. Thus $J$ has no vertices of degree two;
and $T$ is obtained from $J$ by replacing each edge by a path of positive length.

Let $A,B,C$ be subsets of $V(G)$, with $A, B\ne \emptyset$ and disjoint from $C$, and let $S=(A,B,C)$. A {\rm rung} of $S$, or
an {\em $S$-rung},
is an induced path $p_1\cc p_k$ of $G[A\cup B\cup C]$ such that $p_1\in A$, $p_k\in B$ and $p_2\ll p_{k-1}\in C$, and if
$k>0$ then $p_1\notin B$ and $p_k\notin A$. (If $A\cap B\ne \emptyset$ then perhaps $k=0$.)
If every vertex in $A\cup B\cup C$ belongs to an $S$-rung we call $S$ a {\em strip}.
We denote $A\cup B\cup C$ by $V(S)$.
In the later part of the paper, concerned with ``pyramid strip systems'', we will only need strips $(A,B,C)$
with $A\cap B=\emptyset$, but earlier when we look at ``tree strip systems'' we need to allow $A,B$ to intersect.
A strip $(A,B,C)$ is {\em proper} if $A\cap B=\emptyset$.

Let $J$ be a tree with at least three vertices.
A {\em $J$-strip system} $M$ in a graph $G$ means:
\begin{itemize}
\item for each edge $e=uv$ of $J$, a subset $M_{uv}=M_{vu}=M_e$ of $V(G)$
\item for each $v\in V(J)$, a subset $M_v$ of $V(G)$
\end{itemize}
satisfying the following conditions:
\begin{itemize}
\item the sets $M_{e}\; (e \in E(J))$ are pairwise disjoint;
\item for each $u \in V(J)$, $M_u \subseteq \bigcup (M_{uv}: v \in V(J)$ adjacent to $u)$;
\item for each $uv \in E(J)$, $(M_{uv}\cap M_u, M_{uv}\cap M_v, M_{uv}\setminus (M_u\cup M_v))$ is a 
strip (not necessarily proper);
\item if $uv,wx \in E(J)$ with $u,v,w,x$ all distinct, then there are
no edges between $M_{uv}$ and $M_{wx}$;
\item if $uv,uw \in E(J)$ with $v \ne w$, then $M_u \cap M_{uv}$ is complete to
$M_u \cap M_{uw}$, and there are no other edges between $M_{uv}$ and $M_{uw}$.
\end{itemize}
A rung of the strip $(M_{uv}\cap M_u, M_{uv}\cap M_v, M_{uv}\setminus (M_u\cup M_v))$ will be called an {\em $e$-rung}
or {\em $uv$-rung}. (We leave the dependence on $M$ and $J$ to be understood, for the sake of brevity.)
Let $V(M)$ denote the union of the sets $M_e\;(e\in E(J))$.

Let $J$ be a tree, let $M$ be a $J$-strip system in $G$, and let $(\alpha,\beta)$ be a partition of the set of 
leaves of $J$. We say an edge $ab$ of $G$ is a {\em cross-edge}
for $M$ with {\em partition $(\alpha,\beta)$} if:
\begin{itemize}
\item $J$ has no vertex of degree two, and at least three vertices;
\item for every vertex $s\in V(J)$, $s$ has at most one neighbour in $\alpha$, and at most one in $\beta$;
\item for all $e\in E(J)$,  $a,b\notin M_{e}$;
\item $a$ is complete to $\bigcup_{u\in \alpha}M_u$, and $a$ has no other neighbours in $V(M)$;
$b$ is complete to $\bigcup_{u\in \beta}M_u$, and $b$ has no other neighbours in $V(M)$.
\end{itemize}
\begin{figure}[H]
\centering

\begin{tikzpicture}[scale=.8,auto=left]
\tikzstyle{every node}=[inner sep=5pt, fill=white,circle,draw]
\node (a1) at (-4,1.2) {};
\node (a2) at (-3,0) {};
\node (a3) at (-4,-1.2) {};
\node (b1) at (4,1.2) {};
\node (b2) at (3,0) {};
\node (b3) at (4,-1.2) {};
\node (e1) at (-1,2) {};
\node (e2) at (-1,0) {};
\node (f1) at (1,2) {};
\node (f2) at (1,0) {};
\tikzstyle{every node}=[inner sep=1.5pt, fill=black,circle,draw]
\node (A) at (0,2) {};
\node (B) at (0,0) {};

\foreach \from/\to in {a1/e1,f1/b1,a2/e2,f2/b2,a3/b3}
\draw [line width=12pt] (\from) -- (\to);
\foreach \from/\to in {a1/e1,f1/b1,a2/e2,f2/b2,a3/b3}
\draw [line width=11pt, white] (\from) -- (\to);
\foreach \from/\to in {a1/e1,f1/b1,a2/e2,f2/b2,a3/b3}
\draw [dashed] (\from) -- (\to);

\tikzstyle{every node}=[inner sep=6pt, fill=white,circle,draw]
\node (a1) at (-4,1.2) {};
\node (a2) at (-3,0) {};
\node (a3) at (-4,-1.2) {};
\node (b1) at (4,1.2) {};
\node (b2) at (3,0) {};
\node (b3) at (4,-1.2) {};
\node (e1) at (-1,2) {};
\node (e2) at (-1,0) {};
\node (f1) at (1,2) {};
\node (f2) at (1,0) {};

\tikzstyle{every node}=[inner sep=1.5pt, fill=black,circle,draw]
\node (a11) at (-3.8,1.2) {};
\node (a12) at (-4.2,1.2) {};
\node (a21) at (-3,.2) {};
\node (a22) at (-3,-.2) {};
\node (a31) at (-3.8,-1.2) {};
\node (a32) at (-4.2,-1.2) {};
\node (b11) at (3.8,1.2) {};
\node (b12) at (4.2,1.2) {};
\node (b21) at (3,.2) {};
\node (b22) at (3,-.2) {};
\node (b31) at (3.8,-1.2) {};
\node (b32) at (4.2,-1.2) {};
\node (e11) at (-1,2.2) {};
\node (e12) at (-1,1.8) {};
\node (e21) at (-1,.2) {};
\node (e22) at (-1,-.2) {};
\node (f11) at (1,2.2) {};
\node (f12) at (1,1.8) {};
\node (f21) at (1,.2) {};
\node (f22) at (1,-.2) {};

\foreach \from/\to in {a11/a21,a11/a22,a12/a21,a12/a22,a11/a31,a11/a32,a12/a31,a12/a32,a21/a31,a21/a32,a22/a31,a22/a32}
\draw [-] (\from) -- (\to);
\foreach \from/\to in {b11/b21,b11/b22,b12/b21,b12/b22,b11/b31,b11/b32,b12/b31,b12/b32,b21/b31,b21/b32,b22/b31,b22/b32}
\draw [-] (\from) -- (\to);
\foreach \from/\to in {A/e11,A/e12,A/f11,A/f12,B/e21,B/e22,B/f21,B/f22, A/B}
\draw [-] (\from) -- (\to);
\tikzstyle{every node}=[]
\draw (A) node [above]           {$a$};
\draw (B) node [below]           {$b$};

\begin{scope}[shift ={(-8,0)}]
\tikzstyle{every node}=[inner sep=1.5pt, fill=black,circle,draw]
\node (p1) at (-1.5,0) {};
\node (q1) at (-.5,0) {};
\node (r1) at (-.5,1) {};
\node (p2) at (1.5,0) {};
\node (q2) at (.5,0) {};
\node (r2) at (.5,1) {};
\foreach \from/\to in {p1/q1,p1/r1,p2/q2,p2/r2}
\draw [-] (\from) -- (\to);
\draw (p1) to [bend right=40] (p2);

\end{scope}

\end{tikzpicture}

\caption{The smallest possible $J$, and a $J$-strip system with cross-edge} \label{stripsystem}
\end{figure}
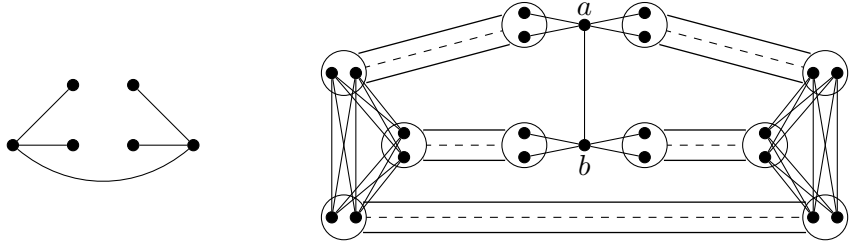

If we are given $J,M$ and $ab$ then we can reconstruct $\alpha, \beta$, so we call $(\alpha, \beta)$ the
{\em corresponding partition}.
If $G$ is an extended tree line-graph $H(T)$ with cross-edge $ab$, where $T$ is a tree, and $J$ is the shape of $T$,
then there is a
$J$-strip system in $G$ with the same cross-edge $ab$, 
defined as follows. Let $(A',B')$ be a bipartition of $T$,
as in the definition of $H(T)$, and let
$\alpha=V(J)\cap A'$, and $\beta=V(J)\cap B'$. 
For each edge $e$ of $J$,
define $M_{e}$ to be the edge-set of the corresponding branch of $T$; and 
for each $u\in V(J)$, let $M_u$ be the set of edges of $T$ incident with $u$.
This defines a $J$-strip system. (Note that some strips might not be proper;
if some branch of $T$ has length one then the $J$-strip system is not proper.)

Let $M$ be a $J$-strip system in $G$ with cross-edge $ab$. 
If $D$ is a subtree of $J$, and we choose an $e$-rung $R_e$ for each $e\in E(J)$, then the subgraph of $G$ induced on
$\bigcup_{e\in E(D)}V(R_e)$, denoted by $R_D$, is the line graph of some tree that has the same shape as $D$. Thus,
$R_D$ depends on the choices of the individual $e$-rungs $R_e$, but we leave this dependence implicit.

Let $M$ be a $J$-strip system in $G$ with cross-edge $ab$ and partition $(\alpha,\beta)$.
We say $X\subseteq V(M)\cup \{a,b\}$
is {\em local} if either:
\begin{itemize}
\item $X\subseteq M_e$ for some $e \in E(J)$; or
\item $X\subseteq M_u$ for some $u\in V(J)$; or
\item $X$ contains $a$ and not $b$, and $X\setminus \{a\}\subseteq M_{u}$ for some leaf $u\in \alpha$; 
or $X$ contains $b$ and not $a$, and $X\setminus \{a\}\subseteq M_{u}$ for some leaf $u\in \beta$.
\end{itemize}
We need a lemma:
\begin{thm}\label{getlocal}
If $X\subseteq V(M)\cup \{a,b\}$ is not local, and $\{a,b\}\not\subseteq X$, then 
there exist $x,y\in X$ such that $\{x,y\}$ is not local.
\end{thm}
\Proof
Suppose first that $a,b\notin X$. Choose $x\in X$, and choose $uv\in E(J)$ such that $x\in M_{uv}$. 
There exists $y\in X\setminus M_{uv}$, and we may assume
that $\{x,y\}$ is local; so we may assume that $x,y\in M_u$. There exists $z\in X\setminus M_u$; and we may assume that
$\{x,z\}$ is local, and so either $z\in M_{uv}$, or $x,z\in M_v$. In either case $\{y,z\}$ is not local, since $J$
is a tree.

Thus we may assume that $a\in X$, and $b\notin X$. Also there exists $x\in X\setminus \{a,b\}$; and we may assume that $\{a,x\}$
is local, and so $x\in M_{u}$ for some $u\in \alpha$.
There exists $y\in X\setminus (M_{u}\cup \{a\})$. Since we may assume that $\{a,y\}$ is local,
$y\in M_{v}$ for some $v\in \alpha$, and so $v\ne u$. 
From the definition of cross-edge, $u,v$ have no common neighbour in $J$, and so
$\{x,y\}$ is not local. 
This proves \ref{getlocal}.~\bbox

We will need two maximizations:
\begin{itemize}
\item
We start with an even-hole-free graph $G$, and an edge $ab$ of $G$, such that there is an extended tree line-graph $H(T)$
that is an induced subgraph of $G$, 
with cross-edge $ab$. Subject to this we choose $T$ with as many branches as possible, that is, such that its shape $J$
has $|E(J)|$ maximum.
\item Then we choose a $J$-strip system $M$ in $G$ with the same cross-edge $ab$, with $V(M)$ maximal.
\end{itemize}
In these circumstances we say that $(J,M)$ is {\em optimal for $ab$}.
Our first goal is to prove:
\begin{thm}\label{treestructplus}
Let $ab$ be an edge of an even-hole-free graph $G$, and let $(J,M)$ be optimal for $ab$.
Let $Z$ be the set of vertices of $G$ adjacent to both $a,b$.
Then for every connected induced subgraph $F$ of $G\setminus (Z\cup V(M))$:
\begin{itemize}
\item if not both $a,b$ have neighbours in $V(F)$, then the set of vertices in $V(M)$
with a neighbour in $V(F)$ is local;
\item if both $a,b$ have neighbours in $V(F)$, then there exists a leaf $t$ of $J$ such that every vertex of $V(M)$
with a neighbour in $V(F)$ belongs to $M_t$.
\end{itemize}
\end{thm}
We break the proof into three steps, \ref{treestructplus1}, \ref{treestructplus2}, and \ref{treestructplus3} below,
depending on the number of $a,b\in N(F)$.

Under
the hypotheses of \ref{treestructplus}, let $(\alpha,\beta)$ be the corresponding partition.
Let us say that a subgraph $F$ is {\em small} if $F$ is connected and $F$ is an induced subgraph of $G\setminus (Z\cup V(M))$;
and a {\em small component} is a component of $G\setminus (V(M)\setminus Z)$. A small subgraph $F$ is {\em $\alpha$-peripheral}
if $X(F)\subseteq M_t$ for some $t\in \alpha$. We define {\em $\beta$-peripheral} similarly; and $F$ is {\em peripheral} if it is
either $\alpha$- or $\beta$-peripheral.
If $F$ is small, the set of vertices in $V(M)$
with a neighbour in $V(F)$ is denoted by $X(F)$.
We begin with:
\begin{thm}\label{treestructplus1}
Under the hypotheses of \ref{treestructplus}, if $F$ is small, and
$a,b\notin N(F)$,
then $X(F)$ is local.
\end{thm}
\Proof
Suppose the theorem is false, and choose a small subgraph $F$
not satisfying the theorem, with $F$ minimal.
By \ref{getlocal}, there exist
$x,y\in X(F)$ such that $\{x,y\}$ is not local, and so $F$ is a path joining these two vertices. Let $F$ have ends $f_1,f_2$.

For $x_1,x_2\in V(M)$, let us say $s\in V(J)$ {\em separates} $x_1,x_2$
if $x_1,x_2\notin M_s$, and $s$ lies on the path of $J$ between $e_1,e_2$, where $x_i\in M_{e_i}\;(i=1,2)$.
\\
\\
(1) {\em If $x_1,x_2\in X(F)$, there is no $s\in V(J)$ that separates $x_1,x_2$.}
\\
\\
Let $x_i\in M_{e_i}\;(i=1,2)$, and suppose that $s\in V(J)$ separates $x_1,x_2$. Then $\{x_1,x_2\}$ is not local,
and so we may assume that $f_1x_1$ and $f_2x_2$ are edges. Choose three leaf-paths $S_1,S_2,S_3$
of $J$, each with one end $s$ and otherwise pairwise vertex-disjoint, with $e_i\in E(S_i)$ for $i = 1,2$.
For $i = 1,2,3$ let $s_i$ be the edge of $S_i$ incident with $s$.
For $i = 1,2,3$ and each $e\in E(S_i)$, 
choose an $e$-rung $R_e$, such that $x_i\in V(R_{e_i})$ for $i = 1,2$. 
For $i = 1,2,3$, let $u_i$ be the end of $R_{s_i}$ in $M_s$.
Then $R_{S_i}$
is an induced path of $G$ from $u_i$ to some $p_i\in N(\{a,b\})$.
We may assume that 
$x_1,x_2$ have been chosen such that for $i = 1,2$ the subpath of $R_{S_i}$ between $x_i, p_i$ is minimal.

Suppose
that there exists $x_3\in X(F)$, in $V(R_{S_1}\cup R_{S_2}\cup R_{S_3})$ and different from and nonadjacent to $x_1,x_2$. Choose $x_3$ such that
the subpath of $R_{S_3}$ between $x_3$ and $p_3$ is minimal.
We claim that $|V(F)|=1$. For if not, we may assume that $x_3$ has a neighbour in $V(F\setminus f_2)$, and since 
$X(F\setminus f_2)$ is local (from the minimality of $F$) and contains $x_1,x_3$, and $x_1,x_3$ are nonadjacent, it follows that 
$X(F\setminus f_2)\subseteq M_{e_1}$, and in particular $x_3$ belongs to $R_{e_1}$. But then there is an induced path 
between the ends of $R_{e_1}$ and contained in $G[V(R_{e_1}\cup (F\setminus f_2))]$, that contains at least one vertex of $F\setminus f_2$,
and the vertices of this path can be added to $M_{e_1}$, contrary to the maximality of $V(M)$. This proves that $|V(F)|=1$.

If $p_1,p_2,p_3\in N(a)$, there is a theta with ends $f_1,a$ and constituent paths
$f_1\dd x_i\dd R_{S_i}\dd p_i\dd a$ for $i = 1,2,3$; and similarly not all $p_1,p_2,p_3\in N(b)$. By exchanging $a,b$ if necessary,
we may assume that two of $p_1,p_2,p_3\in N(a)$; then there is a theta with ends $f_1,a$ with 
constituent paths
$f_1\dd x_i\dd R_{S_i}\dd p_i\dd a$ for the two values of $i$ with $p_i\in N(a)$, and $f_1\dd x_i\dd R_{S_i}\dd p_i\dd b\dd a$
for the value of $i$ with $p_i\in N(b)$.

This proves that $X(F)\cap V(R_{S_3})=\emptyset$, and every vertex of $X(F)\cap V(R_{S_1}\cup R_{S_2})$ is equal or adjacent to one of 
$x_1,x_2$. For $i=1,2$, let $y_i$ be the neighbour of $x_i$ in $R_{S_i}$ between $x_i$ and $u_i$ (this exists, since 
$x_i\notin M_s$.)
The path $R_{S_1\cup S_2}$
can be completed to a hole $H$ by adding $a$ or $b$ or both.
From the minimality of $F$, $X(F\setminus \{f_1,f_2\})=\emptyset$. 
We claim that the only edges between $V(R_{S_1}\cup R_{S_2}\cup R_{S_3})$ and $V(F)$ are the edges $f_1x_1,f_2x_2$ and exactly one
of $f_1y_1,f_2y_2$.
If $|V(F)|=1$ this is true since $f_1$ cannot have two nonadjacent neighbours in $H$, or four neighbours in $H$.
If $f_1\ne f_2$ then from the minimality
of $F$, $f_1$ is nonadjacent to $y_2,x_2$, and $f_2$ is nonadjacent to $x_1,y_1$;
at least one of the pairs $f_1y_1,f_2y_2$ is an edge since otherwise the subgraph induced on $V(H)\cup V(F)$
is a theta,
and not both since otherwise the same subgraph is a prism. (Note that $y_1\ne y_2$ since $x_1,x_2\notin M_s$.)
Thus we may assume that $f_1x_1,f_1y_1,f_2x_2$ are edges, and there are no other edges between 
$V(R_{S_1}\cup R_{S_2}\cup R_{S_3})$ and $V(F)$. If $x_2\notin N(\{a,b\})$, we may assume that at least two of $p_1,p_2,p_3$
are adjacent to $a$, and then there is a theta between $x_2$ and $a$ with constituent paths
$$x_2\dd R_{S_2}\dd u_2\dd u_3\dd R_{S_3}\dd p_3\dd a,$$
$$x_2\dd f_2\dd F\dd f_1\dd x_1\dd R_{S_1}\dd p_1\dd a,$$
$$x_2\dd R_{S_2}\dd p_2\dd a,$$
inserting $b$ before $a$ in one of these paths if necessary.
Thus $e_2$ is a leaf-edge of $J$, and $x_2=p_2\in N(\{a,b\})$, and we may assume that $x_2\in N(b)$. 
We can choose $S_3$ such that it has an end in $\alpha$
(from the definition of a crossedge for a tree strip system), and hence we may assume that $p_3\in N(a)$.
If $p_1\in N(a)$ then the same argument gives a theta, which is impossible; so we may assume that every choice of $S_1$
has an end in $\beta$, and so $e_1$ is also a leaf-edge of $J$. Let $r$ be the end of $e_1$ that is not a leaf of $J$,
and let $t$ be the end of $e_2$ that is not a leaf. 
From the definition of a crossedge, $r\ne t$.
Exactly two vertices of $R_{S_1\cup S_2}$
belong to $M_r$, and they are adjacent, say $r_1,r_2$; and define $t_1,t_2$ similarly, where $r_1,r_2,t_1,t_2$
are in order in $R_{S_1\cup S_2}$ (possibly $r_2=t_1$). By choosing a leaf-path of $J$ with one end $r$ that is edge-disjoint from $S_1,S_2$,
and has an end in $\beta$,
and choosing a rung for each of its edges, we find a path $R$ say of $G[V(M)]$ with ends $r_3, r_0$ say, where
$r_3$ is adjacent to $r_1,r_2$, and $r_0\in N(b)$ and there are no other edges between $V(R)$ and $V(R_1\cup R_2)$,
Define  a path $T$ with ends $t_3,t_0$ similarly, where $t_3$ is adjacent to $t_1,t_2$ and $t_0\in N(b)$, and there
are no edges between $V(R)$ and $V(T)$. There is a near-prism with bases $\{r_1,r_2,r_3\}$, $\{t_1,t_2,t_3\}$
and constituent paths 
$$r_1\dd R_{S_1}\dd y_1\dd f_1\dd F\dd f_2\dd x_2\dd R_{S_2}\dd t_2,$$
$$r_3\dd R\dd r_0\dd b\dd t_0\dd T\dd t_3,$$
$$r_2\dd R_{S_1\cup S_2}\dd t_1,$$
contrary to \ref{subgraphs} (note that possibly $r_2=t_1$). This proves (1).
\\
\\
(2) {\em There is an edge $uv$ of $J$ such that $X(F)\subseteq M_u\cup M_v\cup M_{uv}$.}
\\
\\
Suppose first that for some $uv\in E(J)$, there exists $x\in X(F)\in M_{uv}\setminus (M_u\cup M_v)$.
Then for each $y\in X(F)$, (1) implies that no vertex of $S$ separates $x,y$, and so $y\in M_u\cup M_v\cup M_{uv}$ as required.
Thus we may assume that $X(F)\subseteq \bigcup_{v\in V(J)}M_v$. 
Suppose next that some $x\in X$ belongs to $M_v$ for only one value of $v\in V(J)$. 
Let $y\in X(F)\setminus M_v$, and choose $u\in V(J)$ with $y\in M_u$.
Let $w$ be the neighbour of $v$ in $J$ on the path of $J$ between $v,u(y)$.
Since $w$ does not separate $x,y$, and $x\notin M_w$ from the choice of $x$, it follows that $y\in M_w$; define $w(y)=w$.
If there exist $y_1,y_2\in X(F)\setminus M_v$ with $w(y_1)\ne w(y_1)$, then $v$ separates $y_1,y_2$, contrary to (1).
So there exists a neighbour $w$ of $v$ in $J$, such that $y\in M_w$ for all $y\in X(F)\setminus M_v$; and so the claim holds.

We may therefore assume that every vertex in $X(F)$ belongs to $M_v$ for two vertices $v\in V(J)$. For each 
$x\in X(F)$, let $e(x)$ be the edge $uv$ of $J$ such that $x\in M_u\cap M_v$. If all the edges $e(x)\;(x\in X(F))$
have a common end, then the claim holds; so we may assume that there exist $x_1,x_2\in X(F)$ such that $e(x_1), e(x_2)$
have no common end. Let $e(x_i)=u_iv_i$ for $i = 1,2$; thus $u_1,v_1,u_2,v_2$ are distinct vertices of $J$.
Since no vertex of $J$ separates $x_1,x_2$ by (1), it follows that one of $u_1,v_1$ is $J$-adjacent to one of $u_2,v_2$;
say $v_1,v_2$ are $J$-adjacent, and so $u_1\dd v_1\dd v_2\dd u_2$ is a path of $J$. Suppose there exists $x_3\in X(F)$
such that $x_3\notin M_{v_1}\cup M_{v_2}$. Let $e(x_3)=u_3v_3$ say. Thus $u_3,v_3\ne v_1,v_2$; and we may assume that
$v_1$ lies on the path of $J$ between $v_2$ and $u_3$, by exchanging $x_1,x_2$ if necessary. But then $v_1$
separate $x_2,x_3$, contrary to (1). This proves (2).
\\
\\
(3) {\em There is an edge $uv$ of $J$ such that $X(F)\subseteq M_u\cup M_{uv}$.}
\\
\\
Suppose not. Choose $uv$ as in (2); then there exist $x_1,x_2\in X(F)$ with 
$x_1\in M_u\setminus M_{uv}$ and $x_2\in M_v\setminus M_{uv}$. 
We may assume that $f_1x_1$ and $f_2x_2$ are edges.
From the minimality of $F$, there are no edges between $V(F\setminus f_1)$ and $(M_u\cup M_{uv})\setminus M_{v}$,
and no edges between $V(F\setminus f_2)$ and $(M_v \cup M_{uv})\setminus M_{u}$.

Let $c_1\ll c_k$ be the edges of $J$
incident with $u$, and different from $uv$; and let $d_1\ll d_{\ell}$ be those incident with $v$ and different from $uv$.
Thus $k,\ell\ge 2$. 
If $f_1$ is complete to $M_u\setminus M_{uv}$ and $f_2$ is complete to $M_v\setminus M_{uv}$,
we can add $f_1$ to $M_u$, add $f_2$ to $M_v$, and add $V(F)$ to $M_{uv}$, contrary to the maximality of $V(M)$.
Thus we may assume that $f_1$ has a non-neighbour in $M_u\setminus M_{uv}$; and since $x_1\in M_u\setminus M_{uv}$ and 
$k\ge 2$, we may assume that $x_1\in M_{c_1}\cap M_u$ and $y_1\in M_{c_2}\cap M_u$, and $f_1,y_1$ are nonadjacent. 
Also we may assume $x_2\in M_{d_1}\setminus M_{uv}$.
For $1\le i\le k$ choose a leaf-path $C_i$ of $J$ from $u$ and using $c_i$; and for $1\le i\le \ell$ define $D_i$ similarly;
and choose an $e$-rung $R_e$ for each of their edges $e$, containing $x_1,x_2,y_1$.
If $x_1\notin N(\{a,b\})$, we may assume that at least two of $C_1,C_2,D_1$ have an end in $\alpha$;
and then there is a theta in $G$ with ends $x_1,a$ and constituent paths 
$$x_1\dd R_{C_1}\dd a,$$
$$x_1\dd y_1\dd R_{C_2}\dd a,$$
$$x_1\dd f_1\dd F\dd f_2\dd x_2\dd R_{D_1}\dd a,$$
inserting $b$ into one of these if necessary. Thus we may assume that $x_1\in N(b)$ say. Consequently $c_1$ has an end in $\beta$;
and so $C_2$ can be chosen with an end in $\alpha$. If also $D_1$ can be chosen with an end in $\alpha$
then the same construction still gives a theta; so the leaf of $D_1$ is in $\beta$. 
Hence the leaf of $D_2$ is not in $\beta$, so $f_1$ has no neighbour in
$M_{d_2}$. This restores the symmetry between $u,v$.
Let $R_{uv}$ be a $uv$-rung, with ends   
$r_1\in M_u\cap M_{uv}$ and $r_2\in M_v\cap M_{uv}$. Since $b$ is adjacent to both $x_1,x_2$, it follows that $f_1\ne f_2$,
and for the same reason, $r_1\ne r_2$. From the minimality of $F$, the only edges between $V(F)$ and $V(R_{uv})$
are either $f_1r_1$ or $f_2r_2$; since $x_1\dd r_1\dd R_{uv}\dd r_2\dd x_1$ and $x_1\dd f_1\dd F\dd f_2\dd x_2$ are both odd,
exactly one of these two edges is present, say $f_1r_1$ (without loss of generality, since the symmetry between $u,v$
was restored). But then there is a theta with ends $r_1,x_2$ and constituent paths
$$r_1\dd R_{uv}\dd r_2\dd x_2,$$
$$r_1\dd f_1\dd F\dd f_2\dd x_2,$$
$$r_1\dd y_1\dd R_{C_2}\dd a\dd b\dd x_2,$$
contrary to \ref{subgraphs}. This proves (3).

\bigskip

Choose $uv$ as in (3). Let $c_1\ll c_k$ be the edges of $J$
incident with $u$, and different from $uv$.
Since $X(F)$ is not local, there exists $x_1\in M_u\setminus M_{uv}$
and $x_2\in M_{uv}\setminus M_u$. We may assume that $f_1x_1$ and $f_2x_2$ are edges. From the minimality of $F$,
there are no edges between $V(F\setminus f_2)$ and $M_{uv}\setminus M_u$, and none between $V(F\setminus f_1)$ and
$M_u\setminus M_{uv}$. If $f_1$ is complete to $M_u\setminus M_{uv}$, we can add $f_1$ to $M_u$ and $V(F)$ to $M_{uv}$,
contrary to the maximality of $V(M)$; so we may assume that
$x_1\in M_{c_1}\cap M_u$ and $y_1\in M_{c_2}\cap M_u$, and $f_1,y_1$ are nonadjacent. 
For $1\le i\le k$ choose a leaf-path $C_i$ of $J$ from $u$ and using $c_i$. 
Choose a leaf-path $D$ of $J$ from $u$ and using $uv$. (Possibly $D$ has length one.)
For each edge $e$ of $C_1\ll C_k, D$ choose an $e$-rung $R_e$, where $R_{c_1}$ contains $x_1$, $R_{c_2}$
contains $y_1$, and $R_{uv}$ contains $x_2$. 

Suppose that $x_1\notin N(\{a,b\})$; then we may assume that at least two of $C_1,C_2, D$ have an end in $\alpha$;
and then there is a theta in $G$ with ends $x_1,a$ and constituent paths
$$x_1\dd R_{C_1}\dd a,$$
$$x_1\dd y_1\dd R_{C_2}\dd a,$$
$$x_1\dd f_1\dd F\dd f_2\dd x_2\dd R_{D}\dd a,$$
inserting $b$ into one of these if necessary. Thus we may assume that $x_1\in N(b)$ say. Hence $C_2$ and $D$ can be chosen to 
have an end in $\alpha$, and the same construction still serves to find a theta, a contradiction.
This proves \ref{treestructplus1}.~\bbox

\begin{thm}\label{treestructplus2}
Under the hypotheses of \ref{treestructplus}, if $F$ is small, and
$a\in N(F)$ and $b\notin N(F)$,
then $F$ is $\alpha$-peripheral.
\end{thm}
\Proof
Suppose the theorem is false, and choose a small subgraph $F$
not satisfying the theorem, with $F$ minimal. 
By \ref{getlocal}, there exist
$x,y\in X(F)\cup \{a\}$ such that $\{x,y\}$ is not local, and so $F$ contains a path joining these two vertices; 
and $a$ has a neighbour in this path, by \ref{treestructplus1}, and so $F$ is this path, from the minimality of $F$.
Let $F$ have ends $f_1,f_2$.
\\
\\
(1) {\em Let $D$ be a path of $J$ with distinct ends both in $\beta$, and for each $e\in E(D)$ choose an $e$-rung $R_e$.
Then either $X(F)$
contains no vertices of $R_D$, or it contains exactly two and they are adjacent.}
\\
\\
Let the ends of $D$ be $t_1,t_2\in \beta$.
Since $R_d$ has both ends in $N(b)$, it follows that $a$ has no $G$-neighbours
in $V(R_d)$; and by adding $b$ to $R_D$ we obtain a hole $H$, and so $a$ has a unique $G$-neighbour $b$ in $V(H)$. 
We may assume there exists $y\in V(H)\cap X(F)$; and since $\{a,y\}$ is not local, the minimality of $F$ implies that
$F$ is a path between $a,y$; say $a$ is adjacent to $f_1$ and to no other vertex of $F$, and $y$ is adjacent to $f_2$ and
to no other vertex of $F$. For the same reason, $F\setminus f_2$ is anticomplete to $V(H)$.

If $f_2$ has two nonadjacent vertices in $V(H)$, there are two paths $P_1,P_2$ between $f_2,b$ with interior in $V(H)$,
and with union a hole; but then 
there is a theta with ends $f_2,b$ and constituent paths 
$$f_2\dd F\dd f_1\dd a\dd b,$$
$$f_2\dd P_1\dd b,$$
$$f_2\dd P_2\dd b,$$
a contradiction. 

If $f_2$ has a unique neighbour in $V(H)$, say $x$, and $x$ is nonadjacent to $b$,
then $G[V(H\cup F)]$ is a theta
with ends $x,b$, again a contradiction. 

Suppose next that $f_2$ has a unique neighbour in $V(H)$, say $x$, and
$x$ is adjacent to $b$. Let $x\in M_{t_1}$, say, and let $s_1$ be the neighbour of $t_1$ in $J$.
Since $a\dd b\dd x\dd f_2\dd a$ is not a 4-hole, it follows that $a,f_2$ are not adjacent, and therefore $f_1\ne f_2$,
and so $a\in N(F\setminus f_2)$.
From the minimality of $F$, $X(F\setminus f_2)\cup \{a\}$ is local.
Choose $t_3\in \alpha$ such that $X(F\setminus f_2)\cap M_e=\emptyset$ (this is possible, 
since $|\alpha|\ge 2$ and $X(F\setminus f_2)\cup \{a\}$ is local).
Let $D_3$ be a path of $J$, edge-disjoint from $D$ and with ends $d,t_3$ where $d\in V(D)$. For each $e\in E(D_3)$ 
choose an $e$-rung $R_e$. Let $D_1,D_2$ be the subpaths of $D$ with ends $d$ and $t_1,t_2$ respectively.

If $F$ is anticomplete to $R_{D_3}$, there is a theta with ends $x,a$, and constituent paths
$$x\dd b\dd a,$$
$$x\dd f_2\dd F\dd f_1\dd a,$$
$$x\dd R_{D_1\cup D_3}\dd a,$$
contrary to \ref{subgraphs}. Thus $F$ is not anticomplete to $R_{D_3}$. Now $F$, $R_{D_3}$ are vertex-disjoint, 
since $V(F)$ is disjoint
from $V(M)$ and $V(R_{D_3})\subseteq V(M)$. Let $y\in V(R_{D_3})$ with a neighbour in $V(F)$. If $y$ has a neighbour in 
$V(F\setminus f_2)$, then $\{a,y\}$ is local, from the minimality of $F$; but then $y \in N(a)$ and so $y\in M_{t_3}$,
contrary to the choice of $e_3$. 
Thus $y$ is adjacent to $f_2$ and to no other
vertex of $F$. From the minimality of $F$, $\{x,y\}$ is local; and so either $y\in M_{s_1t_1}$ or $x,y\in M_{s_1}$.
The first is impossible since $s_1t_1$ is not an edge of $D_3$; and so $x,y\in M_{s_1}$. In particular, $d=s_1$, and $y$
is the end of $R_{D_3}$ in $M_d$. But then 
$$y\dd R_{D_2}\dd b\dd a\dd f_1\dd F\dd f_2\dd y$$
is a hole, 
in which $x$ has exactly four neighbours, making a 4-wheel, a contradiction.
This proves (1).

\bigskip

Let $X_1$ be the set of $x\in X(F)\cap V(M)$ such that $x\in M_e$ for some $e\in E(J)$ not incident with any vertex in $\alpha$,
and let $X_2=X(F)\setminus X_1$.
From the minimality of $F$, there are no edges between $V(F\setminus f_2)$ and $X_1$.
\\
\\
(2) {\em $X_1\ne \emptyset$.}
\\
\\
Suppose that $X_1=\emptyset$. Consequently the only edges $e\in E(J)$ with $X(F)\cap M_e\ne \emptyset$ are those 
with an end in $\alpha$.
Suppose that there are distinct $e_1,e_2$, both with an end in $\alpha$, such that $X(F)\cap M_{e_i}\ne \emptyset$ for $i = 1,2$. 
Let $e_i=s_it_i$ where $t_i\in \alpha$ for $i = 1,2$.
Let $D$ be a path of $J$ with both ends in $\beta$, containing $s_1$ and $s_2$. Let $D$ have end-edges $t_1', t_2'$,
where $t_1',s_1,s_2,t_2'$ are in order on $D$.
For $i = 1,2$ let $D_i$ be the subpath of $D$ between $t_i', s_i$; and let $D_3$ be the subpath between $s_1,s_2$.
For each $e\in E(D)\cup \{e_1,e_2\}$ choose an $e$-rung $R_e$,
with $V(R_{e_i})\cap X(F)$ nonempty for $i = 1,2$.
Let the ends of $R_{e_1}, R_{D_1}, R_{D_3}$ in $M_{s_1}$ be $p_1,p_2,p_3$ respectively, and let the ends of
$R_{e_2}, R_{D_2}, R_{D_3}$ in $M_{s_2}$ be $q_1,q_2,q_3$ respectively.
Now $F$ is anticomplete to $R_D$, since $X_1=\emptyset$.
Since $X(F)$ meets both $R_{e_1}, R_{e_2}$, there is an induced path $Q$ between $p_1,q_1$ with interior in 
$V(R_{e_1}\cup R_{e_2}\cup F)$.
There is a near-prism in $G$ with bases
$\{p_1,p_2,p_3\}, \{q_1,q_2,q_3\}$ and constituent paths
$$p_2\dd R_{D_1}\dd b\dd R_{D_2}\dd q_2,$$
$$p_1\dd Q\dd q_1,$$
$$p_3\dd R_{D_3}\dd q_3,$$
a contradiction. 

Consequently there is a unique $e\in E(J)$ such that $X(F)\cap M_e\ne \emptyset$, say $e=st$ where $t\in \alpha$.
Since $X(F)$ does not satisfy the theorem, it follows that
$X(F)\not\subseteq M_t$; let $x\in X(F)\setminus M_t$. Since $\{a,x\}$ is not local, we may assume that
$f_1a$ and $f_2x$ are edges. But then we can add $V(F)$ to $M_e$ and $f_1$ to $M_t$, contrary to the maximality of $V(M)$. 
This proves (2).

\bigskip

For each edge $e\in E(J)$, choose an $e$-rung $R_e$. The subgraph induced on $\bigcup_{e\in E(J)}V(R_e)$
is the line-graph $L(T)$ of a tree $T$, where $T$ has shape $J$, and $E(T)=\bigcup_{e\in E(J)}V(R_e)$.
In particular, $E(T)=V(R_J)$, and $V(J)$ is the set of branch-vertices of $T$.
Let us call such a tree $T$ a {\em realization} of $M$. If $P$ is a subgraph of $T$, then $E(P)$ is a set of vertices of $G$,
and we denote $G[E(P)]$ by $L(P)$ (it is indeed the line-graph of $P$).
\\
\\
(3) {\em For every realization $T$
with $E(T)\cap X_1\ne \emptyset$, there exists $d\in V(T)$ such that
$X_1\cap E(T)$ consists of all edges of $T$ incident with $d$ that belong to branches of $T$ that do 
not have an end-edge in $N(a)$.}
\\
\\
Let $P$
be a path of $T$ with distinct ends, and both end-edges in $N(b)$,
with $E(P)\cap X_1\ne \emptyset$. By (1) there exists $d\in V(P)$
such that $X(F)\cap E(P)$ is the set of edges of $P$ that are $T$-incident with $d$. We will show that $d$ satisfies the claim.
Let $P_1,P_2$ be the two subpaths of $P$
between $d$ and an end of $P$, and let $x_1,x_2$ respectively be the edges of $P_1,P_2$ that are $T$-incident with $d$.
Suppose that $x_3\in E(T)\cap X(F)$; we will show that $x_3$ is incident with $d$ in $T$. We may assume that
$x_3\ne E(P)$. Let $e_3\in E(J)$ with $x_3\in M_{e_3}$. Since $x_3\notin X_2$, 
there is a path of $J$ with both ends in $\beta$ containing $e_3$; and hence there is a path of $T$ containing $x_3$
with both end-edges in $N(b)$. Choose a path $P_3$ of $T$ containing $x_3$ with one end-edge in $N(b)$
and the other in $V(P)$, edge-disjoint from $P$. Let $p$ be the end of $P_3$ in $V(P)$; and let $P_1',P_2'$
be the two subpaths of $P$ between $p$ and the ends of $P$. If $p\ne d$, then $d$ is an internal vertex of one of $P_1',P_2'$,
say $P_1'$; and $X(F)$ contains two nonconsecutive edges of the path $P_1'\cup P_3$, contrary to (1). So $p=d$.
From (1) applied to the path $P_1\cup P_3$, it follows that there is a unique edge of $P_3$ in $X(F)$, and it is
$T$-incident with $d$. This proves that all edges of $T$ in $X(F)$ are $T$-incident with $d$. 

Next we show that every edge of $T$ that is $T$-incident with $d$, and not in a branch of $T$ with end-edge in $N(a)$,
belongs to $X_1$.
Let $y$ be an edge of $T$ that is $T$-incident with $d$, and let $y\in M_e$ say, with no end in $\alpha$.
We must show that $y\in X(F)$. To see this, choose a path $P_3$ of $T$ containing $y$
with one end-edge in $N(b)$ and one end $d$,
edge-disjoint from $P$. From (1) applied to $P_1\cup P_3$ it follows that $y\in X(F)$. 
This proves (3).
\\
\\
(4) {\em Let $T,d$ be as in (3). Then there is a branch $S$ of $T$ with one end $d$ and with an end-edge in $N(a)$,
such that $X_2\cap E(T)\subseteq E(S)$.
In particular $d\in V(J)$, and so $X_1\cap E(T)\subseteq M_d$.}
\\
\\
If $X_2\cap E(T)=\emptyset$, we can assume there is no branch $S$ of $T$ with one end $d$ and with an 
end-edge in $N(a)$ (for otherwise the claim holds); and then by (3),
$X(F)\cap E(T)$ consists of all edges of $T$ incident with $d$, and
the subgraph of $G$ induced on $E(T)\cup V(F)\cup \{a,b\}$ is
an extended tree line-graph $H(T')$
with cross-edge $ab$,
for some tree $T'$ whose shape has more edges than $J$, contrary to the choice of $J$.
Thus we may assume that $X_2\cap E(T)\ne \emptyset$.
Let $t\in \alpha$ with $J$-neighbour $s$, such that the branch, $S$ say, of $T$ with ends $s,t$ contains
an edge in $X(F)$. If $s=d$ for every such choice of $t$, then the claim holds (because there
is at most one leaf of $J$ in $\alpha$ $J$-adjacent to $d$). Thus we may assume that $s\ne d$.
Let $P$ be a path of $T$, including the subpath of $T$ between $s,d$,
and with both end-edges in $N(b)$. Now
$P$ is divided into three subpaths by $s,d$, namely from an end of $P$ to $s$, from $s$ to $d$, and from $d$ to the other end
of $P$. We call these $P_1,P_2,P_3$ respectively. Let $d_1,d_2,d_3$ be the edges of $T$ incident with $s$ that belong to
$E(P_1), E(P_2), E(S)$ respectively. Thus exactly one of $x_1,x_2$ belongs to $E(P_2)$, say $x_1$.
Since there are edges between $V(F)$ and $V(L(S))$, there is an induced path $Q$ between $d_3, f_2$ with 
interior in $V(L(S)\cup F)$. 
Then there is a near-prism in $G$ with bases $\{d_1,d_2,d_3\}$, $\{f_2,x_1,x_2\}$ and
constituent paths
$$ d_1\dd L(P_1)\dd b\dd L(P_3)\dd x_2,$$
$$d_3\dd Q\dd f_2,$$
$$d_2\dd L(P_2)\dd x_1,$$
contrary to \ref{subgraphs}. This proves (4).
\\
\\
(5) {\em Let $T,d,S$ be as in (4), and let $S$ have ends $s,d$ say; then $X_2\subseteq M_{sd}$.}
\\
\\
For each $e\in E(J)$, let $R_e$ be the $e$-rung used to 
define $T$. 
If some vertex $x\in X_2$ belongs to $M_e$ say where $e\in E(J)$, then $e$ has an end in $\alpha$ from the definition
of $X_2$, and if $e\ne sd$, we could replace $R_e$ with an $e$-rung that contains $x$, to obtain
a realization that violates (4). This proves (5).
\\
\\
(6) {\em Let $T,d,S$ be as in (4), and let $S$ have ends $s,d$ say; then $X_1= M_d\setminus M_{sd}$.}
\\
\\
There are at least two edges $e_1,e_2$ of $J$,  $J$-incident with $d$ and with no end in $\alpha$;
let $x_1,x_2$ be the edges of the corresponding branches of $T$ that are $T$-incident with $d$.
We show first that $X_1\subseteq M_d\setminus M_{sd}$.
Let $x\in X_1$, and let 
$x\in M_e$ where $e\in E(J)$. Let $R'_e$ be an $e$-rung containing $x$. 
Let $T'$ be the realization of $M$ obtained by replacing $R_e$ by $R'_e$, and otherwise
using all the same rungs. 
Since $e_1\ne e_2$ we may assume that $e\ne e_2$; and so $x_2,x\in V(T')$.
Hence by (4) applied to $T'$, $e_2,e$ have a common end $d'\in V(J)$, and $x_2,x\in M_{d'}$. 
Also either $e=e_1$ or $x_1\in E(T')$;
and so in either case $e_1$ is incident with $d'$. Consequently $d'$ is the common end of $e_1,e_2$ in $J$,
and so $d'=d$. This proves that $x\in M_d$, and so $X_1\subseteq M_d\setminus M_{sd}$.

Next we show that $M_d\setminus M_{sd}\subseteq X_1$.
To see this, let $y\in M_d\setminus M_{sd}$.
Let $e\in E(J)$ with $y\in M_e$; since $y\notin M_{sd}$ it follows that $e$ has no end in $\alpha$.
Let $R'_e$ be an $e$-rung containing $y$.
Since $y\in M_d$ it follows that $e$ is $J$-incident with $d$.
Let $T'$ be the realization obtained by replacing $R_e$ by $R_e'$. Since $e_1\ne e_2$ we may assume that $e\ne e_2$.
Since $e$ has no end in $\alpha$, there is a path $P'$ of $T'$ with $x_2,y\in E(P')$
and with both end-edges in $N(b)$; and so $X_1$ contains either zero or two consecutive edges in this path, by (1).
Not zero, since $x_2 \in E(P')$; so a unique vertex of $R_e'$ belongs to $X_1$, and that vertex is in $M_d$.
Since $y$ is the only vertex of $R_e'$ in $M_d$, it follows that $y\in X_1$. This proves (6).

\bigskip

From (5) and (6) we can add $f_2$ to $M_d$
and add $f_1$ to $M_s$, and add $V(F)$ to $M_{sd}$, contrary to the maximality of $V(M)$. This proves \ref{treestructplus2}.~\bbox

\begin{thm}\label{treestructplus3}
Under the hypotheses of \ref{treestructplus}, if $F$ is small, and
$a,b\in N(F)$, then $F$ is peripheral.
\end{thm}
\Proof 
We claim first:
\\
\\
(1) {\em $X(F)\subseteq N[\{a,b\}]$.} 
\\
\\
Suppose $x\in V(M)$ has a neighbour in $V(F)$, and $x\notin N(\{a,b\})$. Choose a minimal path $P$ of $F$
such that $x$ and at least one of $a,b$ has a neighbour in $V(P)$. Thus $P$ has one end adjacent to $x$ and the other
to $a$, say. But $a,b$ have no common neighbour in $V(F)$, since $V(F)\cap Z=\emptyset$; and so from the minimality of $P$,
$b$ has no neighbour in $V(P)$. But then $P$ violates \ref{treestructplus2}. This proves (1).
\\
\\
(2) {\em Either $X(F)\subseteq N[a]$ or $X(F)\subseteq N[b]$.}
\\
\\
Suppose not; then there is a vertex $c\in V(M)\cap N[a]$ and $d\in V(M) \cap N[b]$, joined by a path $P$ with interior in $V(F)$.
Choose $c,d$ and $P$ such that $P$ has minimum length. Choose $u\in \alpha$ and $v\in \beta$ with $c\in M_u$ and $d\in M_v$, 
and let $D$ be a path
of $J$ with ends $u,v$. Let $p,q$ be the neighbours in $P$ of $c,d$ respectively. 
Let $c_1\ll c_k$ be the vertices of $N(a)\cap V(P)$ in order on $P$, with $c_1=c$. Note that $c_1\ll c_k$
are not adjacent to $b$ since $V(P)\cap Z=\emptyset$. For $1\le i<k$, let $P_i$ be the subpath of $P$ between $c_i$ and $c_{i+1}$.
Since $a\dd c_i\dd P_i\dd c_{i+1}\dd a$ is a hole, and $b$ is adjacent to $a$ and not to $c_i, c_{i+1}$, it follows that
$b$ has an even number of neighbours in $P_i$. Choose $u'\in \alpha\setminus \{u\}$ and $c'\in M_{u'}$. 
By \ref{treestructplus1} and \ref{treestructplus2}, $p$ has no neighbour
in $M_{u'}$, since $p$ has a neighbour in $M_u$ and $X(p)$ is local;
and by the minimality of $P$, no vertex of $P$ different from $p$ has a neighbour in $M_{u'}$. In particular $c_k,c'$ are nonadjacent.
Let $S'$
be an induced path of $G$ between $c', d$ with interior in $V(M)\setminus N[\{a,b\}]$. Then $a\dd c'\dd S'\dd d\dd P\dd c_k\dd a$
is a hole (note that $c_k$ is not adjacent to $c'$), and $b$ has at least two nonadjacent neighbours in it ($a$ and $d$), and so it has an odd number; and therefore
$b$ has an even number of neighbours in the subpath of $P$ between $c_k, d$. Hence $b$ has an even number of neighbours in $V(P)$
altogether. 
Also, $d\dd P\dd c\dd R_D\dd d$ is a hole,
and $b$ has an even number of neighbours in it, at least one; and it has exactly two and they are adjacent. Consequently
$b$ is adjacent to $q$ and has no other neighbours in $V(P)$ except $d$. Similarly $a$ is adjacent to $c,p$ and has 
no other neighbours in $V(P)$. But then the subgraph induced on $V(R_D)\cup V(P)\cup \{a,b\}$ is a prism, a contradiction.
This proves (2).

\bigskip

From (2) we may assume that $X(F)\subseteq N[a]$. Suppose that there exist distinct $u,u'\in \alpha$ such that 
$X(F)\cap M_u, X(F)\cap M_{u'}\ne\emptyset$. Choose $c\in X(F)\cap M_{u}$ and $c'\in X(F)\cap M_{u'}$, such that there
is an induced path $P$ between $c,c'$ with interior in $V(F)$. Both ends of $P$ are adjacent to $a$; let the 
neighbours of $a$ in $P$ be $c_1\ll c_k$ in order on $P$, where $c_1=c$ and $c_k=c'$. For $1\le i<k$, let $P_i$
be the subpath of $P$ between $c_i, c_{i+1}$. For $1\le i<k$, $a\dd c_i\dd P_i\dd c_{i+1}\dd a$ is a hole, and since
$b$ is adjacent to $a$ and not to $c_i, c_{i+1}$, $b$ has an odd number of neighbours in this hole. Hence it has an even
number in $P_i$ for each $i$, and so an even number in $P$ altogether. Let $D$ be the path of $J$ with ends $u,u'$, and choose
an internal vertex $d\in V(D)$. Let $D_1$ be the subpath of $D$ with ends $d,u$, and let $D_2$ be the subpath
with ends $d,u'$. Let $D_3$
be a path of $J$ between $d, v$ where $v\in \beta$. For each edge $g$ of $D_1\cup D_2\cup D_3$, choose a $g$-rung $R_g$,
with $c\in V(R_{e})$ and $c'\in V(R_{e'})$. For $i = 1,2,3$ let $d_i$ be the end of $R_{D_i}$ in $M_d$.
Then 
$$c\dd R_{D_1}\dd d_1\dd d_2\dd R_{D_2}\dd c'\dd P\dd c$$
is a hole, and $b$ has an even number of neighbours in it; so
it has zero, or exactly two adjacent neighbours. Zero is impossible since then \ref{treestructplus1} and \ref{treestructplus2} would imply that
$X(P)$ is local. Thus $b$ has exactly two neighbours $x,y$ in $V(P)$, and they are adjacent. 
Since $x\notin Z$ it follows that $c,x,y,c'$ are all distinct. Let $c,x,y,d$ be in order
in $P$. Then the subgraph induced on $V(R_{D_1}\cup R_{D_2}\cup R_{D_3}\cup P)$ is a prism, 
with bases 
$\{b,x,y\}, \{d_1,d_2,d_3\}$, and constituent paths
$$d_1\dd R_{D_1}\dd c\dd P\dd x,$$
$$d_2\dd R_{D_2}\dd c'\dd P\dd y,$$
$$d_3\dd R_{D_3}\dd b,$$
a contradiction. This proves \ref{treestructplus3}.~\bbox

From \ref{treestructplus1}, \ref{treestructplus2} and \ref{treestructplus3}, this proves \ref{treestructplus}.

\section{Triangles through the cross-edge}

Next we prove some results about the set called $Z$ in \ref{treestructplus}. We need the following lemma.

\begin{thm}\label{skewpyr}
Let $G$ be even-hole-free, and let $H$ be a hole of $G$, with vertices $h_1\dd h_2\cc h_n\dd h_1$ in order.
Let $a,b\in V(G)\setminus V(H)$ each have at least three neighbours in $V(H)$, and let $\{a,b\}$ be complete to $\{h_1,h_n\}$. If $a,b$
are nonadjacent, then one of $a,b$ is adjacent to $h_{n-1},h_n,h_1$ and to no other vertices in $V(H)$, and the other is adjacent
to $h_n,h_1,h_2$ and to no other vertices in $V(H)$.
\end{thm}
\Proof
Let $P$ be the path $h_2\dd h_3\cc h_{n-1}$, and let $A,B$ be the sets of neighbours of $a,b$
respectively in $V(P)$. Since $G$ has no 4-hole, it follows that $A\cap B=\emptyset$. An {\em $(A,B)$-gap} means a subpath of $P$
with one end in $A$ and the other in $B$, and with no internal vertices in $A\cup B$. 
If there is an $(A,B)$-gap containing both $h_{n-1},h_2$ then the theorem holds, and so 
we may assume not; and hence every $(A,B)$-gap is anticomplete to one of $h_n,h_1$, and therefore has odd length
(because it can be completed to a hole by adding $a,b$ and one of $h_n,h_1$). It follows that no two $(A,B)$-gaps are
anticomplete; because their union with $\{a,b\}$ would induce an even hole. 

There is an $(A,B)$-gap, since
$a,b$ each have at least three neighbours in $V(H)$. 
Choose an $(A,B)$-gap $h_i\cc h_j$ with $i<j$ and $i$ minimum, 
and we may assume that
$h_i\in A$. Hence $b$ is nonadjacent to $h_2\ll h_{j-1}$, and so $b\dd h_1\cc h_j\dd b$ is a hole, and therefore $j$ is even. Moreover, $j-i$
is odd, since $h_i\cc h_j$ is an $(A,B)$-gap; and since $n$ is odd, it follows that $n-i=n+ (j-i)-j$ is even. Consequently $a\dd h_i\cc h_n\dd a$
is not a hole, and so there exists $k\in \{j+1\ll n-1\}$ minimum such that $h_k\in A$.
If $B\cap \{h_i\ll h_k\}=\{h_j\}$, there is a theta with ends $b,h_j$ induced on $\{a,b,h_i\ll h_k\}$,
a contradiction.
Thus one of $h_{j+1}\ll h_{k-1}$ is in $B$, and since no two
$(A,B)$-gaps are anticomplete, it follows that $h_{j+1}\in B$ and $h_{j+2}\ll h_k\notin B$. Since no two $(A,B)$-gaps
are anticomplete, $b$ has no more neighbours in $V(P)$; but then it is the centre of a 4-wheel with hole $H$, a contradiction.
This proves \ref{skewpyr}.~\bbox

Let $G$ be even-hole-free, let $ab\in E(G)$,
and let $(J,M)$ be optimal for $ab$. Let $Z$ be the set of common
neighbours of $a,b$ in $G$.
It would be helpful if $Z$ were a clique, but unfortunately this is not true, even assuming that 
$a$ is splendid. It {\em is} true if both $a,b$ are splendid, but that assumption is too strong for our application 
(to find a bisimplicial vertex, later).
But here is something on those lines, good enough for the application and true without any additional hypothesis.
Let us say that a vertex $y\in Z$ is {\em $a$-external} if 
there is a path from $y$ to $V(M)\setminus N[a]$ containing no neighbours of $a$ except $y$, and we define {\em $b$-external}
similarly. Let us say a vertex $y$ is 
{\em major} if $y\in Z$, and $y$ is both $a$-external and $b$-external. For convenience we write $N[a,b]$ for $N[\{a,b\}]$.

\begin{thm}\label{clique}
Let $ab$ be an edge of an even-hole-free graph $G$,
and let $(J,M)$ be optimal for $ab$. 
Then the set of all major vertices is a clique.
\end{thm}
\Proof
Let $Z$ be the set of common
neighbours of $a,b$ in $G$, and $Y$ the set of major vertices (thus $Y\subseteq Z$).
\\
\\
(1) {\em If $y\in Y$, then either $y$ is complete to one of $N(a)\cap V(M), N(b)\cap V(M)$, or
there is a path from $y$ to $V(M)\setminus N[a,b]$ containing no neighbours
of $a$ or $b$ except $y$.}
\\
\\
We may assume that $X(y)\subseteq N[a,b]$, for otherwise a path of length one satisfies the claim.
Since $y$ is $b$-external, there is a minimal path $P$ with one end $y$, containing no neighbour of $b$ except $y$, 
such that its other end ($p$ say) has a neighbour 
in $V(M)\setminus N[b]$. It follows that $V(P)\cap V(M)=\emptyset$.
Similarly, there is a minimal path $Q$ between $y$ and $q$ say, containing no neighbour of $a$ 
except $y$, where $X(q)\not\subseteq N[a]$.
Thus $a$ might have neighbours in $V(P\setminus y)$, and $b$ might have neighbours in $V(Q\setminus y)$.

If $X(p)\not\subseteq N[a,b]$, then by \ref{treestructplus}, $a$ has no neighbour
in $V(P\setminus y)$ and the claim holds. Thus we may assume that $X(p)\subseteq N[a,b]$, and $X(p)\not\subseteq N[b]$
from the definition of $P$. We claim that $X(p)\subseteq N[a]$. Suppose not; then $p$ has a neighbour in $V(M)\cap N[a]$
and one in $V(M)\cap N[b]$. By \ref{treestructplus}, $p$ is adjacent to both $a,b$, and so $p=y$. Choose
$t_1\in \alpha$ such that $X(p)\cap M_{t_1}\ne \emptyset$, and $t_2\in \beta$ such that 
$X(p)\cap M_{t_2}\ne \emptyset$, and let $D$ be a path of $J$ with ends $t_1,t_2$. For $i = 1,2$ let $e_i\in E(J)$ be incident with $t_i$.
For each $e\in E(D)$,
choose an $e$-rung $R_e$, such that $R_{e_1}$ contains a vertex in $X(p)\cap N[a]$ and $R_{e_2}$ 
contains a vertex in $X(p)\cap N[b]$. Then $p$ has exactly four neighbours in the hole $a\dd R_D\dd b\dd a$, since
$X(y)\subseteq N[a]\cup N[b]$, and so $G$ contains a
4-wheel, a contradiction. This proves that $X(p)\subseteq N[a]$. Similarly $X(q)\subseteq N[b]$.

Let $t_1\in \alpha$ and $t_2\in \beta$, such that $X(p)\cap M_{t_1}\ne\emptyset$,
and $X(q)\cap M_{t_2}\ne \emptyset$. For $i = 1,2$ let $e_i\in E(J)$ be $J$-incident with $t_i$.
Choose $v_1\in X(p)\cap M_{t_1}$ and $v_2\in X(q)\cap M_{t_2}$. Let $D$
be a path of $J$ with ends $t_1,t_2$, and for each $e\in E(D)$ let $R_e$ be an $e$-rung, with $v_i\in V(R_{e_i})$
for $i = 1,2$. Then $R_D$ is an induced path with ends $v_1,v_2$, and with interior anticomplete to $a,b$ and to $V(P\cup Q)$.

By \ref{treestructplus}, there is no path between $v_1,v_2$, with interior disjoint from $V(M)\cup Z$, and so
$V(P\setminus y)\cup \{v_1\}$ is disjoint from and anticomplete to $V(Q\setminus y)\cup \{v_2\}$.
Consequently $v_1\dd P\dd y\dd Q\dd v_2$ is an induced path. Now as we saw above, $p\ne q$ and so at least one of $P,Q$
has length at least one, say $Q$. Thus $b$ has two nonadjacent neighbours in the hole 
$$v_2\dd q\dd Q\dd y\dd P\dd p\dd v_1\dd R_D\dd v_2,$$
and so has an odd number, at least three. They all belong to the path $v_2\dd q\dd Q\dd y$. We may assume that $y$
is not complete to $V(M)\cap N[a]$, so there exists $e_3=s_3t_3$ where $t_3\in \alpha$ and an $e_3$-rung $R_{e_3}$ 
such that $y$ has no neighbour
in $V(R_{e_3})$ (because $X(y)\subseteq N[a,b]$). 
Let $D$
be a path of $J$ with ends $t_2, t_3$, and for each $e\in E(D)$ let $R_e$ be an $e$-rung, with $v_i\in V(R_{e_i})$
for $i = 2,3$. Then the hole
$$v_2\dd q\dd Q\dd y\dd a\dd v_3\dd R_D\dd v_2$$ contains exactly one neighbour of $b$ in addition to those in 
$v_2\dd q\dd Q\dd y$, and so contains an even number, a contradiction. This proves (1).

\bigskip
For each $y\in Y$, let $P_y$ be some minimal path of $G$ between $y$ and its other end (say $p_y$) such that 
$a,b$ have no neighbours in $V(P_y\setminus y)$ and $X(p_y)\not \subseteq N[a,b]$, if there is such a path. If not,
let $P_y$ be the one-vertex path with vertex $y$, and let $p_y=y$.
From the minimality
of $P_y$, $X(P_y\setminus p_y)\subseteq N[a,b]$. (Note that there are two cases when $p_y=y$, the two extremes: when we don't need the path 
$P_y$, because $y$ itself has a neighbour
in $V(M)\setminus N[a,b]$; and when we can't find the path $P_y$, and therefore $y$ is complete to one of $N[a]\cap V(M), N[b]\cap V(M)$
by (1).)
\\
\\
(2) {\em Let $t_1\in \alpha$ and $t_2\in \beta$, and let $D$
be a path of $J$ with ends $t_1,t_2$. Let $y\in Y$. If 
$$X(P_y)\setminus N[a,b]\not\subseteq \bigcup_{e\in E(D)}M_e,$$
there is a vertex $d$ of $D$ and a path $Q$ of $G$ with the following properties:
\begin{itemize}
\item 
$d$ is an internal vertex of $D$, incident with edges $g_1,g_2$ of $D$ say;
\item $Q$ has ends $y,d_3$, where $d_3\in M_d\setminus (M_{g_1}\cup M_{g_2})$; 
\item $V(Q)\subseteq V(P_y)\cup V(M)$; and
\item $Q^*$ is anticomplete to $\bigcup_{e\in E(D)}M_e$, and $V(Q\setminus y)$ is anticomplete to $\{a,b\}$.
\end{itemize}
}
\noindent Since $X(P_y)\setminus N[a,b]\not\subseteq \bigcup_{e\in E(D)}M_e$, 
there exists $e_3\in E(J)\setminus E(D)$ such that $X(P_y)\setminus N[a,b]$ meets $M_{e_3}$. 
Let $C$ be a path of $J$, containing $e_3$ and edge-disjoint from $D$ and with one end in $V(D)$; and choose $e_3, C$
with $C$ minimal. 
Let $d$ be the end of $C$ in $V(D)$.
Choose an $e$-rung $R_e$ for each $e\in E(C)$, choosing $R_{e_3}$ to contain a vertex of $X(P_y)\setminus N[a,b]$. 
Then $R_C$ is an induced path containing a vertex
in $X(p_y)\setminus N[a,b]$, with ends $c,d_3$ say, and $d_3\in M_d\setminus (M_{g_1}\cup M_{g_2})$, where
$g_1,g_2$ are the two edges of $D$ incident with $d$.
Thus $R_C\setminus d_3$ is anticomplete to $\bigcup_{e\in E(D)}M_d$. 
Moreover no vertex of $R_C$ belongs to $N[a,b]$ except
possibly $c$, and in that case $p_y$ has a neighbour in $R_C$ different from $c$.
Choose a minimal subpath $S$ of $R_C$ that has one end $d_3$ and the other
adjacent to $p_y$. Then no vertex of $S$ is adjacent to $a$ or $b$, and so setting $Q$ to be the path
$y\dd P_y\dd p_y\dd S\dd d_3$ satisfies the claim. This proves (2).
\\
\\
(3) {\em Let $t_1\in \alpha$ and $t_2\in \beta$, and let $D$
be a path of $J$ with ends $t_1,t_2$. 
For each $e\in E(D)$ let $R_e$ be an $e$-rung.
For each $y\in Y$, either 
$X(P_y)\cap V(R_D)\ne \emptyset$, or
$$\emptyset\ne X(P_y)\setminus N[a,b]\subseteq \bigcup_{d\in E(D)}M_d.$$
In either case, $X(P_y)\cap \bigcup_{d\in E(D)}M_d$ is nonempty.}
\\
\\
If $X(P_y)\cap V(R_D)\ne \emptyset$ then the claim holds, so we may assume that $X(P_y)\cap V(R_D)= \emptyset$.
Consequently 
$y$ is not complete to either of $N[a]\cap V(M), N[b]\cap V(M)$, and so by (1), $X(P_y)\not\subseteq N[a,b]$. 
Suppose that 
$$X(P_y)\setminus N[a,b]\not \subseteq \bigcup_{e\in E(D)}M_e.$$
Choose $d,Q$ as in (2), and for $i = 1,2$ let $D_i$ be the subpath of $D$ between $d$ and $t_i$.
Let $Q$ have ends $y,d_3$.
Thus $d_3$ has two adjacent neighbours $d_1,d_2$ in $R_D$,
where $d_i\in R_{D_i}$ for $i = 1,2$. 
But then there is a near-prism with bases
$\{d_1,d_2,d_3\}$ and $\{a,b,y\}$, with constituent paths 
$$a\dd R_{D_1}\dd d_1,$$
$$b\dd R_{D_2}\dd d_2,$$
$$y\dd Q\dd d_3,$$
a contradiction.  This proves (3).

\bigskip

Choose distinct $a_1,a_2\in \alpha$ and $b_1,b_2\in \beta$ such that the paths $D_1,D_2$ are vertex-disjoint, where 
for $i = 1,2$, $D_i$ is the path of $J$ with ends $a_i, b_i$. (This is possible since $J$ has at least
two vertices that are not leaves, by hypothesis.)  For $i = 1,2$, let $W_i=\bigcup_{e\in E(D_i)}M_e$. 
We observe that $W_i$ is connected, because $D_i$ has an internal vertex $d$, and $M_d\cap W_i$ is connected,
and every other vertex of $W_i$ can be joined to $M_d\cap W_i$  by a union of rungs.

Suppose that $y_1,y_2\in Y$ are nonadjacent. 
For $i = 1,2$, let us say an induced path $T$ of $G$ between $y_1,y_2$ is {\em $i$-normal} if
for every $v\in V(T^*)\setminus W_i$, there exists $j\in \{1,2\}$ such that
$v\in V(P_{y_j}\setminus y_j)$ and $X(P_{y_j}\setminus y_j)\cap W_i$ is nonempty.
\\
\\
(4) {\em For $i = 1,2$, there is an $i$-normal path.}
\\
\\
Let $i\in \{1,2\}$. For each $j\in \{1,2\}$, (3) implies that $X(P_{y_j})\cap W_i\ne \emptyset$; and so either 
$y_j$ has a neighbour in $W_i$, or 
$X(P_{y_j}\setminus y_j)\cap W_i\ne \emptyset$. Hence there is a path $S_j$
between $y_j$ and a vertex of $W_i$, such that for every $v\in V(S_j)$, either
\begin{itemize}
\item $v\in \{y_j\}\cup W_i$; or
\item $v\in X(P_{y_j}\setminus y_j)$, and $X(P_{y_j}\setminus y_j)\cap W_i\ne \emptyset$. 
\end{itemize}
Since $W_i$ is connected,
it follows that there is an induced path joining $y_1,y_2$ with interior in $V(S_1\cup S_2)\cup W_i$; and this is therefore $i$-normal.
This proves (4).
\\
\\
(5) {\em For $i = 1,2$, let $T_i$ be an $i$-normal path. Then $T_1^*$ is anticomplete to $T_2^*$.}
\\
\\
Suppose not. Since $W_1$ is anticomplete to $W_2$, we may assume (exchanging $D_1,D_2$ or $y_1,y_2$ if necessary)
that there exist $v_1\in V(P_{y_1}\setminus y_1)\cap T_1^*$, and $v_2\in T_2^*$, such that $v_1,v_2$ are equal or adjacent.
Hence 
$X(P_{y_1}\setminus y_1)\cap W_1\ne \emptyset$. By \ref{treestructplus}, $X(P_{y_1}\setminus y_1)$ is local,
and consequently $X(P_{y_1}\setminus y_1)$ is disjoint from
$W_2$; and in particular $v_2\notin W_2$,
and so $v_2\in V(P_{y_2}\setminus y_2)\cap T_2^*$. By the same argument with $T_1, T_2$
exchanged, $X(P_{y_2}\setminus y_2)$ meets $W_2$. But $Q=(P_{y_1}\setminus y_1)\cup (P_{y_2}\setminus y_2)$ is connected
and $X(Q)$ meets both $W_1$ and $W_2$, and so is not local, contrary to \ref{treestructplus}. This proves (5).

\bigskip
From (4), for $i = 1,2$ there is an $i$-normal path $T_i$. By (5), $T_1\cup T_2$
is a hole, and so one of $T_1, T_2$ is odd and one is even; say $T_1$ is odd and $T_2$ is even. For every $1$-normal path $T_1'$,
$T_1'\cup T_2$ is a hole, and so $T_1'$ is odd, and similarly every $2$-normal path is even.
\\
\\
(6) {\em Every $2$-normal path meets both $M_{a_2}, M_{b_2}$. In particular, if $X(P_{y_2}\setminus y_2)$ meets $W_2$ then $y_1$
has no neighbour in $V(P_{y_2})$, and if $X(P_{y_1}\setminus y_1)$ meets $W_2$ then $y_2$
has no neighbour in $V(P_{y_1})$.}
\\
\\
Let $T_2$ be 2-normal.
Since $T_2$ is even, 
$y_1\dd T_2\dd y_2\dd a\dd y_1$ is not a hole; and so $a$ has a neighbour in $T_2^*$, and similarly so does $b$.
But $a$ has no neighbours in $P_{y_1},P_{y_2}$ different from $y_1,y_2$, and the set of neighbours of $a$ in $W_2$
is $M_{a_2}$. Hence $T_2^*$ meets $M_{a_2}$, and similarly it meets $M_{b_2}$. This proves the first claim.
For the second, observe that if $X(P_{y_2}\setminus y_2)$ meets $W_2$ and $y_1$ has a neighbour in $V(P_{y_2})$ then there is
a 2-normal path with interior in $V(P_{y_2})$ and therefore not meeting both (or indeed, either of) $M_{a_2}, M_{b_2}$,
a contradiction. This proves (6).

\bigskip

For each $e\in E(D_1)$ choose an $e$-rung $R_e$.
\\
\\
(7) {\em One of $y_1,y_2$ has no neighbour in $R_{D_1}$.}
\\
\\
Suppose that $y_1,y_2$ both have a neighbour in $V(R_{D_1})$. By \ref{skewpyr}, since $y_1,y_2$ are nonadjacent, 
it follows that one of $y_1,y_2$ is adjacent to $a_1$, and the other to $b_1$, and neither has any more neighbours in 
$V(R_{D_1})$. Since $R_{D_1}$ is even, adding $y_1,y_2$ to $R_{D_1}$ gives a $1$-normal path of even length,
a contradiction. This proves (7).

\bigskip
Henceforth we assume that $y_1$ has no neighbour in $R_{D_1}$. 
\\
\\
(8) {\em $X(P_{y_1}\setminus y_1)\cap W_2=\emptyset$; $X(y_1)\cap W_2\subseteq N[a,b]$;
and $X(P_{y_1})\setminus N[a,b]\not\subseteq W_2$.}
\\
\\
Since $y_1$ has no neighbour in $R_{D_1}$, it follows that $y_1$ is not complete to either of $N[a]\cap V(M), N[b]\cap V(M)$, 
and so by (1), $X(P_{y_1})\not\subseteq N[a,b]$. 
Suppose that $X(P_{y_1}\setminus y_1)\cap W_2\ne \emptyset$.
Consequently $P_{y_1}$ has length at least one, and so $X(y_1)\subseteq N[a,b]$.
Moreover, $X(P_{y_1}\setminus y_1)\cap W_1=\emptyset$, since $X(P_{y_1}\setminus y_1)$ is local.
Since $y_1$ has no neighbour in $R_{D_1}$, it follows that $X(P_{y_1})\cap V(R_{D_1})=\emptyset$. By (3), 
$X(P_{y_1})\setminus N[a,b]\subseteq W_1$. But $X(P_{y_1})\setminus N[a,b]\subseteq X(P_{y_1}\setminus y_1)$,
since $X(y_1)\subseteq N[a,b]$; and so
$$X(P_{y_1})\setminus N[a,b]\subseteq X(P_{y_1}\setminus y_1)\cap W_1=\emptyset,$$
a contradiction. This proves the first claim. For the second, suppose that $y_1$ has a neighbour in $W_2\setminus N[a,b]$.
From the minimality of $P_{y_1}$, $p_{y_1}=y_1$. Consequently $X(P_{y_1})\cap V(R_{D_1})=\emptyset$, and so by (3), 
$X(P_{y_1})\setminus N[a,b]\subseteq W_1$, contradicting that $y_1$ has a neighbour in $W_2\setminus N[a,b]$. This 
proves the second claim. The third claim follows, since we have shown that $X(P_{y_1})\setminus N[a,b]\ne \emptyset$
and is disjoint from $W_2$. This
proves (8).

\bigskip

For each $e\in E(D_2)$, choose an $e$-rung $R_e$, such that $X(P_{y_2})$ meets $R_{D_2}$ (this is possible by (3)).
Since $X(P_{y_1})\setminus N[a,b]\not\subseteq W_2$ by (8), it follows from (3) that $X(P_{y_1})$ meets $R_{D_2}$;
and since $X(P_{y_1}\setminus y_1)\cap W_2=\emptyset$ by (8), it follows that $y_1$ has a neighbour in $R_{D_2}$.
Thus there is a 2-normal path $T_2$ meeting $W_2$ in a subpath of $R_{D_2}$. By (6), both ends of $R_{D_2}$
belong to $T_2$. Consequently a unique vertex of $R_{D_2}$, one of its ends, is adjacent to $y_1$, and a unique vertex of
$R_{D_2}$, its other end, belongs to $X(P_{y_2})$. Let $R_{D_2}$ have ends $s,t$ where $s\in M_{a_2}$
and $t\in M_{b_2}$. Exchanging $a,b$ if necessary, we may assume that $y_1$ is adjacent to $s$ and to no 
other vertex of $R_{D_2}$, and $X(P_{y_2})\cap V(R_{D_2})=\{t\}$. Choose a minimal subpath $P_2$ of $P_{y_2}$
with ends $y_2, p_2$ say, such that $p_2$ is adjacent to $t$. (Possibly $p_2=y_2$.)

By (8), $X(P_{y_1})\setminus N[a,b]\not\subseteq W_2$. By (2) 
there is a vertex $d$ of $D_2$ and a path $Q$ of $G$ with the following properties:
\begin{itemize}
\item
$d$ is an internal vertex of $D_2$, incident with edges $g_1,g_2$ of $D_2$ say;
\item $Q$ has ends $y_1,d_3$, where $d_3\in M_d\setminus (M_{g_1}\cup M_{g_2})$; 
\item $V(Q)\subseteq V(P_{y_1})\cup V(M)$; and
\item $Q^*$ is anticomplete to $W_2$, and $V(Q\setminus y_1)$ is anticomplete to $\{a,b\}$.
\end{itemize}
In particular, $d_3$ has exactly two neighbours in $V(R_{D_2})$, say $d_1,d_2$ where $s,d_1,d_2,t$ are in order in $R_{D_2}$, 
and $d_1,d_2$ are adjacent.
\\
\\
(9) {\em $y_2$ is nonadjacent to $t$, and $V(P_2\setminus y_2)$ is not anticomplete to $V(Q)$, 
and $y_1$ has no neighbour in $V(P_{y_2})$.}
\\
\\
Suppose first that $V(P_2)$ is anticomplete to $V(Q)$. Then there is a near-prism with bases $\{d_1,d_2,d_3\},\{a,s,y_1\}$
and constituent paths 
$$y_1\dd Q\dd d_3,$$
$$a\dd y_2\dd P_2\dd p_2\dd t\dd R_{D_2}\dd d_2,$$
$$s\dd R_{D_2}\dd d_1,$$
contrary to \ref{subgraphs}.

Thus $V(P_2)$ is not anticomplete to $V(Q)$. Suppose next that $V(P_2\setminus y_2)$ is anticomplete to $V(Q)$, and therefore
$y_2$ has a neighbour in $V(Q)$.
Let $Q'$ be a path with ends $y_2, d_3$, where $Q'\setminus y_2$ is a subpath of $Q$. It follows that $y_1$ has
no neighbour in $V(Q')$; for $y_1$ only has one neighbour in $V(Q)$, and that vertex is not adjacent to $y_2$ since otherwise
there would be a 4-hole. If $y_2$ is adjacent to $t$, then there is a near-prism with bases $\{y_2,b,t\}, \{d_1,d_2,d_3\}$
and constituent paths
$$y_2\dd Q'\dd d_3,$$
$$b\dd y_1\dd s\dd R_{D_2}\dd d_1,$$
$$t\dd R_{D_2}\dd d_2.$$
If $y_2$ is not adjacent to $t$, there is a theta in $G$ with ends
$y_2,t$ and constituent paths
$$y_2\dd P_2\dd p_2\dd t,$$
$$y_2\dd b\dd t,$$
$$y_2\dd Q'\dd d_3\dd d_2\dd R_{D_2}\dd t,$$
contrary to \ref{subgraphs}. This proves that $P_2\setminus y_2$ is not anticomplete to $V(Q)$. 
In particular, $P_2$ has length at least one, and so $y_2,t$ are nonadjacent. Hence $t\in X(P_2\setminus y_2)$,
and so by (6), $y_1$ has no neighbour in $V(P_{y_2})$.  This proves (9).

\bigskip

By (9), we may choose $v_1\in V(Q)$ and $v_2\in V(P_2\setminus y_2)$, such that $v_1,v_2$ are equal or adjacent.
Since $y_1$ has no neighbour in $V(P_2)$ by (9), it follows that either $v_1\in V(P_{y_1}\setminus y_1)$ or $v_1\in V(M)$.
Suppose that $v_1\in V(P_{y_1}\setminus y_1)$. Then $F=G[V(P_{y_1}\setminus y_1)\cup V(P_{y_2}\setminus y_2)]$ 
is connected, and disjoint from $V(M)\cup Z$, and $X(F)$ includes both $X(P_{y_1}\setminus y_1)$ and $X(P_{y_2}\setminus y_2)$.
But since $P_{y_1}$ has length at least one, it follows that $X(P_{y_1}\setminus y_1)\setminus N[a,b]$ is nonempty,
and is a subset of $W_1$ by (3). Hence $X(F)$ meets $W_1$, and contains $t$,
and so $X(F)$ is not local, contrary to
\ref{treestructplus}. 

Thus $v_1\in V(M)$. Since $V(Q)$ is anticomplete to $\{a,b\}$, it follows that 
$v_1\in V(M)\setminus N[a,b]$. From the minimality of $P_{y_2}$, no vertex of $P_{y_2}$ except $y_2$ has a neighbour in
$V(M)\setminus N[a,b]$, and so $v_2=p_{y_2}$, and in particular $P_2=P_{y_2}$. 
But 
$X(P_{y_2}\setminus y_2)$ is local, and contains $t$ and $v_1$. Since $t\in M_{t_2}$, and $Q^*$ is anticomplete to $W_2$, 
it follows that 
$v_1,t\in M_{s_2}$, and hence 
$v_1=d_3$ and $t=d_2$. Morover, $V(Q)$ is disjoint from $V(P_2\setminus y_2)$, and the edge $p_{y_2}\dd d_3$ is the only 
edge joining them. (But $y_2$ might have neighbours in $V(Q)$.) Now $y_2$ is nonadjacent to $d_3$, since
otherwise $y_2\dd d_3\dd d_2\dd b\dd y_2$ is a 4-hole. Then 
$$b\dd y_1\dd s\dd R_{D_2}\dd d_1\dd d_3\dd p_{y_2}\dd P_{y_2}\dd y_2\dd b$$
is a hole, and $d_2=t$ has exactly four neighbours in it, namley $d_1,d_3,p_{y_2}$ and $b$, a contradiction.
This proves \ref{clique}.~\bbox

Finally, we have:
\begin{thm}\label{funnies}
Let $ab$ be an edge of an even-hole-free graph $G$,
and let $(J,M)$ be optimal for $ab$.
Let $Z$ be the set of all common neighbours of $a,b$, and let $Y\subseteq Z$ be the set of all major vertices.
If $F$ is a component of $G\setminus (V(M)\cup Z)$, and some vertex in $Z\setminus Y$ has a neighbour in $F$, then
there is a leaf $t$ of $V(J)$, such that 
every vertex in $V(M)$ with a neighbour in $V(F)$ belongs to $M_t$.
\end{thm}
\Proof
Let $z\in Z\setminus Y$ have a neighbour in $V(F)$. Let $X(F)$ be the set of vertices in $V(M)$ with a neighbour in $V(F)$.
If one of $a,b$ has a neighbour in $V(F)$, the claim follows from \ref{treestructplus}, so suppose not. 
If $X(F)\not\subseteq N[a,b]$ has a neighbour in $V(F)$, this contradicts that $z$ is not major.
So $X(F)\subseteq N[a,b]$, and then the claim follows since
$X(F)$ is local  by \ref{treestructplus}.
This proves \ref{funnies}.~\bbox

Let us summarize the previous results. The vertices of $G$ are partitioned into the following sets:
\begin{itemize}
\item The special vertices $a,b$.
\item $V(M)$ (this is further partitioned into strips corresponding to the edges of $J$).
\item The small components. Each small component $F$ satisfies $X(F)\subseteq M_e$ for some $e\in E(J)$ or $X(F)\subseteq M(t)$
for some $t\in V(J)$. Moreover if $N(F)$ contains $a$ or $b$, or a vertex in $Z\setminus Y$,
then $F$ must be peripheral, and if $N(F)$ contains only one of $a,b$,
then $X(F)\subseteq N[a]$ or $N[b]$ correspondingly.
\item The set $Y$ of the major vertices. These form a clique, but we know nothing about their neighbours outside of $Z$.
\item The vertices in $Z\setminus Y$. All their neighbours in $V(M)$ are in $N[a,b]$, and all their neighbours
in small components belong to peripheral small components.
\end{itemize}

If we assume that $a$ is splendid (which will be true in our application), we can simplify the theorem a little; let us
see that next. We need:
\begin{thm}
\label{nosmall}
Let $ab$ be an edge of an even-hole-free graph $G$,
and let $(J,M)$ be optimal for $ab$.
If $a$ is splendid,
there is no small $F$ such that $a$ has a neighbour in $V(F)$.
\end{thm}
\Proof
Let
$Z$ be the set of vertices of $G$ adjacent to both $a,b$.
Suppose that there is such an subgraph $F$, and we may assume that $F$ is small component.
If $b$ has no neighbour in $V(F)$, then since by \ref{treestructplus}
every vertex in $V(M)$ with a neighbour in $V(F)$ belongs to $N[a]$, it follows that
$F$ is a component of $G\setminus N[a]$, contradicting that $a$ is splendid. Thus $b$ has a neighbour in $V(F)$.
For the same reason, some vertex of $V(M)$ nonadjacent to $a$ has a neighbour in $V(F)$; but by \ref{treestructplus},
every such vertex belongs to $B$.

Hence there is an induced path $P$ of $F$ such that $a$ has a neighbour in $V(P)$,
and some vertex in $B$
has a neighbour in $V(P)$.
Let $P$ be minimal with this property.
Let $P$ have ends $p_1,p_2$, where $a$ is adjacent to $p_1$ and to no other vertex of $V(P)$, and some vertex in $B$
($v_2$ say) is adjacent to $p_2$, and no vertex in $B$ has a neighbour in $V(P\setminus p_2)$.
Since $p_1$ is nonadjacent to $b$ (because $p_1\notin Z$) and there is no 4-hole, it follows that
$p_1$ is anticomplete to $B$, and in particular
$p_1\ne p_2$. Let $v_2\in M_{e_2}$, where $e_2\in \beta$.
From \ref{treestructplus}, there is at most one $e\in \alpha$ such that $M_e$ is not anticomplete to $V(P\setminus p_2)$,
and so there exists $d_1\in \alpha$ such that $M_{d_1}$ is anticomplete to $V(P\setminus p_2)$. Since $M_{d_1}$
is anticomplete to $p_2$ by \ref{treestructplus}, it follows that $M_{d_1}$
is anticomplete to $V(P)$.

There is a path $D$ of $J$ with end-edges $d_1,e_2$. Let $R_e$ be an $e$-rung for each $e\in E(J)$, with $v_2\in V(R_{e_2})$;
then $R_D$
is an induced path of $G$ between $a, v_2$, with interior in $V(M)$ and anticomplete to $V(P)$.
Hence $P\cup R_D$ is a hole, and $b$ has two nonadjacent neighbour in $P\cup R_D$, namely $v_2,a$; and since
$G$ has no full star cutset, \ref{bigvertex} applied to $b$ and $P\cup R_D$ implies that $b$ is adjacent to $p_2$
and has no other neighbour in $V(P)$. But then there is a short pyramid with apex $a$ and base $\{b,v_2,p_2\}$,
and constituent paths
$$a\dd b,$$
$$a\dd f_1\dd F\dd f_2,$$
$$a\dd R_D\dd v_2,$$
contradicting that $a$ is splendid. This proves \ref{nosmall}.~\bbox

We deduce an upgraded version of \ref{treestructplus}:
\begin{thm}\label{splendidprism}
Let $G$ be even-hole-free, and $ab$ be an edge of $G$, where $a$ is splendid. Let $(J,M)$ be optimal for $ab$.
Let $Z$ be the set of vertices of $G$ adjacent to both $a,b$, and let $Y$ be the set of major vertices.
Then
\begin{itemize}
\item every vertex in $V(M)$ with a neighbour in $Z\setminus Y$ belongs to $M_t$ for some $t\in \beta$;
and
\item for each $e=st\in \alpha$, $M_s\cap M_t=\emptyset$.
\end{itemize}
Moreover,
for every small subgraph $F$, let $X$ be the set of vertices in $V(M)$
with a neighbour in $V(F)$; then
\begin{itemize}
\item $a$ has no neighbours in $V(F)$; 
\item if $V(F)$ is anticomplete to $\{b\}\cup (Z\setminus Y)$, then either $X\subseteq M_e$ for some $e\in E(J)$ or
$X\subseteq M_t$ for some $t\in V(J)\setminus \alpha$;
\item if either $b$ or some vertex in $Z\setminus Y$ has a neighbour in $V(F)$, then $X\subseteq M_t$ for some $t\in \beta$.
\end{itemize}
\end{thm}
\Proof
Since $a$ is splendid, every vertex in $Z$ is $a$-external, and therefore the vertices in $Z\setminus Y$ are not $b$-external.
In particular, none of them has a neighbour in $V(M)\setminus N[b]$. That proves the first claim.

Suppose that there exists $e=st\in \alpha$ where $t$ is a leaf of $J$, and
$M_s\cap M_t\ne \emptyset$. Let $v\in M_s\cap M_t$. Let $D$ be a path of $J$ containing $s$, with one end
in $\alpha\setminus \{t\}$ and the other in $\beta$.
Choose an $e$-rung $R_e$ for every $e\in E(D)$.
Then the subgraph of $G$ induced on $V(R_D)\cup \{a,b,v\}$ is a short pyramid with apex $a$,
contradicting that $a$ is splendid. This proves the second claim. The third claim, about small sets, follows from 
\ref{treestructplus} and \ref{nosmall}.
This proves \ref{splendidprism}.~\bbox

\section{Graphs with no extended near-prism}

It would be nice if we had a decomposition theorem complementary to the results of the previous sections, describing a decomposition
for even-hole-free graphs that do not contain a extended near-prism. We do not have that; we only have a decomposition theorem for such graphs
that have a splendid vertex. (This is good enough for our purposes, since it is straightforward to show that every minimum
counterexample to \ref{mainthm} has a splendid vertex.) Our next goal is to state and prove this decomposition theorem.

A {\em pyramid strip system} $\mathcal{S} =(a,S_1\ll S_k)$ in $G$ consists of a set of proper strips $S_1\ll S_k$ with $k\ge 3$, 
pairwise vertex-disjoint (that is, the sets
$V(S_1)\ll V(S_k)$ are pairwise disjoint), and a vertex $a$ of $G$ called the {\em apex}, such that,
setting $S_i=(A_i, B_i, C_i)$ for $1\le i\le k$:
\begin{itemize}
\item for $1\le i< j\le k$, $B_i$ is complete to $B_j$, and there are no other edges between $V(S_i)$ and $V(S_j)$;
\item $a$ belongs to none of $V(S_1)\ll V(S_k)$;
\item for $1\le i\le k$, $a$ is complete to $A_i$, and anticomplete to $B_i\cup C_i$.
\end{itemize}
Let $V(\mathcal{S})$ denote $V(S_1)\cup\cdots\cup V(S_k)\cup \{a\}$. For 
an induced subgraph $F$ of $G$ with $V(F) \subseteq V(G)\setminus V(\mathcal{S})$, we say $v\in V(\mathcal{S})$ is an 
{\em attachment} of $F$ if $v$ has a neighbour in $F$, and we define $\mathcal{S}(F)$ to be the set of all attachments of $F$.
A proper strip $S=(A,B,C)$ is {\em indecomposable} if $A\cup C$ is connected, and a pyramid strip system is {\em indecomposable}
if all its strips are indecomposable. 

If $a\in V(G)$ is the apex of a pyramid, then it is also the apex of an indecomposable  pyramid strip system with $k=3$
and with only one rung in each strip. That motivates the following:

\begin{thm}\label{getpyramidstrips}
Let $G$ be even-hole-free, and let $a\in V(G)$ be splendid. Suppose
there is no extended near-prism contained in $G$ such that $a$ is an end of its cross-edge. Let 
$\mathcal{S}=(a,S_1\ll S_k)$ be 
an indecomposable strip system with apex $a$, with strips $S_i=(A_i,B_i,C_i)$ for $1\le i\le k$,
chosen with $V(\mathcal{S})$ maximal. Then 
for each component $F$ of $G\setminus (V(\mathcal{S})\cup N[a]$, either $\mathcal{S}(F)$ is a nonempty subset of $B_1\cup \cdots\cup B_k$, or for some $i\in \{1\ll k\}$, $\mathcal{S}(F)$ is a subset of one of
$V(S_i)$ and has nonempty intersection with $B_i\cup C_i$.
\end{thm}
\Proof
First we observe:
\\
\\
(1) {\em For each component $F$ of $G\setminus (V(\mathcal{S})\cup N[a])$, $\mathcal{S}(F)$ has nonempty intersection 
with $B_i\cup C_i$ for some $i\in \{1\ll k\}$.}
\\
\\
If not, then $F$ is a component of $G\setminus N[a]$, which is impossible since $G\setminus N[a]$ is connected (because
$a$ is splendid). This proves (1).
\\
\\
(2) {\em For each vertex $f$ of $G\setminus (V(\mathcal{S})\cup N[a])$, $\mathcal{S}(f)$ is either a subset of
$B_1\cupcup B_k$ or a subset of $V(S_i)$ for some $i\in \{1\ll k\}$.}
\\
\\
Suppose not. We may assume $f$ has a neighbour in $A_1\cup C_1$, since $\mathcal{S}(F)$ is not a subset of
$B_1\cupcup B_k$. Choose an $S_1$-rung $R_1$ in which $f$ has a neighbour in $A_1\cup C_1$, with ends $a_1\in A_1$
and $b_1\in B_1$.
Suppose also that $f$ has a neighbour in $A_2\cup C_2$, and choose $R_2, a_2,b_2$ similarly. If $f$ has a neighbour in $V(S_3)$, then
there is a theta with ends $f,a$
and constitutent paths
$$f\dd R_1\dd a_1\dd a,$$
$$f\dd R_2\dd a_2\dd a,$$
$$f_2\dd G[V(S_3)]\dd a,$$
contrary to \ref{subgraphs}. Thus $f$ is anticomplete to $V(S_3)\ll V(S_k)$. If $f$ has two nonadjacent neighbours in
$R_1$, there is a theta with ends $a,f$
and constitutent paths
$$f\dd R_1\dd a_1\dd a,$$
$$f\dd R_2\dd a_2\dd a,$$
$$f\dd R_1\dd b_1\dd G[V(S_3)]\dd a,$$
contrary to \ref{subgraphs}. So $f$ has either one or two adjacent neighbours in $R_1$, and similarly it has one or two adjacent
in $R_2$. Since $f$ is not adjacent to both $a_1,a_2$, we may assume by exchanging $S_1,S_2$ if necessary that
$f$ is not adjacent to $a_1$. If $f$ has a unique neighbour $u$ in $R_1$, 
there is a theta with ends $u,a$
and constitutent paths
$$u\dd R_1\dd a_1\dd a,$$
$$u\dd f\dd R_2\dd a_2\dd a,$$
$$u\dd R_1\dd b_1\dd G[V(S_3)]\dd a,$$
contrary to \ref{subgraphs}. Thus $f$ has exactly two adjacent neighbours in $R_1$, say $p,q$, where $a_1,p,q,b_1$ are in
order in $R_1$. If $f$ also has two adjacent neighbours in $R_2$, there is a 4-wheel with centre $f$ and hole
induced on
$V(R_1\cup R_2)\cup \{a\}$, a contradiction. Thus $f$ has a unique neighbour $u$ in $R_2$. If $u\ne a_2$, we obtain
a contradiction as before; and if $u=a_2$, the subgraph induced on $V(R_1\cup R_2)\cup \{a,f\}$ is an extended near-prism, 
and $a$ is an end of its cross-edge, a contradiction.

This proves that $f$ has no neighbour in $A_2\cup C_2$, and similarly none in $A_i\cup C_i$ for $2\le i\le k$.
If $f$ is complete to $B_2\cupcup B_k$, we can add $f$ to $B_1$, contrary to the maximality of $V(\mathcal{S})$. Thus
$f$ has a neighbour in $B_2\cupcup B_k$, and a non-neighbour in this set. Since $k\ge 3$, we may assume
that $f$ has a neighbour $b_2\in B_2$ and a non-neighbour $b_3\in B_3$. But then there is a theta with ends
$b_2,a$ and constituent paths
$$b_2\dd f\dd G[A_1\cup C_1]\dd a,$$
$$b_2\dd G[A_2\cup C_2]\dd a,$$
$$b_2\dd b_3\dd G[A_3\cup C_3]\dd a,$$
contrary to \ref{subgraphs}. This proves (2).

\bigskip
Let us say a subset $X$ of $V(\mathcal{S})$ is {\em local} if $X$ is a subset of $A_1\cupcup A_k$, or of $B_1\cupcup B_k$
or of $V(S_i)$ for some $i\in \{1\ll k\}$. (Note that in (2) we did not include $A_1\cupcup A_k$, but here we do.)
\\
\\
(3) {\em Every subset of $V(\mathcal{S})$ that is not local includes a 2-element subset 
that is not local.}
\\
\\
Suppose $X\subseteq V(\mathcal{S})$ is not local. If there exists $c\in X\cap C_1$, choose $d\in X\setminus V(S_1)$; then
$\{c,d\}$ is not local. So we may assume that $X\cap C_i=\emptyset$ for $1\le i\le k$. There exists
$c\in X\setminus (A_1\cup\cdots\cup A_k)$, say $c\in B_1$. If there exists $d\in X\cap A_i$ where $i\ge 2$ then $\{c,d\}$ is not local,
so we may assume that $X\cap A_i=\emptyset$ for $2\le i\le k$. Since $X\not\subseteq B_1\cupcup B_k$,
there exists $c\in X\cap A_1$; and since $X\not\subseteq V(S_1)$, there exists $d\in X\cap B_i$ for some $i>1$, and then 
$\{c,d\}$ is not local. (The claim also follows from K\"onig's matching theorem.) This proves (3).

\bigskip

Suppose the theorem is false; then from (1) there is a minimal connected induced subgraph $F$ of 
$G\setminus (V(\mathcal{S})\cup N[a])$ such that $\mathcal{S}(F)$ is not local. By (3) there is a 2-element subset  $\{v_1,v_2\}$
of $\mathcal{S}(F)$ that is not local.
From the minimality of $F$, $F$ is the interior of a path joining $v_1,v_2$.
Let $F$ have ends $f_1,f_2$, where $v_i, f_i$ are adjacent for $i = 1,2$. 
\\
\\
(4) {\em No vertex in $F^*$ has a neighbour in $A_1\cupcup A_k$.}
\\
\\
Suppose that some $f_3\in V(F)\setminus \{f_1,f_2\}$ is adjacent to $a_1\in A_1$ say. Let $F_i$ be the subpath of $F$
between $f_i, f_3$ for $i = 1,2$. From the minimality of $F$, each of $\mathcal{S}(F_1),\mathcal{S}(F_2)$ is a subset
of one of $V(S_1), A_1\cupcup A_k$; and since $\mathcal{S}(F)$ is not local, we may assume that $\mathcal{S}(F_1)\subseteq V(S_1)$
and $\mathcal{S}(F_2)\subseteq A_1\cupcup A_k$. 
Moreover, $v_1\notin A_1\cupcup A_k$ and $v_2\notin V(S_1)$. Thus 
$v_1\in B_1\cup C_1$, and we may assume that $v_2\in A_2$. From the minimality of $F$, $\mathcal{S}(F\setminus f_1)$ is local
and hence is a subset of $A_1\cupcup A_k$, and $\mathcal{S}(F\setminus f_2)$ is a subset of
$V(S_1)$ (because they both contains $a_1$). Thus $\mathcal{S}(F\setminus \{f_1,f_2\})\subseteq A_1$, and 
$\mathcal{S}(f_2)\subseteq A_2$ by (2). 
For $i = 1,2$ let $R_i$ be an $S_i$-rung with ends $a_i\in A_i$
and $b_i\in B_i$, containing $v_i$. Thus $v_2=a_2$, and $v_1\ne a_1$, and $a_2$ is the unique neighbour of $f_2$ in $V(R_2)$.

Let $a_1$ have $t$ neighbours in $V(F\setminus f_1)$; thus $t>0$.
Choose a neighbour $c$ of $f_1$ in $V(R_1)$, such that the subpath of $R_1$ between $b_1,c$ is minimal.
Thus $c\ne a_1$. If $c,a_1$ are nonadjacent we can add the interior of the path $c_1\dd F\dd a_1$ to $C_1$,
contrary to the maximality of $V(\mathcal{S})$. So $c,a_1$ are adjacent, and hence $a_1$ has at least $t+1$ neighbours in the
path $b_1\dd R_1\dd c\dd F\dd a_2$. (It would have $t+2$ if $a_1,f_1$ are adjacent, and $t+1$ otherwise.)
This path can be completed to a hole via $a_2\dd R_2\dd b_2\dd b_1$ or via $a_2\dd a\dd G[V(S_3)]\dd b_1$, and the
number of neighbours of $a_1$ in the second hole is one more than in the first. Since there is no even wheel, it follows that
$t=1$, and $f_3$ is the unique neighbour of $a_1$ in $V(F)$;  but then there is a
theta with ends $f_3,c$ and constituent paths
$$f_3\dd F\dd c,$$
$$f_3\dd a_1\dd c_1,$$
$$f_3\dd F\dd a_2\dd R_2\dd b_2\dd b_1\dd R_1\dd c,$$
contrary to \ref{subgraphs}. This proves (4).
\\
\\
(5) {\em If $f_1$ has a neighbour in $A_1\cup C_1$, then $f_2$ has a neighbour in $A_i\cup C_i$ for some $i\in \{2\ll k\}$.}
\\
\\
Suppose not; then $\mathcal{S}(f_2)$ is a subset of $B_1\cupcup B_k$, and we may assume that $f_2$ has a neighbour in $B_2$. 
By (4), no vertex in $A_2$ has a neighbour in $V(F)$, and so
from the minimality of $F$, $\mathcal{S}(F\setminus \{f_1,f_2\})\subseteq B_1$.
If $f_2$ has a nonneighbour $b_3\in B_3$,
there is a theta with ends $b_2,a$ and constituent paths
$$b_2\dd G[A_2\cup C_2]\dd a,$$
$$b_2\dd F\dd f_1\dd R_1\dd a,$$
$$b_2\dd b_3\dd G[A_3\cup C_3]\dd a,$$
contrary to \ref{subgraphs}. So $f_2$ is complete to $B_3$ and similarly to $B_i$ for $3\le i\le k$; and since $k\ge 3$,
it follows by exchanging $S_2,S_3$ that $f_2$ is complete to $B_2$. But then we can add $f_2$ to $B_1$ and $V(F\setminus f_2)$
to $C_1$, contrary to the maximality of $V(\mathcal{S})$. This proves (5).
\\
\\
(6) {\em For $1\le i\le k$, $f_1$ has no neighbour in $C_i$, and does not have both a neighbour in $A_i$ and one in $B_i$.}
\\
\\
Suppose that $f_1$ has either a neighbour in $C_1$, or a neighbour in $A_1$ and one in $B_1$. In the second case,
if the neighbour of $f_1$ in $A_1$ is nonadjacent to the one in $B_1$, we could add $f_1$ to $C_1$, contrary to the maximality
of $V(\mathcal{S})$. Thus in either case, there is an $S_1$-rung $R_1$, such that $f_1$ has either a neighbour in 
$V(R_1)\cap C_1$, or one in $V(R_1)\cap A_1$ and one in $V(R_1)\cap B_1$. Let $R_1$ have ends $a_1\in A_1$ and $b_1\in B_1$.
If $f_1$ has two nonadjacent neighbours in $R_1$, we can add $f_1$ to $C_1$, again contrary to the maximality
of $V(\mathcal{S})$. Thus $f_1$ has either a unique neighbour or exactly two adjacent neighbours in $R_1$. From the
minimality of $F$, $\mathcal{S}(F\setminus f_2)\subseteq V(S_1)$. 

By (5), we may assume that $f_2$ has a neighbour in $A_2\cup C_2$, and so $\mathcal{S}(f_2)$ is a subset of $V(S_2)$ by (2). 
Hence $\mathcal{S}(F\setminus f_1)\subseteq V(S_2)$ by (4).
Consequently $\mathcal{S}(F\setminus \{f_1,f_2\})=\emptyset$. 
The only edges between          
$V(F)$ and $V(\mathcal{S})$ are the edges between $f_1$ and $V(S_1)$, and the edges between $f_2$ and $V(S_2)$.
Choose an $S_2$-rung
$R_2$ in which $f_2$ has a neighbour in $A_2\cup C_2$, with ends $a_2\in A_2$ and $b_2\in B_2$. If $f_2$
has two nonadjacent neighbours in $V(R_2)$ we can add $f_2$ to $C_2$, a contradiction. Thus $f_2$ has one
or exactly two adjacent neighbours in $R_2$. 
Let $f_i$ have $n_i$ neighbours in $V(R_i)$ for $i = 1,2$; thus $n_i\in \{1,2\}$. 

If $n_1=n_2=2$, there is a prism, so we may assume that either $n_1=1$ or $n_2=1$.
If $n_1=1$, let $c$ be the unique neighbour
of $f_1$ in $V(R_1)$  (necessarily $c\in C_1$),
and let $R_3$ be an $S_3$-rung with ends $a_3\in A_3$ and $b_3\in B_3$.
Then there is a theta with ends $c,a$ and constituent paths
$$c\dd R_1\dd a_1\dd a,$$
$$c\dd f_1\dd F\dd f_2\dd R_2\dd a_2\dd a,$$
$$c\dd R_1\dd b_1\dd b_3\dd R_3\dd a_3\dd a,$$
a contradiction. Thus $n_1=2$, and consequently $n_2=1$.

Let $c$ be the unique neighbour of $f_2$ in $V(R_2)$. 
By the same argument with $S_1,S_2$ exchanged, it follows that $c\notin C_2$, and so $c=a_2$.
Let $R_3$ be an $S_3$-rung; then the subgraph induced on $V(R_1\cup R_2\cup R_3\cup F)\cup \{a\}$ is an extended near-prism, 
a contradiction.
This proves (6).

\bigskip

From (6), no vertex of $F$ has a neighbour in $C_1\cupcup C_k$, and since we may assume that $f_1$ has a neighbour in 
$A_1\cup C_1$, it follows from (6) that $\mathcal{S}(f_1)\subseteq A_1$. Since $\{v_1,v_2\}$ is not local,
it follows that $v_2\in B_2\cupcup B_k$, and we may assume that $f_2$ has a neighbour in $B_2$. By (6),
$\mathcal{S}(f_2)\subseteq B_1\cupcup B_k$, contrary to (5).
This proves \ref{getpyramidstrips}.~\bbox

The reader will observe that much of the generality of strip systems was not used in this proof; we never increased the number
of strips, or changed the sets $A_1\ll A_k$. That will come in the next proof, where
we try to enlarge $V(\mathcal{S})$ by adding vertices from $N[a]\setminus V(\mathcal{S})$. 
The {\em parity} of a path or hole is the parity of its length.

For $1\le i\le k$, let $D_i$ be the union of the vertex sets of all components $F$ of
$G\setminus (V(\mathcal{S})\cup N[a])$ such that $\mathcal{S}(F)\cap (A_i\cup C_i)\ne \emptyset$.
For $v\in N[a]\setminus V(\mathcal{S})$, let us say $v$ has:
\begin{itemize}
\item {\em type $\alpha$} if for each $i\in \{1\ll k\}$, either $v$ has a neighbour in $B_i\cup C_i$ or $v$ is complete to $A_i$;
\item {\em type $\alpha'$} if there exists $i\in \{1\ll k\}$ such that $v$ has a neighbour in $D_i$ and none in $B_i\cup C_i$,
and for all $j\in \{1\ll k\}\setminus \{i\}$, $v$ is complete to $A_j$ and anticomplete to  $B_j\cup C_j\cup D_j$ (we also
call this {\em type $\alpha'_i$}; it is ``almost'' a case of type $\alpha$);
\item {\em type $\beta$} if there exists $i\in \{1\ll k\}$ such that $v$ is anticomplete to $A_i\cup B_i\cup C_i$, 
and for all $j\in \{1\ll k\}\setminus \{i\}$,
$v$ has a neighbour in $B_j\cup C_j$ (we also call this {\em type $\beta_i$}).
\end{itemize}
We also need one other type.
In the usual notation, for $v\in N[a]\setminus V(\mathcal{S})$ and $1\le i\le k$, 
let us say $v$ has {\em type $\gamma$} or {\em type $\gamma_i$}, and
$Q$ is the {\em corresponding private path}, if
\begin{itemize}
\item $Q$ is an induced path with one end $v$ and the other $q$ say,  and $V(Q\setminus v)$
is disjoint from $V(\mathcal{S})\cup N[a]$;
\item $q$ has a neighbour in $B$, and $q$ is either complete or anticomplete to $B\setminus B_i$;
\item $v$ is complete to $A_j$ for all 
$j\in \{1\ll k\}\setminus \{i\}$, and $v$ has no neighbours in $B_j\cup C_j\cup D_j$
for $1\le j\le k$; and
\item all edges between $V(\mathcal{S})$ and $V(Q\setminus v)$ are between $q$ and $B$.
\end{itemize}

We will show:
\begin{thm}\label{growstrips}
Let $G$ be even-hole-free, and let $a\in V(G)$ be splendid. Suppose
there is no extended near-prism contained in $G$ such that $a$ is an end of its cross-edge. Let
$\mathcal{S}=(a,S_1\ll S_k)$ be
an indecomposable strip system with apex $a$, with strips $S_i=(A_i,B_i,C_i)$ for $1\le i\le k$,
chosen with $V(\mathcal{S})$ maximal.
For each $v\in N[a]\setminus V(\mathcal{S})$, $v$ has type $\alpha$, $\alpha'$, $\beta$ or $\gamma$.
\end{thm}
\Proof
Let $D_1\ll D_k$ be defined as before.
We begin with:
\\
\\
(1) {\em The sets $D_1\ll D_k$ are pairwise disjoint, and 
every component of $G[B_i\cup C_i\cup D_i]$
contains a vertex of $B_i$.}
\\
\\
This is immediate from \ref{getpyramidstrips}.

\bigskip

Let $H\subseteq I$ be the set of $i\in \{1\ll k\}$ such that $v$ has a neighbour in $B_i\cup C_i$, 
and $J=\{1\ll k\}\setminus H$. 
Let $I\subseteq \{1\ll k\}$ be the set of $i\in \{1\ll k\}$ such that $v$ has a neighbour in $B_i\cup C_i\cup D_i$. 
(Thus $H\subseteq I$.)
\\
\\
(2) {\em If $I\ne \emptyset$, then either:
\begin{itemize}
\item $v$ is complete to $\bigcup_{j\in J}A_j$ (and so $v$ has type $\alpha$), or 
\item $|I|=1$ and $v$ is complete to $\bigcup_{i\notin I}A_i$ (and so $v$ has type $\alpha$ or $\alpha'$),
or 
\item $|J|= 1$, $J=\{j\}$ say,  and $v$ is anticomplete to $A_j$ (and so $v$ has type $\beta_j$). 
\end{itemize}}
\noindent We may assume that $I\ne \emptyset$. Choose $h\in \{1\ll k\}$ as follows:
\begin{itemize}
\item If $H\ne \emptyset$ choose $h\in H$;
\item If $H=\emptyset$ and either $|I|=1$ or $v$ is complete to $A_1\cupcup A_k$, choose $h\in I$;
\item If $H=\emptyset$ and $|I|>1$ and $v$ is not complete to $A_1\cupcup A_k$, choose $h\in I$ such that
$v$ is not complete to $A_j$ for some $j\ne h$.
\end{itemize}
For notational simplicity let us assume $h=1$.
Suppose first that $v$ is complete to $A_j$ for all $j\in J\setminus \{1\}$. If $1\in H$ then
the claim holds, so we may assume that $1\notin H$ and therefore $H=\emptyset$ from the choice of $h$. 
Also, from the choice of $h$, either $|I|=1$ or $v$ is complete to $A_1\cupcup A_k$, and in both cases the claim holds.

Hence we may assume that there exists $j\in J\setminus \{1\}$ such that $v$ is not complete to $A_j$, 
say $j = 2$. 
Choose an induced path $P$ between $v$ and some $b_1\in B_1$ with interior in $C_1\cup D_1$ (this is possible by (1)). 
Choose $a_2\in A_2$ nonadjacent to $v$, and let $R_2$ be an $S_2$-rung containing $a_2$, and let $b_2$
be its end in $B_2$. Now let $a_2'\in A_2$, and define $R_2', b_2'$ similarly. The $R_2, R_2'$ have the same parity,
and so if $v$ is adjacent to $a_2'$ then the holes
$$v\dd P\dd b_1\dd b_2\dd R_2\dd a_2\dd a\dd v,$$
$$v\dd P\dd b_1\dd b_2'\dd R_2'\dd a_2'\dd v$$
have different parity, a contradiction. Thus $v$ is nonadjacent to $a_2'$ for each $a_2'\in A_2$,
and therefore anticomplete to $A_2$.
If $|J\setminus \{1\}|=1$, then $|J|\le 2$, and hence $H\ne \emptyset$, and so $1\in H$ and $|J|=1$.
But then the claim holds. Thus we may assume that 
$|J\setminus \{1\}|\ge 2$; let $3\in J$ say. 
Let $R_3$ be an $S_3$-rung with ends $a_3\in A_3$ and $b_2\in B_3$.
If $v$ is adjacent to $a_3$, then similarly the holes 
$$v\dd P\dd b_1\dd b_2\dd R_2\dd a_2\dd a\dd v,$$
$$v\dd P\dd b_1\dd b_3\dd R_3\dd a_3\dd v$$
have different parity, a contradiction. So $v$ is anticomplete to $\bigcup_{j\in J\setminus \{1\}} A_j$. 
For each $i\in I$, let $P_i$ be an induced path between $v$ and $B_i$ with interior in $C_i\cup D_i$.
Define 
\begin{eqnarray*}
A_0&=&\{v\}\cup \bigcup_{i\in I}A_i;\\
B_0&=& \bigcup_{i\in I}B_i;\\
C_0&=& \bigcup_{i\in I} C_i\cup (V(P_i)\cap D_i);
\end{eqnarray*}
Then $S_0$ is a strip, and $(a,S_i\;(i\in J\cup \{0\}))$ is an indecomposable pyramid strip system contrary to
the maximality of $V(\mathcal{S})$. 
This proves (2).

\bigskip

To complete the proof of the theorem, we therefore may assume that $I=\emptyset$;
so now let $v\in N(a)$
with no neighbour in $B_i\cup C_i\cup D_i$ for $1\le i\le k$. Since $a$ is splendid, $v$ has a neighbour 
$u\in V(G)\setminus N[a]$; and so $u\notin V(\mathcal{S})\cup N[a]$. Let $F$ be the component of 
$G\setminus (V(\mathcal{S})\cup N[a])$ that contains $u$. Since $F$ is contained in none of the sets $D_i$, it follows that
$\mathcal{S}(F)\subseteq B_1\cupcup B_k$. Choose a minimal path $Q$ of $G[V(F)\cup \{v\}]$ with one end $v$
such that the other end, $q$ say, has a neighbour in $B_1\cupcup B_k$. For $1\le i\le k$ let $B_i'\subseteq B_i$
be the set of vertices in $B_i$ adjacent to $q$, and let $B_i''=B_i\setminus B_i'$. Let $A_i'$ be the set of vertices in $A_i$
adjacent to $v$, and $A_i''=A_i\setminus A_i'$. The only edges between $V(\mathcal{S})$ and $V(Q)$ are the edges between $v$
and $\{a\}\cup A_1\cupcup A_k$, and the edges between $q$ and $B_1\cupcup B_k$, since $Q\setminus v$ is a subgraph of $F$
and $\mathcal{S}(F)\subseteq B_1\cupcup B_k$. 

By a {\em rung} we mean an $S_i$-rung for some $i\in \{1\ll k\}$.
For $1\le i\le k$, let us say an $S_i$-rung $R_i$ is {\em crooked} if it has one end in $A_i$ and the other in $B_i'$,
or one end in $A_i'$ and the other in $B_i$; and {\em straight} otherwise. Choose $x,y\in \{0,1\}$ such that
$Q$ has length $x$ modulo 2, and every rung has length $y$ modulo 2.
\\
\\
(3) {\em If $x\ne y$ then no rung is crooked, and either $v$ is complete to $A_1\cupcup A_k$ (and $v$ has type $\alpha$), or 
for some $i$, $v$ is complete to $\bigcup_{j\ne i} A_j$, and anticomplete to $A_i$, and $q$ is complete to 
$\bigcup_{j\ne i} B_j$, and anticomplete to $B_i$ (and so $v$ has type $\gamma_i$, and $Q$ is a private path).}
\\
\\
Suppose that $R_1$ is a crooked $S_1$-rung, with ends $a_1\in A_1$ and $b_1\in B_1$. If $a_1\in A_1''$ and $b_1\in B_1'$ then
$b_1\dd R_1\dd a_1\dd a\dd v\dd Q\dd q\dd b_1$ is an even hole; so $a_1\in A_1'$ and $b_1\in B_1''$.
If there exists $b_2\in B_2'$, then
$b_1\dd R_1\dd a_1\dd v\dd Q\dd q\dd b_2\dd b_1$ is an even hole, a contradiction; so $B_2'\ll B_k'=\emptyset$.
Hence $B_1'\ne \emptyset$; and so for $2\le i\le k$ there is no crooked $S_i$-rung, by the same argument 
with $S_1,S_i$ exchanged, and so 
$A_2'\ll A_k'=\emptyset$. But then we can add $v$ to $A_1$ and $V(Q)\setminus \{v\}$ to $C_1$ (note that the edge
$va_1$ guarantees the indecomposability of the new strip), contrary to the maximality of
$V(\mathcal{S})$. 

Thus every rung is straight.
Suppose that $A_1', A_1''\ne \emptyset$. Let $C_1'$ be the union of all interior of $S_1$-rungs between $A_1', B_1'$,
and let $C_1''$ be the union of all interiors of $S_1$-rungs between $A_1'', B_1''$. Since every $S_1$-rung
is of one of these two types, $C_1'\cup C_1''=C_1$. Since there is no $S_1$-rung with ends in $A_1'$ and $B_1''$,
it follows that $C_1'\cap C_1''=\emptyset$ and $C_1',C_1''$ are anticomplete. For the same reason, the only edges
between $A_1'\cup C_1'$ and $A_1''\cup C_1''$ are between $A_1'$ and $A_2'$. Since $S_1$ is indecomposable,
there is an edge between some $a_1'\in A_1'$ and some $a_1''\in A_1''$. Let $R_1''$ be an $S_1$-rung
with ends $a_1''$ and some $b_1''\in B_1''$. If there exists $a_2\in A_2'$,
let $R_2'$ be an $S_2$-rung with ends $a_2,b_2$; then
$$b_1''\dd R_1''\dd a_1''\dd a_1'\dd v\dd a_2\dd R_2\dd b_2\dd b_1''$$
is an even hole, a contradiction. So $A_2'\ll A_k'=\emptyset$, and since every rung is straight, it follows that
$B_2'\ll B_k'=\emptyset$. But then we can add $v$ to $A_1$ and $V(Q\setminus v)$ to $C_1$,
contrary to the maximality of $V(\mathcal{S})$.

This proves that for each $i\in \{1\ll k\}$, either $A_i'=B_i'=\emptyset$, or $A_i''=B_i''=\emptyset$.
Let $I$ be the set of $i \in \{1\ll k\}$ such that $A_i'\ne \emptyset$. Suppose that $|I|\le k-2$,
say $I=\{i+1\ll k\}$ where $i\ge 3$. Define $S_0=(A_0,B_0,C_0)$, where 
\begin{eqnarray*}
A_0 &=& \{v\}\cup \bigcup_{i\in I}A_i\\
B_0&=& \bigcup_{i\in I} B_i\\
C_0&=& V(Q\setminus v) \cup \bigcup_{i\in I} C_i.
\end{eqnarray*}
Then $(a,S_0,S_1\ll S_i)$ is an indecomposable pyramid strip system, contrary to the maximality of $V(\mathcal{S})$.
So $|I|\ge k-1$. This proves (3).
\\
\\
(4) {\em If $x=y$ then 
there exists $i$ such that $v$ is complete to $\bigcup_{j\ne i} A_j$, and $q$ is anticomplete to
$\bigcup_{j\ne i} B_j$ (and so $v$ has type $\gamma_i$ and $Q$ is a private path).}
\\
\\
Suppose that $R_1$ is a straight $S_1$-rung, with ends $a_1\in A_1$ and $b_1\in B_1$.
If $a_1\in A_1'$ and $b_1\in B_1'$ then $G[V(R_1\cup Q)]$ is an even hole, which is impossible.
Since $R_1$ is straight, it follows that $a_1\in A_1''$ and $b_1\in B_1''$. If there exists $b_2\in B_2'$, then
$b_1\dd R_1\dd a_1\dd a\dd v\dd Q\dd q\dd b_2\dd b_1$ is an even hole, a contradiction; so $B_2'\ll B_k'=\emptyset$.
Hence $B_1'\ne \emptyset$; and so for $2\le i\le k$ there is no straight $S_i$-rung,
by the same argument with $S_1,S_i$ exchanged.
Hence $A_2''\ll A_k''=\emptyset$, and the claim holds.

Thus we may assume that every rung is crooked.
Suppose that $A_1', A_1''\ne \emptyset$. Let $C_1'$ be the union of all interior of $S_1$-rungs between $A_1', B_1''$,
and let $C_1''$ be the union of all interiors of $S_1$-rungs between $A_1'', B_1'$. Since every $S_1$-rung
is of one of these two types, $C_1'\cup C_1''=C_1$. Since there is no $S_1$-rung with ends in $A_1'$ and $B_1'$,
it follows that $C_1'\cap C_1''=\emptyset$ and $C_1',C_1''$ are anticomplete. For the same reason, the only edges 
between $A_1'\cup C_1'$ and $A_1''\cup C_1''$ are between $A_1'$ and $A_2'$. Since $S_1$ is indecomposable,
there is an edge between some $a_1'\in A_1'$ and some $a_1''\in A_1''$. Let $R_1$ be an $S_1$-rung
with ends $a_1''$ and some $b_1'\in B_1'$. If there exists $a_2\in A_2'$,
let $R_2$ be an $S_2$-rung with ends $a_2,b_2$; then 
$$b_1'\dd R_1\dd a_1''\dd a_1'\dd v\dd a_2\dd R_2\dd b_2\dd b_1'$$
is an even hole, a contradiction. So $A_2'\ll A_k'=\emptyset$, and since every rung is crooked, it follows that
$B_2''\ll B_k''=\emptyset$. But then we can add $v$ to $A_1$, $q$ to $B_1$, and $Q^*$ to $C_1$,
contrary to the maximality of $V(\mathcal{S})$.

This proves that for each $i\in \{1\ll k\}$, either $A_i'=B_i''=\emptyset$, or $A_i''=B_i'=\emptyset$.
Let $I$ be the set of $i \in \{1\ll k\}$ such that $A_i'\ne \emptyset$. If $I=\emptyset$,
define $S_0=(\{v\}, \{q\}, Q^*)$, then $(a,S_0, S_1\ll S_k)$ is an indecomposable pyramid strip system, contrary to the
maximality of $V(\mathcal{S})$. So $I\ne \emptyset$. Suppose that $|I|\le k-2$,
say $I=\{i+1\ll k\}$ where $3\le i\le k$.  Define $S_0=(A_0,B_0,C_0)$, where
\begin{eqnarray*}
A_0 &=& \{v\}\cup \bigcup_{i\in I}A_i\\
B_0&=&\{q\}\cup \bigcup_{i\in I} B_i\\
C_0&=& Q^* \cup \bigcup_{i\in I} C_i.
\end{eqnarray*}
Then $(a,S_0,S_1\ll S_i)$ is an indecomposable pyramid strip system, contrary to the maximality of $V(\mathcal{S})$.
So $|I|\ge k-1$ and again the claim holds. This proves (4).

\bigskip
From (3) and (4) it follows that $v$ has type $\gamma_i$, and $Q$
is the corresponding private path.
In view of (2), this proves \ref{growstrips}. ~\bbox

We say $A$ {\em meets} $B$ if $A\cap B\ne \emptyset$.

\begin{thm}\label{getclique}
Let $G$ be even-hole-free, and let $a\in V(G)$ be splendid. Suppose
there is no extended near-prism contained in $G$ such that $a$ is an end of its cross-edge. 
Let $\mathcal{S}=(a,S_1\ll S_k)$ be
an indecomposable strip system with apex $a$, with strips $S_i=(A_i,B_i,C_i)$ for $1\le i\le k$,
chosen with $V(\mathcal{S})$ maximal.
Then $N[a]\setminus V(\mathcal{S})$ is a clique.
\end{thm}
\Proof
For $X,Y\subseteq V(G)$, an
{\em $X-Y$ path} means (in this proof) an induced path $P$ of $G$ with ends $x,y$ say, where $X\cap V(P)=\{x\}$ 
and $Y\cap V(P)=\{y\}$ (possibly $x=y$ and $V(P)=\{x\}$, if $x\in X\cap Y$). 
If $X\subseteq V(G)$, a path of $G$ is said to be {\em within} $X$ if $V(P)\subseteq X$.
Let $B=B_1\cupcup B_k$.
\\
\\
(1) {\em For each $v\in N[a]\setminus V(\mathcal{S})$ there exists $x(v)\in \{0,1\}$ such that for $1\le i\le k$,
every $N(v)-B$ path within $V(S_i)$ has parity $x(v)$.}
\\
\\
Since $v\in N(a)\setminus V(\mathcal{S})$, \ref{growstrips}  implies that there are at least two values of $i\in \{1\ll k\}$
such that $N(v)\cap V(S_i)\ne \emptyset$; and for each such $i$ there is an $N(v)-B$ path within $V(S_i)$. Let 
$N(v)\cap V(S_i)\ne \emptyset$ for $i = 1,2$ say, and for $i = 1,2$ let $P_i$ be an $N(v)-B$ path within $V(S_i)$.
Then $G[V(P_1\cup P_2)\cup \{v\}]$ is a hole, and so $P_1$, $P_2$ have the same parity, say $x(v)\in \{0,1\}$.
We claim that 
for $1\le j\le k$, every $N(v)-B$ path $P$ in $V(S_j)$ has parity $x(v)$. To see this, choose 
$i\in \{1,2\}$ different from $j$; then $G[V(P_i\cup P)\cup \{v\}]$ is a hole, and the claim follows. This proves (1).

\bigskip

In particular, $x(a)$ exists, and so for $1\le i\le k$, all $S_i$-rungs have parity $x(a)$.
Suppose that $u,v\in N(a)\setminus V(\mathcal{S})$ are nonadjacent. 
\\
\\
(2) {\em If $X_1,X_2$ are connected subsets of $V(G)$, disjoint and anticomplete, and 
$u,v$ both have neighbours in $X_i$ for $i = 1,2$, then all $N(u)-N(v)$ paths within $X_1$ have the same parity, and all $N(u)-N(v)$ paths
within $X_2$ have the opposite parity.}
\\
\\
For $i = 1,2$, let $P_i$ be an $N(u)-N(v)$ path $P_i$ within $X_i$; then $G[V(P_1\cup P_2)\cup \{u,v\}]$ is a hole,
and so $P_1,P_2$ have opposite parity. This proves (2).
\\
\\
(3) {\em There do not exist three connected subsets $X_1,X_2,X_3$ of $V(G)$, pairwise disjoint and pairwise anticomplete,
such that for $i = 1,2,3$, $u,v$ both have neighbours in $X_i$.}
\\
\\
This is immediate from (2).
\\
\\
(4) {\em There is at most one $i\in \{1\ll k\}$ such that $N(u)\cap (A_i\cap C_i)=\emptyset$, and the same for $N(v)$.}
\\
\\
Suppose that $N(u)$ is disjoint from $A_i\cup C_i$ for $i = 1,2$. By \ref{growstrips},
$N(a)$ meets at least $k-1$ of $V(S_1)\ll V(S_k)$, so we may assume there exists $b_1\in B_1\cap N(u)$. Choose $b_2\in B_2$,
and for $i = 1,2$ let $R_i$ be an $S_i$-rung containing $b_i$. If $b_2,u$ are adjacent,
there is a short pyramid with apex $a$, with base $\{b_1,b_2,u\}$ and constituent paths
$R_1,R_2$ and the edge $u\dd a$, which is impossible since $a$ is splendid. If $b_2,u$ are nonadjacent,
there is a theta with ends $b_1,a$ and constituent paths $b_1\dd R_1i\dd a, b_1\dd u\dd a$, and $b_1\dd b_2\dd R_2\dd a$,
contrary to \ref{subgraphs}. This proves (4).
\\
\\
(5) {\em $k=3$, and there exists $i\in \{1,2,3\}$ such that $u,v$ both have neighbours in $A_i\cup C_i$.}
\\
\\
Since each $S_i$ is indecomposable, there are only at most two values of $i$
such that $N(u), N(v)$ both meet $A_i\cup C_i$, by (3).
Then both claims follow from (4). This proves (5).

\bigskip

As before, for $i = 1,2,3$, 
let $D_i$ be the union of all components $F$ of
$G\setminus (V(\mathcal{S}\cup N[a])$ such that $\mathcal{S}(F)\cap (A_i\cup C_i)\ne \emptyset$.
From (3), there exists $i\in \{1,2,3\}$ such that not both $u,v$ have neighbours in $A_i\cup C_i\cup D_i$.
\\
\\
(6) {\em If $v$ is anticomplete to $V(S_3)\cup D_3$ then $v$ has type $\beta_3$.}
\\
\\
Suppose not. Certainly $v$ does not have type $\alpha$ or $\alpha'$, 
since it has no neighbour in $V(S_3)\cup D_3$. It does not have type 
$\beta_1$ or $\beta_2$ since it has no neighbour in $B_3\cup C_3$; and not type $\gamma_1, \gamma_2$ since
it is not complete to $A_3$. So $v$ has type $\gamma_3$;
let $Q$ be the corresponding private path, between $v$ and $q$ say, and let $p$ be the neighbour of $v$ in this path. 
Also, since $v$ is complete to $A_1$ and anticomplete to $B_1\cup C_1$, it follows that $x(v)=x(a)$.

For $i = 1,2$, if $N(u)$ meets $A_i\cup C_i\cup D_i$, then there is an $N(u)-\{a\}$ path $R$ within $\{a\}\cup A_i\cup C_i\cup D_i$;
and since its ends are adjacent to $u$, it has odd length. Hence $R\setminus a$ is an $N(u)-N(v)$ path 
(since $v$ is complete to $A_i$ and anticomplete to $B_i\cup C_i$),
and has even length. By (2), $N(u)$ is disjoint from one of $A_1\cup C_1\cup D_1, A_2\cup C_2\cup D_2$, say
$A_2\cup C_2\cup D_2$; and by (4), $N(u)$ meets $A_1\cup C_1$. Let $P_1$ be an even $N(u)-N(v)$ path within $A_1\cup C_1$.

Now $u$ has no neighbour in $A_2\cup C_2\cup D_2$. Suppose that $u$ has a neighbour in the connected set
$C_3\cup B_3\cup B_2\cup V(Q\setminus v)$, and let $T$ be an $N(u)-\{a\}$ path within
$$C_3\cup B_3\cup V(Q\setminus v)\cup B_2\cup C_2\cup A_2\cup \{a\}.$$ 
This path has odd length (because its ends 
are neighbours of $u$), and it contains no neighbour of $v$ except the one in $A_2$ (because $p$ is nonadjacent to $u$).
Consequently the path $T\setminus a$ is an $N(u)-N(v)$-path of even length anticomplete to $P_1$, a contradiction.
So $u$ has no neighbour in $C_3\cup B_3\cup V(Q)\cup B_2$. Since 
$u$ is anticomplete to $V(S_2)\cup D_2$, \ref{growstrips} implies that $u$ has type $\gamma_2$, and in particular,
$u$ is complete to $A_3$ and has no neighbour in $C_3$.
Let $T$ be an $N(v)-\{a\}$ path within $V(Q)\cup V(S_3)\cup \{a\}$; again it has odd length (since its ends
are adjacent to $v$), and $T\setminus a$ is an even $N(u)-N(v)$-path anticomplete to $P_1$,
a contradiction. This proves (6).
\\
\\
(7) {\em There is only one $i\in \{1,2,3\}$ such that both $N(u), N(v)$ meet $A_i\cup C_i\cup D_i$.}
\\
\\
Suppose that $N(u), N(v)$ both meet $A_i\cup C_i\cup D_i$ for $i = 1,2$. Then by (3), one of $u,v$ has no neighbours
in $V(S_3)\cup D_3$, say $v$. By (6), $v$ has type $\beta_3$, and so has a neighbour in $B_1\cup C_1$ and one in $B_2\cup C_2$.
For $i = 1,2$, let $P_i$ be an $N(u)-N(v)$ path within $A_i\cup C_i$.
By exchanging $S_1,S_2$ if necessary, we may assume that $P_1$ has odd length, and so $P_2$ is even. 
Hence there is no $N(u)-N(v)$ path within the connected set $B_2\cup B_3\cup C_2\cup C_3\cup D_2\cup D_3$,
because we could combine it with one of $u\dd a\dd v$ and $u\dd P_1\dd v$ to make an even hole. Since $v$ has a neighbour
in this set, $u$ does not. So $u$ does not have type $\beta$. By (6), $u$ has a neighbour $a_3\in A_3$. 
Let $R_3$ be an $S_3$-rung containing $a_3$, and for $i = 1,2$, let $R_i$ be an $N(v)-B_i$
path within $B_i\cup C_i$. For $i = 1,2,3$, let $b_i$ be the end of $R_i$ in $B_i$.
Thus $R_1, R_2$ both have parity $x(v)$. For $i = 1,2$, let $Q_i$ be the induced path $R_i\dd b_i\dd b_3\dd R_3$. Thus
$Q_2$ is an $N(u)-N(v)$ path, but $Q_1$ might not be.
Now $Q_1, Q_2$ have the same parity. Since 
$Q_2$ is anticomplete to $P_1$ it follows that $Q_2$ is even, and hence $Q_1$ is even; and since $Q_1$ is anticomplete to $P_2$,
it follows that $Q_1$ is not an $N(u)-N(v)$ path. But it has one end in $N(v)$ and no other vertex in $N(v)$; and its other
end is in $N(u)$. Consequently some internal vertex is in $N(u)$, and so $u$ has a neighbour in $V(R_1)$.

If $u$ has a unique neighbour $t\in V(R_1)$, there is a theta with ends $t,v$ and constituent paths
$$t\dd R_1\dd v,$$
$$t\dd u\dd a\dd v,$$
$$t\dd R_1\dd b_1\dd b_2\dd R_2\dd v,$$
contrary to \ref{subgraphs}. (Note that $t,v$ are nonadjacent since $u,v$ have no common neighbour nonadjacent to $a$.)
If $u$ has two nonadjacent neighbours in $V(R_1)$, there is a theta with ends $u,v$ and constituent paths
$$u\dd R_1\dd v,$$
$$u\dd a\dd v,$$
$$u\dd R_1\dd b_1\dd b_2\dd R_2\dd v,$$
contrary to \ref{subgraphs}. If $u$ has exactly two adjacent neighbours $p,q\ in V(R_1)$, where $v,p,q,b_1$ 
are in order in $R_1$, there is a near-prism with
bases $\{v,p,q\}$ and $\{b_1,b_2,b_3\}$ and constituent paths 
$$p\dd R_1\dd v\dd R_2\dd b_2,$$
$$u\dd R_3\dd b_3,$$
$$q\dd R_1\dd b_1,$$
contrary to \ref{subgraphs}. 
This proves (7).

\bigskip

In view of (4), (5) and (7), we may assume that $u,v$ both have neighbours in $A_1\cup C_1$; $v$ has a neighbour
in $A_2\cup C_2$ and none in $A_3\cup C_3\cup D_3$, and $u$ has a neighbour in $A_3\cup C_3$ and none in $A_2\cup C_2\cup D_2$.
\\
\\
(8) {\em $u$ has no neighbour in $B_2$, and $v$ has no neighbour in $B_3$.}
\\
\\
Suppose that $v$ has a neighbour in $B_3$, say $b_3$, and so $x(v)=0$. Let $R_3$ be an $S_3$-rung with ends $a_3, b_3$. 
The path $a\dd a_3\dd R_3\dd b_3$ is odd, since its ends
are neighbours of $v$, and so $x(a) = 0$.

Suppose first that $x(u)=0$. There is an $N(u)-N(v)$ path with one end $b_3$ and otherwise contained in $A_3\cup C_3$.
Its length has parity $x(u)$, and it is anticomplete to $P_1$, where $P_1$ is an $N(u)-N(v)$ path within $A_1\cup C_1$;
so $P_1$ has odd length by (2). Hence there is no $N(u)-N(v)$ path within the connected set    
$B_2\cup C_2\cup D_2\cup B_3\cup C_3\cup D_3$, and so $u$ is anticomplete to this set. By (6) $u$ has type $\beta_2$,
a contradiction since $u$ has no neighbour in $B_3\cup C_3$. 

This shows that $x(u)=1$, and hence $u$ has no neighbour in $B$. 
Let $R_1$ be an $S_1$-rung with ends $a_1\in A_1$ and $b_1\in B_1$, that contains
a neighbour of $u$, and let $T$ be an $N(u)-B_1$ subpath of $R_1$. Thus $T$ has parity $x(u)$ and hence is odd, 
and so $a_1\notin V(T)$ since $x(a)=0$.
Consequently $u$ has a neighbour in $R_1^*$. Since the connected sets $\{a\}$, $R_1^*$ and
$V(S_3)$ are pairwise anticomplete, (3) implies that $v$ has no neighbour in $R_1^*$. But the path $T\dd b_1\dd b_3$
is even, and anticomplete to $\{a\}$; and so this path is not an $N(u)-N(v)$ path, and so $v$ has a neighbour in 
$T$, and therefore $v,b_1$ are adjacent. Since $R_1$ is even, and $v$ has no neighbour in $R_1^*$, it follows
that $v,a_1$ are not adjacent. But then there is a short pyramid with apex $a$, base $\{v,b_1,b_3\}$,
and constituent paths 
$$a\dd a_1\dd R_1\dd b_1,$$
$$a\dd v,$$
$$a\dd a_3\dd R_3\dd b_3,$$
contradicting that $a$ is splendid. This proves (8).

\bigskip

Thus $u$ has no neighbour in $V(S_2)\cup D_2$, and $v$ has no neighbour in $V(S_3)\cup D_3$. By (6), $u$ has type $\beta_2$
and $v$ has type $\beta_3$. Since $v$ has a neighbour in $B_2\cup C_2$, there is an $S_2$-rung $R_2$ with ends $a_2\in A_2$
and $b_2\in B_2$, such that $v$ has a neighbour in $R_2$ different from $a_2$. Choose an $S_3$-rung $R_3$ with ends $a_3, b_3$
similarly for $u$. Now $v$ has two nonadjacent neighbours in the hole 
$$a\dd a_2\dd R_2\dd b_2\dd b_3\dd R_3\dd a_3\dd a,$$
and hence it has at least three, and an odd number; and they all belong to $R_2$ except $a$. Similarly 
$R_3$ contains a positive even number of neighbours of $u$. Also, the hole 
$$v\dd R_2\dd b_2\dd b_3\dd R_3\dd a_3\dd a\dd v$$
is odd, and so $x(v)\ne x(a)$, and similarly $x(u)\ne x(a)$. 
\\
\\
(9) {\em Every $S_1$-rung contains an even number of neighbours of $v$, and an even number of neighbours of $u$.}
\\
\\
Let $R_1$ be an $S_1$-rung with ends $a_1\in A_1$ and $b_1\in B_1$. Since 
$$a\dd a_1\dd R_1\dd b_1\dd b_2\dd R_2\dd a_2\dd a$$
is a hole, and the path $R_2\dd a_2\dd a$ contains an odd number at least three of neighbours of $v$, 
and the total cannot be even and at least three, it follows that there is an 
even number of neighbours of $v$ in $R_1$.  Similarly $R_1$ contains an even number of neighbours of $u$.
This proves (9).
\\
\\
(10) {\em For every $N(u)-N(v)$ path $P_1$ within $A_1\cup C_1$, $P_1$ has even length, and either $V(P_1)\subseteq A_1$, or 
one end of $P_1$ belongs to $A_1$ and its other vertices belong to $C_1$. In particular, $A_1\cap V(P)\ne \emptyset$.}
\\
\\
There is an $N(u)-N(v)$ path $Q$ within $B_2\cup C_2\cup B_3\cup C_3$, and it
is anticomplete to $\{a\}$ and so odd; and it is also anticomplete to $P_1$,
and so $P_1$ is even. Now
$u\dd P_1\dd v\dd Q\dd u$ is a hole $H$ say, and the neighbours of $a$ in it are $u, v$, and all vertices of $V(P_1)\cap A_1$.
Since $a$ is splendid and therefore $V(G)\setminus N[a]$ is connected, \ref{bigvertex} implies that either
\begin{itemize}
\item $a$ is complete to $H$; or 
\item the subgraph induced on the set of vertices of $H$ adjacent to $a$ is a path; or
\item $a$ has exactly three neighbours in $H$, and two of them are adjacent.
\end{itemize}
The first is impossible since $a$ is not complete to $V(Q)$. The second implies that $V(P_1)$ is complete to $a$, that is,
$V(P_1)\subseteq A_1$; and the third implies that one end of $P_1$ belongs to $A_1$ and the others belong to $C_1$.
This proves (10).
\\
\\
(11) {\em No $S_1$-rung meets both $N(v_1)$ and $N(v_2)$.}
\\
\\
Let $R_1$ be an $S_1$-rung with ends $a_1\in A_1$ and $b_1\in B_1$. By (10), not both $N(u), N(v)$ meet $R_1^*$, so we may 
assume that $N(v)\cap V(R_1)= \{a_1,b_1\}$ (since it has even cardinality by (9)). Thus $b_1\notin N(u)$, and so $N(u)$
meets $R_1^*$ by (9). Since $u$ has an even number of neighbours in $V(R_1)$, and $v\dd a_1\dd R_1\dd b_1\dd v$ is a hole,
and there is no even wheel and no theta, it follows that $u$ has exactly two neighbours in $R_1$ and they are adjacent. 
But then the subgraph induced on $V(R_1)\cup \{u,v,a\}$ is a near-prism, contrary to \ref{subgraphs}. 
This proves (11).
\\
\\
(12) {\em There is no $N(u)-N(v)$ path within $A_1\cup C_1$ with one end in $A_1$ and all other
vertices in $C_1$.}
\\
\\
Suppose there is such a path, $P$ say.
Let $P$ have ends $p\in A_1\cap N(u)$ and $q\in N(v)$ (possibly $p=q$), with $V(P)\setminus \{p\}\subseteq C_1$.
If $p=q$, an $S_1$-rung with one end $p$ contradicts (11); so $p\ne q$. 
Let $R_1$ be an $S_1$-rung with ends $a_1\in A_1$ and $b_1\in B_1$, containing $q$. The path $p\dd P\dd q\dd R_1\dd b_1$
includes an $S_1$-rung with one end in $N(u)$, and therefore contains another neighbour of $u$ by (9). 
This does not belong to 
$V(P)$, so it belongs to $V(R_1)$; and so $V(R_1)$ meets both $N(u)$ and $N(v)$, contrary to (11). This proves (12).

\bigskip

From (10) and (12), every $N(u)-N(v)$ path within $A_1\cup C_1$ is within $A_1$.
Choose $P_1$ as in (10) to have as few vertices in $A_1$ as possible.
It follows that 
$V(P_1)\subseteq A_1$. Let $P_1$ have ends $p,q$, where $p$ is adjacent to $u$ and $q$ to $v$.
From (11) $p\ne q$.
Let $R_1$ be an $S_1$-rung with one end $p$, and let $b_1$ be the end of $R_1$ in $B_1$. 
By (11), $v$ has no neighbour in $V(R_1)$. Now $V(P_1\setminus p)$ is disjoint from $V(R_1\setminus p)$;
suppose these two sets are anticomplete.
Then $q\dd P_1\dd p\dd R_1\dd b_1$ 
is an $N(v)-B_1$ path, and so
it has parity $x(v)$. But its parity is the same as that of $R_1$, since $P_1$ is even; and so $x(v)=x(a)$, a contradiction.
Hence $V(P_1\setminus p)$ is not anticomplete to $V(R_1\setminus p)$.

Suppose that $V(P_1\setminus p)$ is not anticomplete to $R_1^*$. Since every $N(u)-N(v)$ path within $A_1\cup C_1$ is within $A_1$,
it follows that
no vertex of $R_1^*$ is adjacent to $u$. But from (11), at least two vertices of $R_1$
are adjacent to $u$, and so $b_1$ is adjacent to $u$. Since $V(P_1\setminus p)$ is not anticomplete to $V(R_1\setminus p)$,
there is an $S_1$-rung with one end $b_1$ and the other in $V(P_1\setminus p)$, and this $S_1$-rung therefore contains
a unique neighbour of $u$, contrary to (9).

Thus $V(P_1\setminus p)$ is anticomplete to $R_1^*$, and so $b_1$ has a neighbour $r\in V(P_1\setminus p)$. By (9), $u$
has a neighbour in $V(R_1\setminus p)$, and so there is an induced path $Q$ between $u,b_1$ with interior in $R_1^*$.
Hence $Q$ has parity $x(u)+1$, and since the path $r\dd b_1$ is an $S_1$-rung and so has parity $x(a)\ne x(u)$,
it follows that $a\dd u\dd Q\dd b_1\dd r\dd a$ is an even hole, a contradiction. This proves \ref{getclique}.~\bbox

\section{Using the decomposition theorems}

Let $S=(A,B,C)$ be a strip in a graph $G$, and let $a\in V(G)\setminus V(S)$ be complete to $A$ and anticomplete to 
$B\cup C$. 
Let $D$ be the union of all the vertex sets of all components $F$ of $G\setminus (V(S)\cup N[a])$ such that $F$ is 
not anticomplete to $A\cup C$, and let $Z$ be the set of all vertices in $V(G)\setminus V(S)$ that are adjacent or equal
to $a$ and have a neighbour in $A\cup C\cup D$.
For $v\in Z$, a {\em backdoor} for $v$ is an induced path $R$ of $G$ with ends
$v,b$ say, such that
$R^*$ is anticomplete to $V(S)\cup D$, and 
$b$ is complete to $B$ and has no neighbours in $A\cup C\cup D$.
We say $(S,a, D, Z)$ is  a {\em completed strip} if
\begin{itemize}
\item $S$ is proper;
\item $Z$ is a clique; and
\item every vertex in $Z$ has a backdoor.
\end{itemize}
We will see that both our decomposition theorems yield completed strips; and
completed strips are good for finding bisimplicial vertices by induction, because of the following.
\begin{thm}\label{striptobip}
Let $G$ be even-hole-free, such that \ref{mainthm} holds for all graphs with fewer vertices than $G$. Let $(S, a,D,Z)$
be a completed strip in $G$, where $S=(A,B,C)$. Let there be at least three vertices in $G$ that are not in $A\cup C\cup D$
and have no neighbour in this set.
Then some vertex in $A\cup C\cup V(F)$ is bisimplicial in $G$.
\end{thm}
\Proof
For each $z\in Z$, let $R_z$ be a backdoor for $z$. 
Let $Z_1$ be the set of all $z\in Z$ such that $R_z$
has odd length, and $Z_2$ the set for which $R_z$ has even length.
\\
\\
(1) {\em If $v\in A\cup C\cup D$, then every neighbour of $v$ in $G$ belongs to $V(S)\cup D\cup \{a\}\cup Z$.}
\\
\\
Suppose $u\in V(G)$ is adjacent to $v$, and $u\notin V(S)\cup D\cup \{a\}\cup Z$. Thus $u$ is not adjacent to $a$, since $u\notin Z$
and $v\in A\cup C\cup D$. If $v\in A\cup C$ then $u\in D$ from the definition of $D$; and if $v\in D$, let $v\in V(F)$ where $F$
is a component of $G\setminus V(S)$ such that $F$ is anticomplete to $a$
and not anticomplete to $A\cup C$; then $u$ also belongs to $V(F)$ and hence to $D$, a contradiction. This proves (1),
\\
\\
(2) {\em If $z\in Z_2$, every induced path between $z,B$ with interior in $A\cup C\cup D$ is even.}
\\
\\
Let $P$ be an  induced path between $z$ and some $b'\in B$ with interior in $A\cup C\cup D$; then $V(P\cup R_z)$ induces
an odd hole, and since $R_z$
is even it follows that $P$ is even. This proves (2).
\\
\\
(3) {\em If $z\in Z_1$, every induced path between $z$ and $B$ with interior in $Z_2\cup A\cup C\cup D$ is odd.}
\\
\\
Let $P$ be an  induced path between $z$ and some $b'\in B$ with interior in $Z_2\cup A\cup C\cup D$. 
If $Z_2\cap V(P)=\emptyset$, then $V(P\cup R_z)$
induces an odd hole, and since $R_z$ is odd it follows that $P$ is odd. So we may assume that there exists $z_2\in V(P)\cap Z_2$.
Since $Z_1\cup Z_2$ is a clique, $z_2$ is unique, and is the neighbour of $z_1$ in $P$. Thus $P\setminus z_1$ is an induced path
between $z_2$ and $b$ with interior in $A\cup C\cup D$, and so is even by (2); and so $P$ is odd. This proves (3).

\bigskip

Let $G'$ be the graph obtained from $G[V(S)\cup D\cup Z]$ by adding two new vertices $b,c$,
where $b$ is complete to $B\cup Z_1$ and $c$ is complete to $Z\cup \{b\}$.
We claim that $G'$ is even-hole-free. To see this, 
suppose that $H$ is an even hole in $G'$. Since $G$ is even-hole-free, $H$ contains at least one of $b,c$; and 
if $H$ contains $c$ then it also contains $b$ since the other $G'$-neighbours of $c$ are a clique. 
Thus $b\in V(H)$. If both $H$-neighbours of $b$ belong to $B$, then there is an induced subgraph of the even-hole-free graph
$G[V(S)\cup D\cup Z\cup V(R_a)]$ isomorphic to $H$, which is impossible. Thus 
$b$ is $H$-adjacent to some vertex $z_1\in Z_1\cup \{c\}$.
Since $b$ is $G'$-complete to $Z_1\cup \{c\}$, only one vertex of $H$ belongs to this set. Consequently
the other $H$-neighbour of $b$ belongs to $B$, and $|V(H)\cap B|=1$. If $z_1\in Z_1$ then $c\notin V(H)$ and
$H\setminus b$ is an even induced path of $G$ between $z_1$ and $B_1$
with interior in $Z_2\cup A\cup C\cup D$, contrary to (3).
Thus $z_1=c$, and hence $V(H)\cap Z_1=\emptyset$, and the other $H$-neighbour of $c$ is some $z_2\in Z_2$.
But then $H\setminus \{b,c\}$ is an odd induced path of $G$ between $z_2, B$ with interior in $A\cup C\cup D$, contrary to (2).
This proves that $G'$ is even-hole-free.

Now $A\ne \emptyset$, and so 
$bc$ is a non-dominating clique of $G'$, since $S$ is proper. But $|V(G')|< |V(G)|$, since every vertex of $G'$ except $b,c$ belongs to $V(G)$
and is not anticomplete to $A\cup C\cup D$.
From the inductive hypothesis, there is a vertex $v\in V(G')\setminus N_{G'}[b,c]$
that is bisimplicial in $G$. Consequently $v\in A\cup C\cup D$. Since $v$ is nonadjacent to $b$, all edges of $G'$ with both ends in $N_{G'}(v)$ are edges
of $G$. But all neighbours of $v$ in $G$ are neighbours of $v$ in $G'$, by (1); and so $v$ is bisimplicial in $G$.
This proves \ref{striptobip}.~\bbox

In order to prove \ref{mainthm}, we will show:
\begin{thm}\label{induct}
Let $G$ be an even-hole-free graph, such that \ref{mainthm} holds for all graphs with fewer vertices than $G$.
Let $K$ be a non-dominating clique in $G$ with $|K|\le 2$.
Then some vertex in $V(G)\setminus N[K]$ is bisimplicial.
\end{thm}
We divide the proof into four parts. First we need:

\begin{thm}\label{induct1}
Let $G$ be an even-hole-free graph, such that \ref{mainthm} holds for all graphs with fewer vertices than $G$.
Let $K$ be a non-dominating clique in $G$ with $|K|\le 2$, and let $a\in V(G)\setminus N[K]$ be splendid, and such that
there is an extended near-prism in $G$ with cross-edge $ab$ for some $b$.
Then some vertex in $V(G)\setminus N[K]$ is bisimplicial.
\end{thm}
\Proof
Choose 
a tree $J$ and a $J$-strip system $M$ in $G$ with the same cross-edge $ab$, with $(J,M)$ optimal for $ab$.
Let $Z$ be the set of all vertices adjacent to both $a,b$, and $Y$ the set of major vertices.
Let $(\alpha,\beta)$ be the corresponding partition. For each $e=st\in E(J)$ with $t\in \alpha$, let $D_e$
be the union of the vertex sets of all components of $G\setminus (V(M)\cup Z)$
that are not anticomplete to $M_e\setminus M_s$. By \ref{splendidprism}, if $F'$ is such a component then $a,b$
have no neighbour in $F'$, and every vertex in $V(M)$ with a neighbour in $F'$ belongs to $M_e$. 
\\
\\
(1) {\em For each edge $e=st$ of $J$ with $t\in \alpha$, there is a bisimplicial vertex of $G$ in
$(M_e\setminus M_s)\cup D_e$, where $D_e$ is the union of the vertex sets of all components of $G\setminus (V(M)\cup Z)$
that are anticomplete to $a$ and not anticomplete to $M_e\setminus M_s$.}
\\
\\
Let $A=M_t\cap M_e$, $B=M_s\cap M_e$, $C=M_e\setminus (M_s\cup M_t)$ and $D=D_e$; then $S=(A,B,C)$ is a strip, and it is 
proper, by \ref{splendidprism}.
Let $Z'$ be the set of all vertices in $V(G)\setminus V(S)$ that are adjacent or equal
to $a$ and have a neighbour in $A\cup C\cup D$. We claim that all vertices in $Z'$ are major. Let $z\in Z'$. Then $\{z\}$
is not small, since $a$ has a neighbour in $\{z\}$, and so $b,z$ are adjacent; and hence $z\in Z$. Since $z$ has a neighbour in $V(S)$,
and $b$ has no neighbour in $V(S)$, it follows that $z$ is $b$-external; and since $a$ is splendid, every vertex in $N(a)$
is $a$-external. This proves that $z\in Y$. Consequently $Z'$ is a  clique, by \ref{clique}.

Choose $t'\in \beta$, and let $P$ be a path of $J$ with ends $s,t'$. Choose an $f$-rung $R_f$
for each $f\in E(P)$. Let $u,v$ be the ends of $R_P$, where $u\in M_{t'}$ and $v\in M_s$. For each $z\in Z'$, since
$z$ is adjacent to $b$, there is a path from $z$ to $v$ with interior in $V(R_P)\cup \{b\}$; and 
this is a backdoor for $z$ since 
$v$ is complete to $B$ and anticomplete to $A\cup C$. 

Now $D$ is the union of the vertex sets of all components $F$ of $G\setminus (V(M)\cup Z)$
that are not anticomplete to $M_e\setminus M_s$. By \ref{splendidprism}, for each such $F$, $a$ has no neighbour in $V(F)$;
and so 
$D$ is the union of the vertex sets of all components $F$ of $G\setminus (V(S)\cup N[a])$ such that $F$ is 
not anticomplete to $A\cup C$. Hence $(S,a, D, Z')$ is a completed strip, and there are at least three vertices of $G$
that are anticomplete to $A\cup C\cup D$, namely $b$ and at two vertices of $M_{e'}$ (the latter has at least two vertices, since
the corresponding strip is proper by \ref{splendidprism}). From \ref{striptobip}, there is a bisimplicial vertex of $G$ in
$A\cup C\cup D$. This proves (1).

\bigskip

Choose edges $e=st$ and $e'=s't'$ of $J$ where $t,t'\in \alpha$ are distinct; then by (1), there are bisimplicial vertices 
$v\in (M_e\setminus M_s)\cup D_e$, and $v'\in (M_{e'}\setminus M_{s'})\cup D_{e'}$, defining $D_e, D_{e'}$ as in (1).
Suppose they both belong to $N[K]$. Now for $k\in K$, $k$ is not adjacent to $a$ since $a\in V(G)\setminus N[K]$ by hypothesis;
and so $k\notin Z$. We may choose $k\in K$ adjacent or equal to $v$, and so $k$ is not anticomplete to $(M_e\setminus M_s)\cup D_e$.
Consequently $k\in M_e\cup D_e$. Similarly there exists $k'\in K\cap (M_{e'}\cup D_{e'})$. But $M_e\cup D_e$ is anticomplete to
$M_{e'}\cup D_{e'}$ by \ref{treestructplus}, a contradiction. This proves that one of $v,v'$ is anticomplete to $K$, and so
satisfies the theorem. This proves \ref{induct1}.~\bbox

Second, we need:
\begin{thm}\label{induct2}
Let $G$ be an even-hole-free graph, such that \ref{mainthm} holds for all graphs with fewer vertices than $G$.
Let $K$ be a non-dominating clique in $G$ with $|K|\le 2$, and let $a\in V(G)\setminus N[K]$ be splendid. Suppose that
there is no extended near-prism in $G$ such that $a$ is an end of its cross-edge, and there is a pyramid in $G$ with apex $a$.
Then some vertex in $V(G)\setminus N[K]$ is bisimplicial.
\end{thm}
\Proof
From \ref{growstrips} since there is a pyramid with apex $a$, and all its constituent paths have length at least two (because $a$ is
splendid), there is an indecomposable strip system with apex $a$. Let $\mathcal{S}=(a,S_1\ll S_k)$ be
an indecomposable strip system with apex $a$, with strips $S_i=(A_i,B_i,C_i)$ for $1\le i\le k$,
chosen with $V(\mathcal{S})$ maximal. In the notation of \ref{growstrips}, 
for $1\le i\le k$, let $D_i$ be the union of the vertex sets of all components $F$ of
$G\setminus (V(\mathcal{S}\cup N[a])$ such that $\mathcal{S}(F)\cap (A_i\cup C_i)\ne \emptyset$.
\\
\\
(1) {\em For $1\le i\le k$, there is a bisimplicial vertex of $G$ in $A_i\cup C_i\cup D_i$.}
\\
\\
Let $1\le i\le k$, $i = 1$ say; 
and let $Z$ be the set of all $z\in N(a)\setminus V(\mathcal{S})$ such that $z$ has a neighbour in $A_1\cup C_1\cup D_1$. Thus $Z$
is a clique by \ref{getclique}. We need to show that each $z\in Z$ has a backdoor. By \ref{growstrips}, $z$ has type 
$\alpha$, $\alpha'$, $\beta$ or $\gamma$, and hence for some $2\le j\le k$, $z$ has a neighbour in $V(S_j)$. 
Choose an $S_j$-rung $R$
in which $z$ has a neighbour, with an end $b\in B_j$ say; then a path between $z,b$ with interior in $V(R)$ provides a backdoor.
Thus each $z\in Z$ has a backdoor; and there are at least three vertices in $G$ that are anticomplete to $A_1\cup C_1\cup D_1$,
for instance all vertices of $A_2\ll A_k$ and $B_2\ll B_k$. By \ref{striptobip}, this proves (1).

\bigskip

Since $|K|\le 2$ and $k\ge 3$, we may assume that $K$ is disjoint from $S_1\cup D_1$. Let $v\in A_1\cup C_1\cup D_1$ be bisimplicial.
Since $K$ is anticomplete to $a$, it follows
from \ref{growstrips} that $K$ is anticomplete to $v$, and so $v$ satisfies the theorem. This proves \ref{induct2}.~\bbox

Third, we need:

\begin{thm}\label{induct3}
Let $G$ be an even-hole-free graph, such that \ref{mainthm} holds for all graphs with fewer vertices than $G$.
Let $K$ be a non-dominating clique in $G$ with $|K|\le 2$, and let $a\in V(G)\setminus N[K]$ be splendid. Suppose that
there is no pyramid in $G$ with apex $a$. Then $a$ is bisimplicial.
\end{thm}
\Proof
We begin with:
\\
\\
(1) {\em There do not exist distinct $y_1,y_2,y_3\in N(a)$, pairwise nonadjacent.}
\\
\\
Suppose such $y_1,y_2,y_3$ exist. Now
$G\setminus N[a]$ is connected, and $y_1,y_2,y_3$ all have neighbours in it, since $a$ is splendid. Let $S$ be a minimal connected
induced subgraph of $G\setminus N[a]$ such that $y_1,y_2,y_3$ all have neighbours in $S$. No two of $y_1,y_2,y_3$
have a common neighbour in $V(S)$, since such a vertex would make a 4-hole with $a$ and two of $y_1,y_2,y_3$. Consequently
$|V(S)|\ge 2$, and so there are at least two vertices $x\in V(S)$ such that $S\setminus x$ is connected. Choose two such vertices
$x_1,x_2$ say. From the minimality of $S$, for $i = 1,2$ one of $y_1,y_2,y_3$ has no neighbour in $V(S)\setminus \{x_i\}$,
and so we may assume that for $i = 1,2$, $x_i$ is the unique neighbour of $y_i$ in $V(S)$.
Let $P=p_1\cc p_k$ be an induced path
of $S$ with $p_1=x_1$ and $p_k=x_2$. Now $y_3$ might or might not have neighbours in $V(P)$. Let $Q=q_0\dd q_1\cc q_{\ell}$
be a minimal path in $G[S\cup \{y_3\}]$ where $q_0=y_3$ and $q_{ell}$ has a neighbour in $V(P)$.
(Thus if $y_3$ has a neighbour in $V(P)$ then $\ell=0$.) If $q_{\ell}$ has a unique neighbour $p_i\in V(P)$, there is a theta in $G$
with ends $a,p_i$ and constituent paths
$$a\dd y_1\dd P \dd p_i,$$
$$a\dd y_2\dd P\dd p_i,$$
$$a\dd y_3\dd Q\dd p_i,$$
contrary to \ref{subgraphs}. 

Suppose that $q_{\ell}$ has two nonadjacent neighbours in $V(P)$. Then $\ell=0$ by the minimality of $S$ (because
if $\ell>0$, we could delete from $S$ a vertex of $P$ between the first and last neighbour of $q_{\ell}$ in $P$). 
Let $H$ be the hole induced on $V(P)\cup \{a,y_1,y_2\}$. Then $y_3$ is adjacent to $a$ and not to its neighbours in $H$; and $y_3$
has two other neighbours in $V(H)$, nonadjacent to each other. By \ref{bigvertex}, $G$ admits a full star cutset, contrary to
\ref{nofullstar}.

Thus $q_{\ell}$ has exactly two neighbours in $V(P)$ and they are adjacent, say $p_i, p_{i+1}$.
But then there is a pyramid
with apex $a$, base $\{q_{\ell},p_i,p_{i+1}\}$ and constituent paths
$$a\dd Q\dd q_{\ell},$$
$$a\dd y_1\dd P\dd p_i,$$
$$a\dd y_2\dd P\dd y_{i+1},$$
a contradiction. This proves (1).

\bigskip

We suppose that $a$ is not bisimplicial, and so the graph complement of $G[N(a)]$ is not bipartite, and hence has an induced odd cycle.
It has no induced cycle of length at least six, since $G[N(a)]$ has no 4-hole; and none of length three by (1).
Thus it has an induced cycle of length five, and hence so does $G[N(a)]$. 
Let $v_1\cc v_5\dd v_1$ be a 5-hole of $G$ where $v_1\ll v_5$ are adjacent to $a$.
Choose a connected subgraph $S$ with $V(S)\cap N(a)=\emptyset$, minimal such that at least
four of $v_1\ll v_5$ have a neighbour in $V(S)$.
\\
\\
(2) {\em If $u,v\in \{v_1\ll v_5\}$ are nonadjacent then they have no common neighbour in $V(S)$.}
\\
\\
Because if $s\in V(S)$ is adjacent to both $u,v$ then $s\dd u\dd a\dd v\dd s$ is a 4-hole. This proves (2).
\\
\\
(3) {\em If $P=p_1\cc p_k$ is a path of $S$ such that $p_1v_2$ and $p_kv_4$ are edges, then one of $v_1,v_5$
has a neighbour in $\{p_1\ll p_k\}$.}
\\
\\
Suppose not, and choose $k$ minimum. Thus $v_2\dd p_1\dd p_k\dd v_4$ is an induced path. If $v_3$ is nonadjacent to $p_1\ll p_k$
then there is a theta with ends $v_2,v_4$ and constituent paths
$$v_2\dd v_3\dd v_4,$$
$$v_2\dd v_1\dd v_5\dd v_4,$$
$$v_2\dd p_1\cc p_k\dd v_4,$$
contrary to \ref{subgraphs}. So $v_3$ is adjacent to at least one of $p_1\ll p_k$. Let $v_3$ be adjacent to $n\ge 1$ of $p_1\ll p_k$.
If $n$ is odd then there is an even wheel with centre $v_3$ and hole $a\dd v_2\dd p_1\cc p_k\dd v_4\dd a$;
and if $n$ is even there is an even wheel with centre $v_3$ and hole $v_1\dd v_2\dd p_1\cc p_k\dd v_4\dd v_5\dd v_1$,
in both cases contrary to \ref{subgraphs}. This proves (3).

\bigskip

From (2) it follows that $|V(S)|\ge 2$.
Let $X$ be the set of vertices $x\in V(S)$  such that $S\setminus x$ is
connected. For each $x\in X$, let $T(x)$ be the set of $v\in \{v_1\ll v_5\}$ such that $x$ is the unique neighbour
of $v$ in $V(S)$. The minimality of $S$ implies that $T(x)\ne \emptyset$ for each $x\in X$, and (2) implies that $T(x)$ is a clique.
\\
\\
(4) {\em Exactly four of $v_1\ll v_5$ have a neighbour in $V(S)$.}
\\
\\
Suppose $v_1\ll v_5$ all have a neighbour in $V(S)$. From the minimality of $S$, it follows that $|T(x)|\ge 2$ for each $x\in X$,
and since the sets $T(x)\;(x\in X)$ are pairwise disjoint, it follows that $|X|\le 2$. Since $|V(S)|\ge 2$ and $S$ is connected, it follows
that $S$ is a path of length at least one, and $X$ consists of the ends of $S$. Let $S$ have vertices $s_1\cc s_k$ in order.
Now $T(s_1)$ is a clique,
so we may assume that $T(s_1)=\{v_1,v_2\}$. Since $T(s_1),T(s_k)$ are disjoint, similarly we may assume that $T(s_k)=\{v_4,v_5\}$.
Thus each of $v_1,v_2,v_4,v_5$ has a unique neighbour in $V(S)$, and $v_3$ has at least one such neighbour.
But then there is a 4-hole with centre $v_3$ and hole $s_1\dd S\dd v_3\dd a\dd v_1\dd s_1$, contrary to \ref{subgraphs}. This proves (4).

\bigskip

We may therefore assume that $v_3$ has no neighbour in $V(S)$.
If $x\in X$, then $T(x)\ne \{v_2\}$, since otherwise $S\setminus x$ would contain a path in which $v_1,v_4$ have neighbours and $v_2,v_3$
do not, contrary to (3). Similarly $T(x)\ne \{v_4\}$, and $T(x)\ne \{v_2,v_4\}$ since $T(x)$ is a clique.
Thus $T(x)$ contains one of $v_1,v_5$.
Hence $|X|=2$, and so $S$ is a path $s_1\cc s_k$ say, where $v_1\in T(s_1)$ and $v_5\in T(S_k)$. If both $v_2,v_4$
have a neighbour in $S^*$, there is a theta with ends $v_2,v_4$ and constituent paths
$$v_2\dd v_3\dd v_4,$$
$$v_1\dd v_1\dd v_5\dd v_4,$$
$$v_2\dd G[S^*]\dd v_4,$$
contrary to \ref{subgraphs}. From the symmetry we may therefore assume that $v_2$ has no neighbour in $S^*$. Also by (2), $v_2$
is nonadjacent to $s_k$, so $v_2\in T(s_1)$. Let $v_4$ have $n$ neighbours in $V(S)$. If $n$ is even then
there is an even wheel with centre $v_4$ and hole $a\dd v_2\dd s_1\cc s_k\dd v_5\dd a$, and if $n$ is odd and $n>1$
 then there is an even wheel
with centre $v_4$ and hole $v_1\dd s_1\cc s_k\dd v_5$. Thus $n=1$. Let $s_i$ be the unique neighbour of $v_4$ in $V(S)$.
If $i=k$, there is a prism with bases $\{v_1,v_2,s_1\}$, $\{v_4,v_5,s_k\}$ and constituent paths
$$v_1\dd v_5,$$
$$v_2\dd v_3\dd v_4,$$
$$s_1\cc s_k,$$
contrary to \ref{subgraphs}.
If $i<k$ there is a theta with ends $s_i,v_5$ and constituent paths
$$s_i\cc s_k\dd v_5,$$
$$ s_i\cc s_1\dd v_1\dd v_5,$$
$$s_i\dd v_4\dd v_5,$$
contrary to \ref{subgraphs}. This proves \ref{induct3}.~\bbox

Finally, the fourth part of the proof of \ref{induct}; we will show:

\begin{thm}\label{induct4}
Let $G$ be even-hole-free, and let $K$ be a non-dominating clique in $G$ with $|K|\le 2$. Suppose that
\ref{mainthm} holds for all graphs with fewer vertices than $G$, but there is no bisimplicial vertex of $G$
in $V(G)\setminus N[K]$. Then there is a splendid vertex in $V(G)\setminus N[K]$.
\end{thm}
\Proof
If $K\ne \emptyset$ let $Z$ be the set of all vertices in $V(G)\setminus K$ that are complete to $K$, and if $K=\emptyset$
let $Z=\emptyset$.
Choose $a\in V(G)\setminus N[K]$ with as few neighbours in $Z$ as possible; and subject to that, with degree as small as possible.
We claim that $a$ is splendid. By \ref{nofullstar} we may assume that $G$ admits no full star cutset, and so
for every vertex $v$, the subgraph induced on $V(G)\setminus N[v]$ is connected. In particular, this holds when $v=a$, which is the
first requirement to be splendid.
\\
\\
(1) {\em Every vertex in $N(a)$ has a neighbour in $V(G)\setminus N[a]$.}
\\
\\
Suppose that $v\in N(a)$ has no neighbour in $V(G)\setminus N[a]$. Then every neighbour of $v$ belongs to $N[a]$, and in particular,
$v\notin N[K]$, and every vertex in $Z$ adjacent to $v$ is also adjacent to $a$, and the degree of $v$ is at most that of $a$.
From the choice of $a$, equality holds, and so $a,v$ have the same neighbours (except for $a,v$ themselves).
Let $G'=G\setminus v$. Since $K$ is non-dominating in $G'$, the inductive hypothesis implies that there exists $u\in V(G')\setminus N_{G'}[K]$
that is bisimplicial in $G'$. If $u=a$, then since $v$ is adjacent to every neighbour of $a$, it follows that $a$
is bisimplicial in $G$; so we may assume that $u,v,a$ are all distinct.
If $u,v$ are nonadjacent, then $u$ is bisimplicial in $G$.
If $u,v$ are adjacent, then $u,a$ are adjacent, and since $v,a$ have the same neighbours in $N[u]$,
it follows that $u$ is bisimplicial in $G$. In each case this is impossible. This proves (1).

\bigskip

Suppose there is a short pyramid in $G$ with apex $a$; with base $\{b_1,b_2,b_3\}$ say, and constituent paths $R_1,R_2,R_3$
where $R_i$ has ends $a, b_i$ for $i = 1,2,3$, and $R_3$ has length one. Thus $R_1,R_2$ have length at least three. For $i = 1,2$
let $y_i$ be the neighbour of $a$ in $R_i$. Let $S$ be the set of vertices of $G$ nonadjacent to both $a,b_3$.
\\
\\
(2) {\em If $P=p_1 \cc p_k$ is a path with $p_1\ll p_k\in S$, of minimum length such that $p_1$ has a neighbour in $R_1^*\setminus \{y_1\}$ and
$p_k$ has a neighbour in $V(R_2)$, then $p_1$ has exactly
two adjacent neighbours in $V(R_1)$ and $y_2$ is the unique neighbour of $p_k$ in $V(R_2)$, and these three edges are the only
edges between
$\{p_1\ll p_k\}$ and $V(R_1\cup R_2\cup R_3)$.}
\\
\\
From the minimality of $k$, none of $p_2\ll p_k$ has a neighbour in $R_1^*\setminus \{y_1\}$, but they might be adjacent to $b_1$ or $y_1$.
Also none of $p_1\ll p_{k-1}$ has a neighbour in $V(R_2)$. (Note that possibly $k=1$.)
Suppose that $p_k$ has two nonadjacent neighbours in $V(R_2)$. Then there is a
theta with ends $p_k,a$ and constituent paths
$$p_k\dd R_2\dd a,$$
$$p_k \dd R_2\dd b_2\dd b_3\dd  a$$
$$p_k\dd (P\cup R_1\setminus b_1)\dd a,$$
contrary to \ref{subgraphs}. If $p_k$ has exactly two neighbours $x,y$ in $R_2$
and they are adjacent (and $a,x,y,b_2$ are in this order in $R_2$, say), there is a near-prism with bases
$\{b_1,b_2,b_3\}$ and $\{p_k,x,y\}$, with
constituent paths
$$x\dd R_2\dd a\dd b_3,$$
$$p_k\dd (P\cup R_1\setminus a)\dd b_1,$$
$$y\dd R_2\dd b_2,$$
contrary to \ref{subgraphs}.
Thus $p_k$ has a unique neighbour $u$ say in $V(R_2)$. If $u\ne y_2$, there is a theta with ends $u,a$ and constituent paths
$$u\dd R_2\dd a,$$
$$u\dd R_2\dd b_2\dd b_3\dd a,$$
$$u\dd (P\cup R_1\setminus b_1)\dd a,$$
contrary to \ref{subgraphs}. So $u=y_2$.
Hence $p_k$ is not adjacent to $y_1$, because otherwise there would be a 4-hole $p_k\dd y_2\dd a\dd y_1\dd p_k$.
If $b_1$ is the unique neighbour of $p_k$ in $V(R_1)$, there is a theta with ends $y_1,b_1$ and constituent paths
$$y_1\dd a\dd R_1\dd b_1,$$
$$y_1\dd R_2\dd b_2\dd b_1,$$
$$y_1\dd p_k\dd b_1,$$
contrary to \ref{subgraphs}. So if $p_k$ has a neighbour in $V(R_1)$ then $k=1$. Thus the only edges between $\{p_1\ll p_k\}$
and $V(R_1\cup R_2\cup R_3)$ are the edges between $p_1$ and $V(R_1)$, and the edge $p_ky_2$. If $p_1$ has two nonadjacent neighbours in
$R_1$, say $x,y$ where $a,x,y,b_1$ are in order in $R_1$, then there is a theta with ends $p_1,a$ and constituent paths
$$p_1\dd x\dd R_1\dd a,$$
$$p_1\dd y\dd R_1\dd b_1\dd b_3\dd a,$$
$$p_1\dd P\dd p_k\dd y\dd a,$$
contrary to \ref{subgraphs}. If $p_1$ has a unique neighbour say $v$ in $V(R_1)$, then since $v\ne y_1$ (because by hypothesis,
$p_1$ has a neighbour in $R_1^*\setminus \{y_1\}$), there is a theta with ends $v,a$ and constituent paths
$$v\dd R_1\dd a,$$
$$v\dd R_1\dd b_1\dd b_3\dd a,$$
$$v\dd P\dd y_2\dd a,$$
contrary to \ref{subgraphs}. So $p_1$ has exactly two neighbours in $V(R_1)$ and they are adjacent. Ths proves (2).
\\
\\
(3) {\em There is no path $p_1\ll p_k$ with $p_1\ll p_k\in S$, such that $p_1$ has a neighbour in $R_1^*\setminus \{y_1\}$
and $p_2$ has a neighbour in $R_2^*\setminus \{y_2\}$.}
\\
\\
Suppose $P=p_1\ll p_k$ is such a path, chosen with $k$ minimum.
Note that $y_1,y_2,b_1,b_2$ may have neighbours in the interior of $P$, but from the minimality of $k$,
$p_1\ll p_{k-1}$ have no neighbours in $R_2^*\setminus \{y_2\}$, and $p_2\ll p_k$ have no neighbours in $R_1^*\setminus \{y_1\}$.
Choose $i\in \{1\ll k\}$ minimum such that $p_i$ has a neighbour in $V(R_2)$. From (2) applied to the path
$p_1\cc p_i$, it follows that $p_1$ has exactly two neighbours in $V(R_1)$, say $x_1,y_1$, and they are adjacent, and $y_2$
is the unique neighbour of $p_i$ in $V(R_2)$, and these three edges are the only edges between
$\{p_1\ll p_i\}$ and $V(R_1\cup R_2\cup R_3)$. In particular $i<k$. Choose $j\in \{1\ll k\}$ maximum such that $p_j$ has a neighbour in
$V(R_1)$; then similarly $p_k$ has exactly two neighbours in $V(R_2)$, say $x_2,y_2$, and they are adjacent, and $y_1$
is the unique neighbour of $p_j$ in $V(R_1)$, and these three edges are the only edges between
$\{p_j\ll p_k\}$ and $V(R_1\cup R_2\cup R_3)$. Thus $j>i$, and since $1\le i<j\le k$ it follows that $k\ge 2$.
Let $Q$ be the path $p_i\dd p_{i+1}\cc p_j$.
Thus the only edges between $\{p_1\ll p_k\}$ and $V(R_1\cup R_2\cup R_3)$
are edges between $p_1$ and $V(R_1)$, edges between $p_k$ and $V(R_2)$, the edges $p_iy_2, p_jy_1$, and edges between
$Q^*$ and $\{y_1,y_2,b_1,b_2\}$. If $b_1$ has a neighbour in $Q^*$, there is a theta
with ends $b_1,y_1$ and constituent paths
$$b_1\dd R_1\dd y_1,$$
$$b_1\dd b_3\dd a\dd y_1,$$
$$b_1\dd Q\dd y_1,$$
contrary to \ref{subgraphs}. So $b_1$ has no neighbour in $\{p_2\ll p_k\}$, and similarly $b_2$ has no neighbour in $\{p_1\ll p_{k-1}\}$.
If $y_1,y_2$ both have neighbours in $P^*$, there is a theta with ends $y_1,y_2$ and constituent paths
$$y_1\dd G[P^*]\dd y_2,$$
$$y_1\dd a\dd y_2,$$
$$y_1\dd R_1\dd b_1\dd b_2\dd R_2\dd y_2,$$
contrary to \ref{subgraphs}. Thus we may assume that $y_2$ has no neighbour in $P^*$, and in particular $i=1$. Consequently
$p_1,y_1$ are nonadjacent, since $p_1\dd y_1\dd a\dd y_2\dd p_1$ is not a 4-hole. Then there is a theta with ends
$p_1,y_1$ and constituent paths
$$p_1\dd R_1\dd y_1,$$
$$p_1\dd R_1\dd b_1\dd b_3\dd a\dd y_1,$$
$$p_1\dd P\dd y_1,$$
contrary to \ref{subgraphs}. This proves (3).

\bigskip

For $i = 1,2$, let $S_i$ be the component of $G[S]$ that contains $R_i\setminus \{a,y_i,b_i\}$.
So $S_1,S_2$ are nonempty since $R_1,R_2$ have length at least three; and $S_1,S_2$ are distinct by (3).
For $i = 1,2$, let $B_i$ be the set of vertices
adjacent to $b_3$ and not to $a$, with a neighbour in $S_i$. So $b_i\in  B_i$ for $i = 1,2$. If there exists
$v\in B_1\cap B_2$, there is a theta with ends $v,a$ and constituent paths
$$v\dd b_3\dd a,$$
$$v\dd S_1\dd y_1\dd a,$$
$$v\dd S_2\dd y_2\dd a,$$
contrary to \ref{subgraphs}. So $B_1\cap B_2=\emptyset$.

The only vertices of $G$ not in $V(S_1)$ but with a neighbour in $V(S_1)$ belong to $B_1\cup N[a]$.
From the inductive hypothesis, applied to the graph $G'=G[S_1\cup B_1\cup N[a]\cup \{b_3\}]$, since the edge $ab_3$ is non-dominating in $G'$, it follows that some vertex in $S_1$ is bisimplicial in $G'$ and hence in $G$. Since there is no bisimplicial vertex
of $G$ in $V(G)\setminus N[K]$, it follows that $N[K]\cap S_1\ne \emptyset$, and similarly $N[K]\cap S_2\ne \emptyset$.
But $K\cap N[a]=\emptyset$ from the choice of $a$; and so $K\cap (V(S_i)\cup B_i)\ne \emptyset$ for $i = 1,2$.
Since the sets $V(S_1), B_1,B_2, V(S_2)$ are pairwise disjoint, and there are no edges between $B_2\cup V(S_2)$ and $V(S_1)$,
it follows that $K\cap S_1=\emptyset$, and so $K\cap B_1\ne \emptyset$; and similarly $K\cap B_2\ne \emptyset$.  In particular $|K|=2$.
Let $K\cap B_i=b_i'$ for $i = 1,2$.

We recall that $Z$ is the set of all vertices adjacent to both $b_1',b_2'$, and so $b_3\in Z$. Now $a,b_3$ are adjacent. But
$y_1\notin N[K]$ (because if $y_1, b_i'$ are adjacent then there is a 4-hole $y_1\dd b_i'\dd b_3\dd a\dd y_1$), and $y_1,b_3$
are nonadjacent. From the choice of $a$, $y_1$ has at least as many neighbours in $Z$ as does $a$; and since $b_3$
is adjacent to $a$ and not to $y_1$, there exists $z\in Z$ adjacent to $y_1$ and not to $b_3$. Since $z\dd y_1\dd a\dd b_3\dd z$
is not a 4-hole, $z,b_3$ are nonadjacent. Since $b_2'\in B_2$ and hence $b_2'\notin B_1$, and $b_2'$ is adjacent to $z$,
it follows that $z\notin V(S_1)$.
But then there is a theta with ends $b_1', y_1$ and constituent paths
$$b_1'\dd z\dd y_1,$$
$$b_1'\dd b_3\dd a\dd y_1,$$
$$b_1'\dd S_1\dd y_1,$$
contrary to \ref{subgraphs}. This proves \ref{induct4}.~\bbox

From \ref{induct1}, \ref{induct2}, \ref{induct3} and \ref{induct4}, this completes the proof of \ref{induct}, and hence of \ref{mainthm}.

\section*{Acknowledgement}
We would like to thank Rong Wu for bringing the error in \cite{bisimplicial} to our attention, which gave rise to this paper.

\end{document}